\documentclass[final,leqno]{siamltex}

\usepackage{amssymb}
\usepackage{amsmath}
\usepackage{url}
\usepackage{enumerate}
\usepackage{graphicx}
\usepackage{enumerate}
\usepackage{subfig}

\newtheorem{remark}[theorem]{Remark}

\def\f{{\infty}}

\def\tr{\textrm{tr}}
\newcommand{\Prob}{\mathbb{P}}
\newcommand{\E}{\mathbb{E}}
\def\qed
{\hfill\vbox{\hrule width 0.5em\nointerlineskip\hbox to
0.5em{\vrule height 0.5em \hfill\vrule height
0.5em}\nointerlineskip\hrule width 0.5em}}

\title{$\epsilon$-Nash Mean Field Game Theory for Nonlinear Stochastic Dynamical Systems with Major and Minor Agents\thanks{A brief version of this paper was presented at the 51th IEEE CDC Conference, Maui, HI, Dec. 2012 \cite{NourianCainesCDC2012}.}}

\author{Mojtaba Nourian\thanks{Department of Electrical and Electronic Engineering, The University of Melbourne, VIC 3010, Australia (email: {\tt mojtaba.nourian@unimelb.edu.au}). M. Nourian's work was performed at the Centre for Intelligent Machines (CIM) and Department of Electrical \& Computer Engineering, McGill University, Montreal, QC H3A 2A7, and GERAD, Montreal, Canada.} \and Peter E.
Caines\thanks{CIM and Department of Electrical \& Computer Engineering, McGill University,
Montreal, QC H3A 2A7, and GERAD, Montreal, Canada (email: {\tt peterc@cim.mcgill.ca}).}}
\begin{document}

\maketitle

\begin{abstract} This paper studies a large population dynamic game involving nonlinear stochastic
dynamical systems with agents of the following mixed types: (i) a major agent, and (ii) a population of $N$ minor agents where $N$ is very large. The major and minor (MM) agents are coupled via both: (i) their
individual nonlinear stochastic dynamics, and (ii) their individual finite time horizon nonlinear
cost functions. This problem is approached by the so-called $\epsilon$-Nash Mean Field Game ($\epsilon$-NMFG)
theory. A distinct feature of the mixed agent MFG problem is that even asymptotically (as the
population size $N$ approaches infinity) the noise process of the major agent causes random
fluctuation of the mean field behaviour of the minor agents. To deal with this, the overall
asymptotic ($N \rightarrow \infty$) mean field game problem is decomposed into: (i) two non-standard
stochastic optimal control problems with random coefficient processes which yield forward adapted
stochastic best response control processes determined from the solution of (backward in time)
stochastic Hamilton-Jacobi-Bellman (SHJB) equations, and (ii) two stochastic coefficient
McKean-Vlasov (SMV) equations which characterize the state of the major agent
and the measure determining the mean field behaviour of the minor agents. This yields to a Stochastic Mean Field Game (SMFG) system which is in contrast to the deterministic mean field game system of the standard MFG problems with only minor agents. Existence and uniqueness of the solution to the SMFG system (SHJB and SMV equations) is established by a fixed point argument in the Wasserstein space of random probability measures. In the case
that minor agents are coupled to the major agent only through their cost functions, the $\epsilon_N$-Nash equilibrium property of the SMFG best responses is shown for a finite $N$ population system where $\epsilon_N=O(1/\sqrt N)$.
\end{abstract}



\pagestyle{myheadings}
\thispagestyle{plain}
\markboth{M. Nourian and P. E. Caines}{Mean Field Game Theory Involving Major and Minor Agents}

\section{Introduction} An important class of games is that of dynamic games with a very large number
of minor agents in which each agent interacts with the average (or so-called mean field) effect of
other agents via couplings in their individual dynamics and individual cost functions. A minor agent is an agent which, asymptotically as the population size goes to infinity, has
a negligible influence on the overall system while the overall population's effect on it is
significant. Stochastic dynamic games with mean field couplings arise in fields such as wireless
power control \cite{HuangCaines_03c}, consensus dynamics \cite{Mnourian_IEEE11_Con}, flocking \cite{nourian2011mean}, charging control of plug-in electric vehicles \cite{ma2012optimal}, synchronization of coupled nonlinear oscillators \cite{yin2010synchronization}, crowd dynamics \cite{dogbe2010modeling} and economics \cite{Weintraub_08,gueant2011mean}.

For large population stochastic dynamic games with mean field couplings and no major agent, the $\epsilon$-Nash Mean Field
Game ($\epsilon$-NMFG) (or Nash Certainty Equivalence (NCE)) theory was originally developed as a decentralized
methodology in a series of papers by Huang together with Caines and Malham\'e, see \cite{HuangCaines_03c,HuangCaines_TAC07} for the $\epsilon$-NMFG linear-quadratic-Gaussian (LQG) framework, and \cite{HuangMalhameCaines_06,huang2007invariance,Peter2009bode} for a general formulation of nonlinear McKean-Vlasov type $\epsilon$-NMFG problems. For this class of game problems a
closely related approach has been independently developed by Lasry and Lions
\cite{LasryLionsStationair2006,LasryLionsHorizonfini2006,Lasry07,gueant2011mean} where the term Mean
Field Games (MFG) was initially used. For models of many firm industry dynamics, Weintraub et. al.
proposed the notion of oblivious equilibrium by use of mean field approximations
\cite{Weintraub_05}. The $\epsilon$-NMFG framework for LQG systems is extended to systems of agents with ergodic (long time
average) costs in \cite{LiZhang_TAC08}, while Kolokoltsov et. al. extend the $\epsilon$-NMFG theory to
general nonlinear Markov systems \cite{kolokoltsov2011mean}. The extension of the $\epsilon$-NMFG framework so as to model the collective system dynamics which include large population of leaders and followers, and an unknown (to the followers) reference trajectory for the leaders is studied in \cite{Mnourian_LF11}. The reader is referred to the survey paper \cite{buckdahn2011some} for some of the research on MFG theory up to 2011.

The central idea of the $\epsilon$-NMFG theory is to specify a certain equilibrium relationship between
the individual strategies and the mass effect (i.e., the overall effect of the population on a given
agent) as the population size goes to infinity \cite{HuangCaines_TAC07}. Specifically, in the
equilibrium: (i) the individual strategy of each agent is a best response to the infinite
population mass effect in the sense of a so-called $\epsilon$-Nash equilibrium, and (ii) the set of
strategies collectively replicates the mass effect, this being a dynamical game theoretic fixed
point property. The defining property of the $\epsilon$-NMFG equilibrium with individual strategies
$\{u_i^o: 1 \leq i \leq N\}$ requires that for any given $\epsilon > 0$, there exists $N(\epsilon)$
such that for any population size $N(\epsilon) \leq N$, when any agent $j$, $1 \leq j \leq N$,
distinct from $i$ employs $u_j^o$, then agent $i$ can benefit at most $\epsilon$ by unilaterally
deviating from his strategy $u_i^o$, and this holds for all $1 \leq i \leq N$. The estimates in \cite{HuangCaines_03c,HuangMalhameCaines_06,HuangCaines_TAC07} show $\epsilon=O(1/\sqrt N)$ while distinct estimates are obtained in the framework of \cite{kolokoltsov2011mean}.  

A stochastic maximum principle for control problems of mean field type is studied in
\cite{andersson2011maximum} where the state process is governed by a stochastic differential
equation (SDE) in which the coefficients depend on the law of the SDE. The reader is referred to
\cite{buckdahn2009mean1,buckdahn2009mean2} for the analysis of forward--backward stochastic
differential equations (FBSDEs) of mean field type and their related partial differential equations.

Recently, Huang \cite{HuangSIAM2010} introduced a large population LQG dynamic game model with mean field couplings which involves not only a large number of multi-class minor agents but also a major
agent with a significant influence on minor agents (see
\cite{hart1973values,haimanko2000nonsymmetric,neyman2002values} for static cooperative games of
agents with different influences or so-called mixed agents). Since all minor agents respond to
the same major agent, the mean field behaviour of minor agents in each class is directly impacted by the major agent and hence is a random process \cite{HuangSIAM2010}. This is in contrast to the situation in the standard MFG models with only minor agents. A state-space augmentation approach for the
approximation of the mean field behaviour of the minor agents is taken in order to Markovianize the problem and hence to obtain $\epsilon$-NMFG equilibrium strategies \cite{HuangSIAM2010}. An extension of the model in \cite{HuangSIAM2010} to the systems of agents with Markov jump parameters in their dynamics and random parameters in their cost functions is studied in \cite{wang2011distributed} in a discrete-time setting. See also \cite{2011KC_CDC} for the extension of the model in \cite{HuangSIAM2010} to the case of systems with egoistic and altruistic agents.

The model of \cite{HuangSIAM2010} with finite classes of minor agents is extended in
\cite{NguyenHuang_CDC11} to the case of minor agents parameterized by an infinite set of dynamical
parameters where the state augmentation trick cannot be applied to obtain a finite dimensional
Markov model. Due to the LQ structure of the problem an appropriate representation for the mean
field behaviour of the minor agents as a random process is assumed which depends linearly on the
random initial state and Brownian motion of the major agent. Appropriate
approximation of the model by LQG control problems with random parameters in the dynamics and costs
yields non-Markovian forward adapted $\epsilon$-NMFG strategies resulting from backward stochastic
differential equations (BSDEs) obtained by a stochastic maximum principle \cite{NguyenHuang_CDC11}.

In this paper we extend the LQG model for major and minor (MM) agents \cite{HuangSIAM2010} to the
case of a nonlinear stochastic dynamic games formulation of controlled McKean-Vlasov (MV) type
\cite{HuangMalhameCaines_06}. Specifically, we consider a large population dynamic game involving
nonlinear stochastic dynamical systems with agents of the following mixed types: (i) a major agent,
and (ii) a population of $N$ minor agents where $N$ is very large. The MM agents are coupled via both: (i) their
individual nonlinear stochastic dynamics, and (ii) their individual finite time horizon nonlinear
cost functions. 

Applications of the major and minor formulation may be found in charging control of plug-in electric vehicles \cite{Ma_PEV_CDC11,ma2012optimal}, economic and social opinion models with an influential leader (e.g., \cite{dring2009boltzmann}), and power markets involving large consumers and large utilities together with many domestic consumers represented by smart meter agents and possibly large numbers of renewable energy based generators \cite{2012KMC_CDC}.

A distinctive feature of the mixed agent MFG problem is that even asymptotically (as the population size $N$ approaches infinity) the noise process of the major agent causes random fluctuation of the mean field behaviour of the minor agents \cite{HuangSIAM2010,NguyenHuang_CDC11}. 

The main contributions of the paper are as follows:
\begin{itemize}
\item The overall asymptotic ($N \rightarrow \infty$) mean field game problem is decomposed into: (i) two non-standard Stochastic Optimal Control Problems (SOCPs) with random coefficient processes which yield forward adapted stochastic best response control processes determined from the solution of (backward in time) stochastic Hamilton-Jacobi-Bellman (SHJB) equations, and (ii) two stochastic coefficient McKean-Vlasov (SMV) equations which characterize the state of the major agent and the measure determining the mean field behaviour of the minor agents. This yields to a Stochastic Mean Field Game (SMFG) system which is in contrast to the deterministic mean field game system of the standard MFG problems with only minor agents.   

\item Existence and uniqueness of the solution to the SMFG system (SHJB and SMV equations) is established by a fixed point argument in the Wasserstein space of random probability measures. 

\item In the case that minor agents are coupled to the major agent only through their cost functions, the $\epsilon_N$-Nash equilibrium property of the SMFG best responses is shown for a finite $N$ population system where $\epsilon_N=O(1/\sqrt N)$. 

\item As a particular but important case, the results of Nguyen and Huang \cite{NguyenHuang_CDC11} for major and minor agent MFG LQG systems with homogeneous population are retrieved in Appendix G in \cite{MNourian_SiamAppendix}. 

\item Finally, the results of this paper are illustrated with a major and minor agent version of a game model of the synchronization of coupled nonlinear oscillators \cite{yin2010synchronization} (see Appendix H in \cite{MNourian_SiamAppendix}). 
\end{itemize}

It is to be emphasized that the non-standard nature of the SOCPs in (i), which consists of the coupling through the SMV equations in (ii), arises from a distinct feature of the problem formulation. The source of this non-standard nature is the game structure whereby the minor agents are (through the Principle of Optimality) optimizing with respect to the future stochastic evolution of the major agent's state which is partly a result of that agent's future best response control actions. This feature vanishes in the non-game theoretic setting of one controller with one cost function with respect to the trajectories of all the system components (the classical SOCPs), moreover it also vanishes in the infinite population limit of the standard $\epsilon$-NMFG models with no major agent. This is true for both completely and partially observed SOCPs. The nonstandard feature of the SOCPs here give rise to the analysis of systems with (non necessarily Markovian) stochastic parameters. Here, as in
\cite{NguyenHuang_CDC11,yong2002leader}, the theory of BSDEs (see in particular
\cite{bismut1976linear,pardoux1990adapted,peng1992stochastic,peng1993backward}) is used in the
resulting stochastic dynamic game theory. More specifically, we utilize techniques
from \cite{peng1992stochastic} which applies the Principle of Optimality to a stochastic nonlinear control problem with random coefficients; this leads to a formulation of a SHJB equation by use of (i) a semi-martingale representation for the corresponding stochastic value function, and (ii) the It\^{o}-Kunita formula. An application of Peng results to portfolio-consumption optimization under habit formation in complete markets is studied in \cite{egglezos2007aspects}.

The organization of the paper is as follows. Section \ref{S2:mm:Prob} is dedicated to the problem formulation. A McKean-Vlasov approximation for major and minor agent system is studied in Section \ref{S4:mm:MFCT}. Section \ref{S3:mm:SHJB} presents a preliminary nonlinear SOCP with random parameters. The SMFG system of equations of the MM agents is given in Section \ref{S5:mm:SMF}, and the existence and uniqueness of its solution is established in Section \ref{S6:mm:Exis}. The $\epsilon$-Nash equilibrium property of the resulting SMFG control laws is studied in Section \ref{S8:mm:Nash}. Finally, Section \ref{S9:mm:Conc} concludes the paper.

\subsection{Notation and Terminology} The following notation will be used throughout the paper. Let $\mathbb{R}^n$ denote the $n$-dimensional real Euclidean space with the standard Euclidean norm $|\cdot|$ and the standard Euclidean inner product $\big<\cdot,\cdot\big>$. The transpose of a vector (or matrix)
$x$ is denoted by $x^T$. $\tr(A)$ denotes the trace of a square matrix $A$. Let $\mathbb{R}^{n \times m}$ be the Hilbert space consisting of all $(n \times m)$-matrices with the inner product $<A,B>:=\tr(AB\textsuperscript{T})$ and the norm $|A|:=<A,A>^{1/2}$. The set of non-negative real numbers is denoted by $\mathbb{R}_{+}$. $T \in [0,\infty)$ is reserved to denote the terminal time. The integer $N$ is reserved to designate the population size of the minor agents. The superscript $N$ for a process (such as state, control or cost function) is used to indicate the dependence on the population size $N$. We use the subscript 0 for the major agent $\mathcal A_0$ and an integer valued subscript for an individual minor agent $\{\mathcal A_i: 1 \leq i \leq N\}$. At time $t \geq 0$, (i) the states of agents $\mathcal A_0$ and $\mathcal A_i$ are respectively denoted by $z_0^N(t)$ and $z_i^N(t)$, $1 \leq i \leq N$, and (ii) for the system configuration of minor agents $(z_1^N(t), \cdots,z_N^N(t))$ the empirical distribution $\delta_t^N$ is defined as the normalized sum of Dirac's masses, i.e., $\delta_t^N:=(1/N) \sum_{i=1}^N \delta_{z_i^N(t)}$ where $\delta_{(\cdot)}$ is the Dirac measure. $C(S)$ is the set of continuous functions and $C^k(S)$ the set of $k$-times continuously differentiable functions on $S$. The symbol $\partial_t$ denotes the partial derivative with respect to variables $t$. We denote $D_x$ and $D^2_{xx}$ as the gradient and Hessian operators with respect to the variable $x$. These are respectively denoted by $\partial_x$ and $\partial^2_{xx}$ when applied to a function defined on a one-dimensional domain. Let $(\Omega, \mathcal F,\{\mathcal F_t\}_{t\geq 0},\Prob)$ be a complete filtered probability space. $\E$ denotes the expectation. The conditional expectation with respect to the $\sigma$-field $\mathcal V$ is denoted by $\E_{\mathcal V}$. For an Euclidean space $H$ we denote by $L^2_{\mathcal G}([0,T];H)$ the space of all $\{\mathcal G_t\}_{t \geq 0}$-adapted $H$-valued processes $f(t,\omega)$ such that
$\E\int_0^T|f(t,\omega)|^2dt<\infty$. We use the notation $(\E_{\omega} h)(z) := \int h(z,\omega) \Prob_\omega(d\omega)$ for any function $h(z,\omega)$ and sample point $\omega \in \Omega$. Finally, note that we may not display the dependence of random variables or stochastic processes on the sample point $\omega \in \Omega$.

\section{Problem Formulation} \label{S2:mm:Prob} We consider a dynamic game involving: (i) a major
agent $\mathcal A_0$, and (ii) a population of $N$ minor agents $\{\mathcal A_i: 1 \leq i \leq
N\}$ where $N$ is very large. We assume homogenous minor agents although the modelling may be generalized to the case of
multi-class heterogeneous minor agents \cite{HuangMalhameCaines_06,HuangSIAM2010} (see \cite{Mnourian_MTNS12}).

The dynamics of the agents are given by the following controlled It\^{o} stochastic differential equations on $(\Omega, \mathcal F,\{\mathcal F_t\}_{t
\geq 0},\Prob)$:
\begin{align}
\label{CNP:Major:GenDyn}
& dz_0^N(t) =\frac{1}{N} \sum_{j=1}^N f_0[t,z_0^N(t),u_0^N(t),z_j^N(t)]dt \\
& \hspace{1cm} + \frac{1}{N} \sum_{j=1}^N \sigma_0[t,z_0^N(t),z_j^N(t)] dw_0(t),\quad z_0^N(0)=
z_0(0), \quad 0 \leq t \leq T, \notag \\
\label{CNP:Minor:GenDyn}
& dz_i^N(t) = \frac{1}{N} \sum_{j=1}^N f[t,z_i^N(t),u_i^N(t),z_0^N(t),z_j^N(t)] dt \\
& \hspace{1cm} +\frac{1}{N} \sum_{j=1}^N \sigma [t,z_i^N(t),z_0^N(t),z_j^N(t)] dw_i(t),\quad z_i^N(0)=
z_i(0), \quad 1 \leq i \leq N, \notag
\end{align}
with terminal time $T \in (0,\infty)$ where (i) $z_0^N: [0,T] \rightarrow \mathbb{R}^n$ is the state
of the major agent $\mathcal A_0$ and $z_{i}^N: [0,T] \rightarrow \mathbb{R}^n$ is the state of the
minor agent $\mathcal A_i$; (ii) $u_0^N: [0,T] \rightarrow U_0$ and $u_i^N:
[0,T] \rightarrow U$ are respectively the control inputs of $\mathcal A_0$ and
$\mathcal A_i$; (iii) $f_{0}:[0,T]\times\mathbb{R}^n\times U_0\times\mathbb{R}^n \rightarrow
\mathbb{R}^n$, $\sigma_0:[0,T]\times\mathbb{R}^n\times\mathbb{R}^n \rightarrow \mathbb{R}^{n\times
m}$, $f:[0,T]\times\mathbb{R}^n\times U\times\mathbb{R}^n\times\mathbb{R}^n\rightarrow\mathbb{R}^n$ and $\sigma:[0,T]\times\mathbb{R}^n\times\mathbb{R}^n\times\mathbb{R}^n\rightarrow\mathbb{R}^{n\times m}$; (iv) the set of initial states is given by $\{z_j^N(0)= z_j(0):0
\leq j \leq N\}$, and (v) the sequence $\{(w_{j}(t))_{t\geq0}: 0 \leq  j \leq N\}$ denotes $N+1$
mutually independent standard Brownian motions in $\mathbb{R}^m$. We denote the filtration
$\mathcal F_t$ as the $\sigma$-field generated by the initial states and the Brownian motions up to
time $t$, i.e., $\mathcal F_t:= \sigma\{z_j(0),w_j(s): 0 \leq j \leq N, 0 \leq s \leq t\}$. We also set
$\mathcal F_t^{w_0} = \sigma\{z_0(0),w_0(s): 0 \leq s \leq t\}$. These filtrations are augmented by all the $\Prob$-null sets in $\mathcal F$.

For $0 \leq j \leq N$ denote $u_{-j}^N:=\{u_0^N, \cdots,u_{j-1}^N,u_{j+1}^N, \cdots, u_N^N\}$. The
objective of each agent is to minimize its finite time horizon nonlinear cost function given by
\begin{align}
\label{CNP:Major:GenCost} & J_{0}^{N}(u_0^N;u_{-0}^N) := \E \int_0^T \Big(\frac{1}{N}\sum_{j=1}^N
L_0[t,z_0^N(t),u_0^N(t),z_j^N(t)]\Big)dt, \\
\label{CNP:Minor:GenCost} & J_{i}^{N}(u_i^N;u_{-i}^N) := \E \int_0^T \Big(\frac{1}{N}\sum_{j=1}^N
L[t,z_i^N(t),u_i^N(t),z_0^N(t),z_j^N(t)]\Big)dt,
\end{align}
for $1 \leq i \leq N$, where $L_0:[0,T]\times\mathbb{R}^n\times
U_0\times\mathbb{R}^n\rightarrow\mathbb{R}_{+}$ and $L(z_i,u_i,z_0,x):[0,T]\times\mathbb{R}^n\times
U\times\mathbb{R}^n\times\mathbb{R}^n\rightarrow\mathbb{R}_{+}$ are the nonlinear cost-coupling
functions of the major and minor agents. For $0 \leq j \leq N$, we indicate the dependence of $J_j$
on $u_j^N$, $u_{-j}^N$ and the population size $N$ by $J_j^N(u_j^N;u_{-j}^N)$.

We note that in the modelling (\ref{CNP:Major:GenDyn})-(\ref{CNP:Minor:GenCost}) the major agent $\mathcal A_0$ has a significant influence on minor agents while each minor agent has an asymptotically negligible impact on other agents in a large $N$ population system. The major and minor agents are coupled via both: (i) their individual nonlinear stochastic dynamics (\ref{CNP:Major:GenDyn})-(\ref{CNP:Minor:GenDyn}), and (ii) their individual finite time horizon nonlinear cost functions (\ref{CNP:Major:GenCost})-(\ref{CNP:Minor:GenCost}). 

We note that the coupling terms may be written as functionals of the empirical distribution $\delta_{(\cdot)}^N$ by the formula $\int_{\mathbb{R}^n} \phi(x) \delta_{t}^N(dx) = (1/N) \sum_{i=1}^N \phi(x_i(t))$ for a bounded continuous function $\phi$ in $\mathbb{R}^n$.

\begin{remark} \label{RemarkforLQGcase} Under suitable conditions, the results of this paper may be adapted to deal with cost-couplings of the form: 
\begin{align*}
&L_0[t,z_0^N(t),u_0^N(t),z_j^N(t),\frac{1}{N}\sum_{j=1}^N z_j^N(t)],~  L[t,z_i^N(t),u_i^N(t),z_0^N(t),z_j^N(t),\frac{1}{N}\sum_{j=1}^Nz_j^N(t)],
\end{align*}
in (\ref{CNP:Major:GenCost})-(\ref{CNP:Minor:GenCost}).
\end{remark}

\subsection{Assumptions} \label{Sub:ass} Let the empirical distribution of $N$ minor agents' initial states be
defined by $F_N (x) = (1/N) \sum_{i=1}^N 1_{\{\E z_i(0)<x\}}$, where $1_{\{\E z_i(0)<x\}}=1$ if
$\E z_i(0)<x$, and $1_{\{\E z_i(0)<x\}}=0$ otherwise. We enunciate the following assumptions:

({\bf A1}) The initial states $\{z_j(0): 0 \leq j \leq N\}$ are
$\mathcal F_0$-adapted random variables mutually independent and independent of all Brownian motions $\{(w_{j}(t))_{t\geq0}:0 \leq  j \leq N\}$, and there exists a constant $k$ independent of $N$ such that $\sup_{0 \leq j \leq N} \E |z_j(0)|^2 \leq k <\infty$.

({\bf A2}) $\{F_N: N \geq 1\}$ converges to a probability distribution $F$ weakly, i.e., for any bounded and continuous function $\phi$ on $\mathbb{R}^n$ we have $\lim_{N \rightarrow \infty} \int_{\mathbb{R}^n} \phi(x)dF_N(x) = \int_{\mathbb{R}^n} \phi(x)dF(x).$

({\bf A3}) $U_0$ and $U$ are compact metric spaces.

({\bf A4}) The functions $f_0[t,x,u,y]$, $\sigma_0[t,x,y]$, $f[t,x,u,y,z]$ and $\sigma[t,x,y,z]$ are continuous and bounded with respect to all their parameters, and Lipschitz continuous in $(x,y,z)$. In addition, their first order derivatives (w.r.t. $x$) are all uniformly continuous and bounded with respect to all their parameters, and Lipschitz continuous in $(y,z)$.

({\bf A5}) $f_0[t,x,u,y]$ and $f[t,x,u,y,z]$ are Lipschitz continuous in $u$.

({\bf A6}) $L_0[t,x,u,y]$ and $L[t,x,u,y,z]$ are continuous and bounded with respect to all their parameters, and Lipschitz continuous in $(x,y,z)$. In addition, their first order derivatives (w.r.t. $x$) are all uniformly continuous and bounded with respect to all their parameters, and Lipschitz continuous in $(y,z)$.

({\bf A7}) ({\it Non-degeneracy Assumption}) There exists a positive constant $\alpha$ such that
\begin{align*}
& \sigma_0[t,x,y]\sigma_0^T[t,x,y] \geq \alpha I, \quad
\sigma[t,x,y,z]\sigma^T(t,x,y,z) \geq \alpha I, \qquad \forall~(t,x,y,z),
\end{align*}
where $\sigma_0$ and $\sigma$ are given in (\ref{CNP:Major:GenDyn}) and (\ref{CNP:Minor:GenDyn}). 

\section{McKean-Vlasov Approximation for Mean Field Game Analysis} \label{S4:mm:MFCT} Motivated by the analysis in Section I.1 of \cite{sznitman1991topics} and in Section 8.1 of \cite{HuangMalhameCaines_06}, we take a probabilistic approach to establish the following asymptotic properties: (i) The influence of any minor agent $\mathcal A_i$ on any other minor agent $\mathcal A_j$ is asymptotically negligible as the population size $N$ goes to infinity, and (ii) In the limit, the effect of the mass of agents on a given minor agent $\mathcal A_i$ is that of the behaviour of a mass of predictable generic agents. This is in the form of a single mean field function in the LQG case \cite{HuangCaines_03c,HuangCaines_TAC07} or a predictable state probability distribution in the nonlinear case \cite{HuangMalhameCaines_06,huang2007invariance,Lasry07}.

Let $\varphi_0(\omega,t,x):\Omega \times [0,T] \times \mathbb{R} \rightarrow U_0$ and
$\varphi(\omega,t,x):\Omega \times [0,T] \times \mathbb{R} \rightarrow U$ be two arbitrary $\mathcal
F_t^{w_0}$-measurable stochastic processes for which we introduce the following assumption:

({\bf H4}) $\varphi_0(\omega,t,x)$ and $\varphi(\omega,t,x)$ are Lipschitz continuous in $x$,
and $\varphi_0(\omega,t,0)\in L_{\mathcal F_t^{w_0}}^2([0,T];U_0)$ and $\varphi(\omega,t,0) \in L_{\mathcal F_t^{w_0}}^2([0,T];U)$.

We assume that $\varphi_0(t,x):=\varphi_0(\omega,t,x)$ and
$\varphi(t,x):=\varphi(\omega,t,x)$ are respectively used by the major and minor agents as their
control laws in (\ref{CNP:Major:GenDyn}) and (\ref{CNP:Minor:GenDyn}) (i.e.,
$u_0=\varphi_0$ and $u_i=\varphi$ for $1 \leq i \leq N$). Then we have the following
closed-loop equations with random coefficients:
\begin{align*}
& d\hat z_0^N(t) =\frac{1}{N} \sum_{j=1}^N f_0[t,\hat z_0^N(t),\varphi_0(t,\hat
z_0^N(t)),\hat z_j^N(t)]dt \\
& \qquad + \frac{1}{N} \sum_{j=1}^N \sigma_0[t,\hat z_0^N(t),\hat z_j^N(t)] dw_0(t), \quad
\hat z_0^N(0)=z_0(0), \quad  0 \leq t \leq T,\\
& d\hat z_i^N(t) = \frac{1}{N} \sum_{j=1}^N f[t,\hat z_i^N(t),\varphi(t,\hat
z_i^N(t)),\hat z_0^N(t),\hat z_j^N(t)] dt \\
& \qquad +\frac{1}{N} \sum_{j=1}^N \sigma [t,\hat z_i^N(t),\hat z_0^N(t),\hat z_j^N(t)]
dw_i(t), \quad \hat z_i^N(0)=z_i(0), \quad 1 \leq i \leq N.
\end{align*}
Under ({\bf A4})-({\bf A5}) and ({\bf H4}) there exists a unique solution $\big(z_0^N(\cdot), \cdots, z_N^N(\cdot)\big)$ to the above system (see Theorem 6.16, Chapter 1 of \cite{yong1999stochastic}, page 49).

We now introduce the McKean-Vlasov (MV) system
\begin{align*}
& d\bar z_0(t) = f_0[t,\bar z_0(t),\varphi_0(t,\bar z_0(t)),\mu_t] dt +
\sigma_0[t,\bar z_0(t),\mu_t]dw_0(t), \quad 0
\leq t \leq T, \\
& d\bar z(t)= f[t,\bar z(t),\varphi(t,\bar z(t)),\bar z_0(t),\mu_t] dt + \sigma [t,\bar
z(t),\bar z_0(t),\mu_t] dw(t), 
\end{align*}
with initial condition $(\bar z_0(0),\bar z(0))$, where for an arbitrary function $g\in C(\mathbb{R}^s)$ for appropriate $s$, and probability distribution $\mu_t$ in $\mathbb{R}^n$ we set
\begin{align*}
& g[t,z,\varphi,z_0,\mu_t]=\int_{\mathbb{R}^n} g[t,z,\phi,z_0,x]\mu_t(dx),
\end{align*}
when the indicated integral converges. In using the MV system it is assumed that the infinite population of minor agents can be modelled by the collection of sample paths of individual agents subject to their individual initial conditions and their individual Brownian sample paths. 

In the above MV system $\big(\bar z_0(\cdot), \bar z(\cdot),\mu_{(\cdot)}\big)$ is a ``consistent solution'' if $\big(\bar z_0(\cdot), \bar z(\cdot)\big)$ is a solution to the above MV system, $\mu_{t}$, $0 \leq t \leq T$, is the conditional law of $\bar z(t)$ given $\mathcal F_t^{w_0}$ (i.e., $\mu_t := \mathcal L \big(\bar z(t)|\mathcal F_t^{w_0}\big)$). 

Under ({\bf A4})-({\bf A5}) and ({\bf H4}) it can be shown by a fixed point argument that there exists a unique solution $\big(\bar z_0(\cdot), \bar z(\cdot),\mu_{(\cdot)}\big)$ to the above system (see Theorem 1.1 in \cite{sznitman1991topics} or Theorem \ref{CNP:Analysis:Minor:SMV} below).

We also introduce the equations
\begin{align*}
& d\bar z_0(t) = f_0[t,\bar z_0(t),\varphi_0(t,\bar z_0(t)),\mu_t]dt +
\sigma_0[t,\bar z_0(t),\mu_t]dw_0(t), \quad 0 \leq t \leq T, \\
& d\bar z_i(t) = f[t,\bar z_i(t),\varphi(t,\bar z_i(t)),\bar z_0(t),\mu_t]dt +
\sigma[t,\bar z_i(t),\bar z_0(t),\mu_t] dw_i(t), \quad 1 \leq i \leq N,
\end{align*}
with initial conditions $\bar z_j(0)=z_j(0)$, $0 \leq j \leq N$, which can be viewed as
$N$ independent samples of the MV system above. We develop a decoupling result below such that each $\hat z_i^N$, $1 \leq i \leq N$, has the natural limit $\bar z_i$ in the infinite population limit (see Theorem 12 in \cite{HuangMalhameCaines_06}). 

The proof of the following theorem, which is based on the Cauchy-Schwarz inequality, Gronwall's lemma and the conditional independence of minor agents given $\mathcal F_t^{w_0}$, is given in Appendix A in \cite{MNourian_SiamAppendix}. 

\begin{theorem}\label{Theorem:MCT}[McKean-Vlasov Convergence Result] Assume ({\bf A1}), ({\bf A3})-({\bf A5}) and ({\bf H4}) hold.
Then we have
\begin{align}
\label{MCT:Result}
& \sup_{0 \leq j \leq N} \sup_{0 \leq t \leq T} \E|\hat z^N_j(t)-\bar z_j(t)|=O(1/\sqrt N),
\end{align}
where the right hand side may depend upon the terminal time $T$. \qed
\end{theorem}

\section{A Preliminary Nonlinear Stochastic Optimal Control Problem with Random
Coefficients} \label{S3:mm:SHJB} Let $(W(t))_{t \geq 0}$ and $(B(t))_{t \geq 0}$ be mutually
independent standard Brownian motions in $\mathbb{R}^m$, with $\mathcal F_t^{W,B} :=
\sigma\{W(s), B(s): s \leq t\}$ and $\mathcal F_t^W := \sigma\{W(s): s \leq t\}$ where both are augmented by all the $\Prob$-null sets in $\mathcal F$. 

We now consider the following single agent nonlinear
stochastic optimal control problem (SOCP) on $(\Omega, \mathcal F,\{\mathcal
F_t\}_{t\geq 0},\Prob)$:
\begin{align}
& dz(t,\omega) = f[t,\omega,z,u] dt + \sigma[t,\omega,z] dW(t) + \varsigma[t,\omega,z]
dB(t), \quad 0 \leq t \leq T, \label{CNP:GenAgeDyn} \\
& \inf_{u \in \mathcal U} J(u) := \inf_{u \in \mathcal U} \E \Big[\int_0^T L[t,\omega,z(t),u(t)]
dt\Big],
\label{CNP:GenAgeCos}
\end{align}  
where the coefficients $f, \sigma, \varsigma$ and $L$ are random depending on $\omega \in
\Omega$ explicitly. In (\ref{CNP:GenAgeDyn})-(\ref{CNP:GenAgeCos}): (i) $z:[0,T] \times
\Omega \rightarrow \mathbb{R}^n$ is the state of the agent with $\mathcal F_0^{W,B}$-adapted random
initial state $z(0)$ such that $\E|z(0)|^2 <\infty$; (ii) $u:[0,T] \times
\Omega \rightarrow U$ is the control input where $U$ is a compact metric space; (iii)
the functions $f:[0,T]\times\Omega\times\mathbb{R}^n\times U
\rightarrow \mathbb{R}^n$, $\sigma, \varsigma:[0,T]\times\Omega\times\mathbb{R}^n\rightarrow
\mathbb{R}^{n \times m}$ are $\mathcal F_t^W$-adapted stochastic processes; (iv) the admissible control set
$\mathcal U$ is taken as $\mathcal U:=\big\{u(\cdot) \in U: u(t)$ is adapted to $\sigma$-field
$\mathcal F_t^{W,B}$ and $\E\int_0^T |u(t)|^2 dt <\infty\big\}$. We introduce the following assumptions (see \cite{peng1992stochastic}).

({\bf H1}) $f[t,x,u]$ and $L[t,x,u]$ are a.s. continuous in $(x,u)$ for each $t$, a.s. continuous in $t$ for
each $(x,u)$, $f[t,0,0] \in L_{\mathcal F_t}^2([0,T];\mathbb{R}^n)$ and $L[t,0,0] \in L_{\mathcal
F_t}^2([0,T];\mathbb{R}_+)$. In addition, they and all their first derivatives (w.r.t. $x$) are a.s. continuous and bounded.

({\bf H2}) $\sigma[t,x]$ and $\varsigma[t,x]$ are a.s. continuous in $x$ for each $t$, a.s. continuous in $t$ for each $x$ and $\sigma[t,0]$, $\varsigma[t,0] \in L_{\mathcal F_t}^2([0,T];\mathbb{R}^{n \times m})$. In addition, they and all their first derivatives (w.r.t. $x$) are a.s. continuous and bounded.

({\bf H3}) ({\it Non-degeneracy Assumption}) There exist non-negative constants $\alpha_1$ and $\alpha_2$ such that
\begin{align*}
& \sigma[t,\omega,x]\sigma^T[t,\omega,x] \geq \alpha_1 I, \quad \varsigma[t,\omega,x]\varsigma^T(t,\omega,x) \geq \alpha_2 I,  \quad a.s., \quad \forall (t,\omega,x),
\end{align*}
where $\alpha_1$ or $\alpha_2$ (but not both) can be zero.

The value function for the SOCP (\ref{CNP:GenAgeDyn})-(\ref{CNP:GenAgeCos}) is defined by (see
\cite{peng1992stochastic})
\begin{align}
\label{PengValueFun}
& \phi\big(t,x(t)\big) = \inf_{u \in \mathcal U} \E_{\mathcal F_t^W} \int_t^T L[s,\omega,z(s),u(s)] ds,
\end{align}
where $x(t)$ is the initial condition for the process $z(\cdot)$. We note that $\phi\big(t,x(t)\big)$ is an $\mathcal F_t^W$-adapted process which is sample path continuous a.s. under the assumptions ({\bf H1})-({\bf H2}). We assume that there exists an optimal control law $u^o \in \mathcal U$ such that
\begin{align*}
& \phi\big(t,x(t)\big) = \E_{\mathcal F_t^W} \int_t^T L[s,\omega,x(s),u^o(s,\omega,x(s))] ds,
\end{align*}
where $x(\cdot)$ is the closed-loop solution when the control law $u^o$ is applied. By the Principle of Optimality, it can be shown that the process
\begin{align} \label{MRT:1}
& \zeta(t) := \phi\big(t,x(t)\big) + \int_0^t L[s,\omega,x(s),u^o(s,x(s))] ds,
\end{align}
is an $\{\mathcal F_t^W\}_{ 0 \leq t \leq T}$-martingale (see \cite{boel1977optimal}). Next, by the martingale representation theorem (see Theorem 5.7, Chapter 1, \cite{yong1999stochastic}) along the optimal solution $x(\cdot)$ there exists an $\mathcal F_t^W$-adapted process $\psi\big(\cdot,x(\cdot)\big)$ such that
\begin{align} \label{MRT:2}
&  \zeta(t) =  \phi\big(0,x(0)\big) + \int_0^t \psi^T(s,x(s))dW(s), \qquad t \in [0,T].
\end{align}     
From (\ref{MRT:1})-(\ref{MRT:2}) and the fact that $\phi(T,x(T))=0$, it follows that
\begin{align*}
& \zeta(T) = \int_0^T L[s,\omega,x(s),u^o(s,x(s))] ds = \phi\big(0,x(0)\big) + \int_0^T \psi^T(s,x(s))dW(s),
\end{align*}
which gives 
\begin{align} \label{MRT:3}
& \phi(0,x(0)) = \int_0^T L[s,\omega,x(s),u^o(s,x(s))] ds - \int_0^T \psi^T(s,x(s))dW(s).
\end{align}
Hence, combining (\ref{MRT:1})-(\ref{MRT:3}) yields   
\begin{align}
\label{SMR}
& \phi\big(t,x(t)\big) =  \int_t^T L[s,\omega,x(s),u^o(s,x(s))] ds - \int_t^T \psi^T\big(s,x(s)\big) dW(s) \\
& \hspace{1.4cm} =: \int_t^T\Gamma\big(s,x(s)\big) ds - \int_t^T \psi^T\big(s,x(s)\big) dW(s), \quad t \in [0,T], \notag
\end{align}
where $\phi\big(s,x(s)\big)$, $\Gamma\big(s,x(s)\big)$ and $\psi\big(s,x(s)\big)$ are $\mathcal F_s^W$-adapted
stochastic processes (see the assumed semi-martingale representation form (3.5) in \cite{peng1992stochastic}).

Using the extended It\^o-Kunita formula (see Appendix B in \cite{MNourian_SiamAppendix}) and the Principle of Optimality, Peng \cite{peng1992stochastic} showed that since $\phi(t,x)$ can be expressed in the semi-martingale form
(\ref{SMR}), and if $\phi(t,x)$, $\psi(t,x)$, $D_x\phi(t,x)$, $D^2_{xx}\phi(t,x)$ and $D_x \psi(x,t)$ are a.s. continuous in $(x,t)$, then the pair $\big(\phi(s,x),\psi(s,x)\big)$
satisfies the following backward in time stochastic Hamilton-Jacobi-Bellman (SHJB) equation: 
\begin{align}
& \label{CNP:SHJB}
-d\phi(t,\omega,x)=\Big[H[t,\omega,x,D_x\phi(t,\omega,x)]+\big<\sigma[t,\omega,x],D_x\psi(t
,\omega,x)\big>\\ 
& \qquad +\frac{1}{2}\tr\big(a[t,\omega,x]D_{xx}^2 \phi(t,\omega,x)\big)\Big]dt
-\psi^T(t,\omega,x)dW(t,\omega), \quad \phi(T,x) =0,\notag
\end{align}
where $(t,x) \in [0,T] \times\mathbb{R}^n$,
$a[t,\omega,x]:=\sigma[t,\omega,x]\sigma^T[t,\omega,x]+\varsigma[t,\omega,x]\varsigma^T[t,\omega,x]$
, and the stochastic
Hamiltonian $H:[0,T] \times\Omega\times\mathbb{R}^n\times\mathbb{R}^n\rightarrow \mathbb{R}$ is
given by
\begin{align*}
& H[t,\omega,x,p]:=\inf_{u\in \mathcal U}\big\{\big<f[t,\omega,x,u], p\big> + L[t,\omega,x,u]\big\}.
\end{align*}
We note that the appearance of the term $\big<\sigma[t,\omega,x],D_x\psi(t,\omega,x)\big>$ in
equation (\ref{CNP:SHJB}) corresponds to the Brownian motion $W(\cdot)$ in the
extended It\^o-Kunita formula (\ref{ItoKunita}) for the composition of $\mathcal F_t^{W}$-adapted stochastic
processes $\phi(t,\omega,x)$ and $z(t,\omega)$ given in (\ref{SMR}) and (\ref{CNP:GenAgeDyn}),
respectively. 

The solution to the backward in time SHJB equation (\ref{CNP:SHJB}) is a unique forward in
time $\mathcal F_t^W$-adapted pair $(\phi,\psi)(t,x) \equiv
\big(\phi(t,\omega,x),\psi(t,\omega,x)\big)$ (see \cite{peng1992stochastic,yong1999stochastic}). We omit the proof of the following theorem which closely resembles that of Theorem 4.1 in \cite{peng1992stochastic}.

\begin{theorem} \label{CNP:PenTheorem} Assume ({\bf H1})-({\bf H3}) hold. Then the SHJB equation (\ref{CNP:SHJB}) has a unique solution $(\phi(t,x),\psi(t,x))$ in $\big(L_{\mathcal F_t}^2([0,T];\mathbb{R}),L_{\mathcal F_t}^2([0,T];\mathbb{R}^m)\big)$. \qed
\end{theorem}

The forward in time $\mathcal F_t^W$-adapted optimal control process of the SOCP
(\ref{CNP:GenAgeDyn})-(\ref{CNP:GenAgeCos}) is given by (see \cite{peng1992stochastic}) 
\begin{align}
\label{CNP:OC}
& u^o(t,\omega,x) :=\arg\inf_{u\in U} H^u[t,\omega,x,D_x\phi(t,\omega,x),u] \\
& \qquad \qquad ~~ = \arg\inf_{u\in U}\big\{\big<f[t,\omega,x,u],D_x\phi(t,\omega,x)\big> +
L[t,\omega,x,u]\big\}. \notag
\end{align} 

By a verification theorem approach, Peng \cite{peng1992stochastic} showed that if the unique solution $(\phi,\psi)(t,x)$ of the SHJB equation (\ref{CNP:SHJB}) satisfies:

(i) for each $t$, $(\phi,\psi)(t,\cdot)$ is a $C^2(\mathbb{R}^n)$ map from $\mathbb{R}^n$ into $\mathbb{R} \times
\mathbb{R}^m$,

(ii) for each $x$, $(\phi,\psi)(t,x)$ and $(D_x \phi,D^2_{xx} \phi, D_x \psi)(t,x)$ are
continuous $F_t^W$-adapted stochastic processes, then $\phi(x,t)$ coincides with the value function (\ref{PengValueFun}) of the SOCP
(\ref{CNP:GenAgeDyn})-(\ref{CNP:GenAgeCos}).

\section{The Major and Minor Agent Stochastic Mean Field Game System} \label{S5:mm:SMF} In the formulation
(\ref{CNP:Major:GenDyn})-(\ref{CNP:Minor:GenCost}) all minor agents are reacting to the same major
agent and hence the major agent has non-negligible influence on the mean field behaviour of the
minor agents. In other words, the noise process of the major agent $w_0$ causes random fluctuation
of the mean-field behaviour of the minor agents and makes it stochastic (see the discussion in
Section 2 of \cite{HuangSIAM2010} for the major and minor agent MFG LQG model).

In this section, we first construct two auxiliary stochastic optimal control problems (SOCP) with random coefficients for the major and a generic minor agent in Sections \ref{S5:mm:SMF:major} and \ref{S5:mm:SMF:minor}, respectively. Then, we present the stochastic mean field system for the major and minor agents game formulation (\ref{CNP:Major:GenDyn})-(\ref{CNP:Minor:GenCost}) via the mean field game consistency condition in Section \ref{MFCC}. 

\subsection{Stochastic Optimal Control Problem of the Major Agent} \label{S5:mm:SMF:major}

By the McKean-Vlasov convergence result in Theorem \ref{Theorem:MCT} which indicates that a single minor agent's statistical properties can effectively approximate the empirical distribution produced by all minor agents, we may approximate the empirical distribution of minor agents $\delta^N_{(\cdot)}$ with a stochastic probability measure $\mu_{(\cdot)}$ which depends on the noise process of the major agent $w_0$.

In this section, let $\mu_t(\omega)$, $0 \leq t \leq T$, be an exogenous nominal minor agent stochastic measure process such that $\mu_0(dx):=dF(x)$ where $F$ is defined in ({\bf A6.2}). Note that in Section \ref{MFCC} $\mu_t(\omega)$ will be characterized via the mean field game consistency condition as the random measure of minor agents' mean field behaviour.

We define the following SOCP
(\ref{CNP:GenAgeDyn})-(\ref{CNP:GenAgeCos}) with $\mathcal F_t^{w_0}$-adapted random coefficients
from the major agent's model (\ref{CNP:Major:GenDyn}) and (\ref{CNP:Major:GenCost}) in the infinite
population limit:
\begin{align}
& ~ dz_0(t) = f_0[t,z_0(t),u_0(t),\mu_t(\omega)]dt +
\sigma_0[t,z_0(t),\mu_t(\omega)]dw_0(t,\omega), \quad z_0(0),
\label{CNP:MajorinfDyn} \\
& \inf_{u_0 \in \mathcal U_0} J_0 (u_0) := \inf_{u_0 \in \mathcal U_0} \E \Big[\int_0^T
L_0[t,z_0(t),u_0(t),\mu_t(\omega)] dt\Big], \label{CNP:MajorinfCos}
\end{align}
where we explicitly indicate the dependence of the random measure $\mu_{(\cdot)}$ on the sample
point $\omega \in \Omega$. 

Step I ({\it Major Agent's Stochastic Hamilton-Jacobi-Bellman (SHJB) Equation}):

The value function of the major agent's SOCP (\ref{CNP:MajorinfDyn})-(\ref{CNP:MajorinfCos}) is defined by
\begin{align}
\label{CNP:ValueFunction:Major}
& \phi_0\big(t,x(t)\big) = \inf_{u_0 \in \mathcal U_0} \E_{\mathcal F_t^{w_0}} \int_t^T L_0[s,z_0(s),u_0(s),\mu_s(\omega)]ds,
\end{align}
where $x(t)$ is the initial condition for the process $z_0(s)$ (see (\ref{PengValueFun})). As in Section \ref{S3:mm:SHJB}, $\phi_0\big(t,x(t)\big)$ has the form (see (\ref{SMR}))
\begin{align*}
& \phi_0\big(t,x(t)\big) = \int_t^T\Gamma_0\big(s,x(s)\big) ds - \int_t^T \psi^T_0\big(s,x(s)\big) dw_0(s), \quad t \in [0,T],
\end{align*}
where $\phi_0\big(s,x(s)\big)$, $\Gamma_0\big(s,x(s)\big)$ and $\psi_0\big(s,x(s)\big)$ are $\mathcal F_s^{w_0}$-adapted stochastic processes. If $\phi_0(t,x)$, $\psi_0(t,x)$, $D_x\phi_0(t,x)$, $D^2_{xx}\phi_0(t,x)$ and $D_x \psi_0(x,t)$ are a.s. continuous in $(x,t)$, then the pair $\big(\phi_0(s,x),\psi_0(s,x)\big)$ satisfies the following stochastic Hamilton-Jacobi-Bellman (SHJB) equation:
\begin{align}
\label{CNP:SHJB1:Major}
& -d\phi_0(t,\omega,x)=\Big[H_0[t,\omega,x,D_x\phi_0(t,\omega,x)]+\big<\sigma_0[t,x,\mu_t(\omega)],
D_x\psi_0(t,\omega,x)\big>\\ 
& ~ +\frac{1}{2}\tr\big(a_0[t,\omega,x]D_{xx}^2 \phi_0(t,\omega,x)\big)\Big]dt
-\psi^T_0(t,\omega,x)dw_0(t,\omega), \quad \phi_0(T,x) =0, \notag
\end{align}
where $(t,x) \in [0,T] \times\mathbb{R}^n$,
$a_0[t,\omega,x]:=\sigma_0[t,x,\mu_t(\omega)]\sigma^T_0[t,x,\mu_t(\omega)]$, and the stochastic
Hamiltonian $H_0:[0,T] \times\Omega\times\mathbb{R}^n\times\mathbb{R}^n\rightarrow \mathbb{R}$ is
given by
\begin{align*}
& H_0[t,\omega,x,p]:= \inf_{u \in \mathcal U_0}\big\{\big<f_0[t,x,u,\mu_t(\omega)], p\big> +
L_0[t,x,u,\mu_t(\omega)]\big\}.
\end{align*}
The solution to the backward in time SHJB equation (\ref{CNP:SHJB1:Major}) is a forward in time
$\mathcal F_t^{w_0}$-adapted pair
$\big(\phi_0(t,x),\psi_0(t,x)\big)\equiv\big(\phi_0(t,\omega,x),\psi_0(t,\omega,x)\big)$ (see \cite{peng1992stochastic}).

We note that the appearance of the term $\big<\sigma_0[t,x,\mu_t(\omega)],
D_x\psi_0(t,\omega,x)\big>$ in equation (\ref{CNP:SHJB1:Major}) corresponds to the major agent's Brownian motion $w_0(\cdot)$ in the extended It\^o-Kunita formula (\ref{ItoKunita}) for the composition of $\mathcal F_t^{w_0}$-adapted processes $\phi_0(t,\omega,x)$ and $z_0(t,\omega)$ in (\ref{CNP:MajorinfDyn}).

The best response process of the major agent's SOCP (\ref{CNP:MajorinfDyn})-(\ref{CNP:MajorinfCos}) is given by
\begin{align}
\label{CNP:BR1:Major}
& u^o_0(t,\omega,x) \equiv u^o_0(t,x|\{\mu_{s}(\omega)\}_{0\leq s\leq T}) :=\arg\inf_{u_0\in U_0}
H_0^{u_0}[t,\omega,x,u_0,D_x\phi_0(t,\omega,x)] \\
& ~ \qquad \equiv \arg\inf_{u_0\in U_0} \big\{\big<f_0[t,x,u_0,\mu_t(\omega)],
D_x\phi_0(t,\omega,x)\big> + L_0[t,x,u_0,\mu_t(\omega)]\big\}, \notag
\end{align}
where the infimum exists a.s. here and in all analogous infimizations in the chapter due to the continuity of all functions appearing in $H_0^{u_0}$ and the compactness of $U_0$. It should be noted that the stochastic best response control $u^o_0$ is a forward in time $\mathcal F_t^{w_0}$-adapted process which depends on the Brownian motion $w_0$ via the stochastic
measure $\mu_t(\omega)$, $0 \leq t \leq T$. The notation in (\ref{CNP:BR1:Major}) indicates that $u^o_0$ at time $t$ depends upon the stochastic measure $\mu_{s}(\omega)$ on the whole interval $0 \leq s \leq T$.

Step II ({\it Major Agent's Stochastic Coefficient McKean-Vlasov (SMV) Equation}): By substituting the best response control process
$u^o_0$ (\ref{CNP:BR1:Major}) into the major agent's dynamics (\ref{CNP:MajorinfDyn}) we get the
following stochastic McKean-Vlasov (SMV) dynamics with random coefficients:
\begin{align}
\label{CNP:Major:CL}
& ~ dz_0^o(t,\omega) = f_0[t,z_0^o,u^o_0(t,\omega,z_0^o),\mu_t(\omega)]dt
+ \sigma_0[t,z_0^o,\mu_t(\omega)]dw_0(t,\omega),
\end{align}
with $z_0^o(0)=z_0(0)$, where $f_0$ and $\sigma_0$ are random processes via the stochastic measure $\mu$ and $u^o_0$.

\subsection{Stochastic Optimal Control Problem of the Generic Minor Agent} \label{S5:mm:SMF:minor}
As in Section \ref{S5:mm:SMF:major} let $\mu_{t}$, $0 \leq t \leq T$, be the exogenous nominal minor agent stochastic measure process approximating the empirical distribution produced by all minor agents in the infinite population limit such that $\mu_0(dx)=dF(x)$ where $F$ is defined in ({\bf A6.2}). We let $z_0^o(\cdot)$ be the solution to the major agent's SMV equation (\ref{CNP:Major:CL}).

We define the following SOCP (\ref{CNP:GenAgeDyn})-(\ref{CNP:GenAgeCos}) with $\mathcal
F_t^{w_0}$-adapted random coefficients from the $i^{\textrm{th}}$ generic minor agent's model
(\ref{CNP:Minor:GenDyn}), (\ref{CNP:Minor:GenCost}) in the infinite population limit:
\begin{align}
\label{CNP:MinorinfDyn}
& ~ dz_i(t) = f[t,z_i(t),u_i(t),z_0^o(t,\omega),\mu_t(\omega)] dt +
\sigma[t,z_i(t),z_0^o(t,\omega),\mu_t(\omega)] dw_i(t), \\
\label{CNP:MinorinfCos}
& \inf_{u_i \in \mathcal U} J_i(u_i) := \inf_{u_i \in \mathcal U} \E \Big[\int_0^T
L[t,z_i(t),u_i(t),z_0^o(t,\omega),\mu_t(\omega)]dt\Big], \quad z_i(0), 
\end{align}
where we explicitly indicate the dependence of the solution to the major agent's SMV equation $z_0^o(\cdot)$ and the nominal minor agent's random measure $\mu_{(\cdot)}$ on the sample point $\omega \in \Omega$.

Step I ({\it Generic Minor Agent's Stochastic Hamilton-Jacobi-Bellman (SHJB) Equation}):

The value function of the generic minor agent's SOCP (\ref{CNP:MinorinfDyn})-(\ref{CNP:MinorinfCos}) is defined by
\begin{align}
\label{CNP:ValueFunction:Minor}
& \phi_i\big(t,x(t)\big) = \inf_{u_i \in \mathcal U_0} \E_{\mathcal F_t^{w_0}} \int_t^T L[s,z_i(s),u_i(s),z_0^o(s,\omega),\mu_s(\omega)]ds,
\end{align}
where $x(t)$ is the initial condition for the process $z_i(\cdot)$. As in Section \ref{S3:mm:SHJB}, $\phi_i\big(t,x(t)\big)$ has the form (see (\ref{SMR}))
\begin{align*}
& \phi_i\big(t,x(t)\big) = \int_t^T\Gamma_i\big(s,x(s)\big) ds - \int_t^T \psi^T_i \big(s,x(s)\big)  dw_0(s), \quad t \in [0,T],
\end{align*}
where $\phi_i\big(s,x(s)\big)$, $\Gamma_i\big(s,x(s)\big)$ and $\psi_i\big(s,x(s)\big)$ are $\mathcal F_s^{w_0}$-adapted stochastic processes. If $\phi_i(t,x)$, $\psi_i(t,x)$, $D_x\phi_i(t,x)$ and $D^2_{xx}\phi_i(t,x)$ are a.s. continuous in $(x,t)$, then the pair $\big(\phi_i(s,x),\psi_i(s,x)\big)$ satisfies the following backward in time stochastic Hamilton-Jacobi-Bellman (SHJB) equation (see (\ref{CNP:SHJB})):
\begin{align}
\label{CNP:SHJB1:Minor}
&
-d\phi_i(t,\omega,x)=\Big[H[t,\omega,x,D_x\phi_i(t,\omega,x)]+\frac{1}{2}\tr\big(a[t,\omega,x]
D_{xx}^2 \phi_i(t,\omega,x)\big)\Big]dt \\
& \qquad \qquad \quad \qquad -\psi^T_i(t,\omega,x)dw_0(t,\omega), \quad \phi_i(T,x) =0,  \notag
\end{align}
where $(t,x) \in [0,T] \times\mathbb{R}^n$,
$a[t,\omega,x]:=\sigma[t,x,z_0^o(t,\omega),\mu_t(\omega)]\sigma^T[t,x,z_0^o(t,\omega),\mu_t(\omega)]
$, and the stochastic Hamiltonian $H:[0,T]
\times\Omega\times\mathbb{R}^n\times\mathbb{R}^n\rightarrow \mathbb{R}$ is
given by
\begin{align*}
& H[t,\omega,x,p]:= \inf_{u \in \mathcal U}\big\{\big<f[t,x,u,z_0^o(t,\omega),\mu_t(\omega)], p\big>
+ L[t,x,u,z_0^o(t,\omega),\mu_t(\omega)]\big\}.
\end{align*}
The solution to the backward in time SHJB equation (\ref{CNP:SHJB1:Minor}) is a forward in time
$\mathcal F_t^{w_0}$-adapted
pair $\big(\phi_i(t,x),\psi_i(t,x)\big)\equiv\big(\phi_i(t,\omega,x),\psi_i(t,\omega,x)\big)$ (see \cite{peng1992stochastic}). We note that since the coefficients of the SOCP (\ref{CNP:MinorinfDyn})-(\ref{CNP:MinorinfCos}) are $\mathcal F_t^{w_0}$-adapted random processes we have the major agent's Brownian motion $w_0$ in (\ref{CNP:SHJB1:Minor}) which allows us to seek for a forward in time adapted solution to the backward in time SHJB equation (\ref{CNP:SHJB1:Minor}).

It is important to note that in (\ref{CNP:SHJB1:Minor}) unlike the major agent's SHJB equation
(\ref{CNP:SHJB1:Major}) we do not have the term
$\big<\sigma[t,x,z_0^o(t,\omega),\mu_t(\omega)]D_x\psi_i(t,\omega,x)\big>$ since the coefficients in the
minor agent's model (\ref{CNP:MinorinfDyn})-(\ref{CNP:MinorinfCos}) are $\mathcal F_t^{w_0}$-adapted
random processes depending upon the major agent's Brownian motion $(w_0)$ which is independent
of the minor agent's Brownian motion $(w_i)$ (see the extended It\^o-Kunita formula (\ref{ItoKunita})).

As in Section \ref{S5:mm:SMF:major}, the stochastic best response process of the minor agent's SOCP
(\ref{CNP:MinorinfDyn})-(\ref{CNP:MinorinfCos}) is
\begin{align}
\label{CNP:BR1:Minor}
& u^o_i(t,\omega,x) \equiv u^o_i(t,x|\{z_0^o(s,\omega),\mu_{s}(\omega)\}_{0\leq s\leq
T}) :=\arg\inf_{u\in U}
H^u[t,\omega,x,u,D_x\phi_i(t,\omega,x)] \\
& ~  \equiv \arg\inf_{u\in U} \big\{\big<f[t,x,u,z_0^o(t,\omega),\mu_t(\omega)],
D_x\phi_i(t,\omega,x)\big> + L[t,x,u,z_0^o(t,\omega),\mu_t(\omega)]\big\}, \notag
\end{align}
where the infimum exists a.s. here and in all analogous infimizations in the chapter due to the continuity of all functions appearing in $H^{u}$ and the compactness of $U$. It should be noted that the stochastic best response process of the generic minor agent $u^o_i$ is a forward in time $\mathcal
F_t^{w_0}$-adapted random process which depends on the Brownian motion $w_0$ via the major agent's state $z_0^o(t,\omega)$ and the stochastic measures $\mu_t(\omega)$, $0 \leq t \leq T$. The notation in (\ref{CNP:BR1:Minor}) indicates that $u^o_i$ at time $t$ depends upon $z_0^o(s,\omega)$ and $\mu_{s}(\omega)$ on the whole interval $0 \leq s \leq T$.

Step II ({\it Minor Agent's Stochastic Coefficient McKean-Vlasov (SMV) and Stochastic Coefficient Fokker-Planck-Kolmogorov (SFPK) Equations}): By substituting the best response control process
$u^o_i$ (\ref{CNP:BR1:Minor}) into the minor agent's dynamics (\ref{CNP:MinorinfDyn}) we get the
following stochastic McKean-Vlasov (SMV) dynamics with random coefficients:
\begin{align}
\label{CNP:Minor:CL}
& ~ dz_i^o(t,\omega,\omega') = f[t,z_i^o,u^o_i(t,\omega,z_i),z_0^o(t,\omega),\mu_t(\omega)]dt \\
& \hspace{3cm} + \sigma[t,z_i^o,z_0^o(t,\omega),\mu_t(\omega)]dw_i(t,\omega'),
\quad z_i^o(0)=z_i(0), \notag
\end{align}
where $f$ and $\sigma$ are random processes via $z_0^o$, $\mu$, and the best response control process $u^o_i$ which all depend on the Brownian motion of the major agent $(w_0)$. 

Based on the McKean-Vlasov approximation in Section \ref{S4:mm:MFCT}, the generic agent's statistical properties can effectively approximate the empirical distribution produced by all minor agents in a large population system.
Hence, we obtain a new stochastic measure $\hat \mu_t(\omega)$ for the mean field behaviour of
minor agents as the conditional law of the generic minor agent's process $z_i^o(t,\omega)$ given $\mathcal F_t^{w_0}$.
We characterize $\hat \mu_t(\omega)$, $0 \leq t \leq T$, by
$P(z_i^o(t,\omega)\leq\alpha|\mathcal F_t^{w_0})=\int_{-\infty}^\alpha \hat \mu(t,\omega,dx)$
a.s. for all $\alpha \in \mathbb{R}^n$ and $0 \leq t \leq T$, with $\hat
\mu_0(dx)=\mu_0(dx)=dF(x)$ where $F$ is defined in ({\bf A6.2}).

An equivalent method to characterize the SMV of the generic minor agent is to express
(\ref{CNP:Minor:CL}) in the form of stochastic Fokker-Planck-Kolmogorov (SFPK) equation with random
coefficients:
\begin{align}
\label{CNP:SFPK1:Minor}
& d\hat p(t,\omega,x)=
\Big(-\big<D_x,f[t,x,u^o_i(t,\omega,x),z_0^o(t,\omega),\mu_t(\omega)]\hat p(t,\omega,x)\big> \\
& \qquad \qquad \qquad +\frac{1}{2} \tr
\big<D^2_{xx},a[t,\omega,x]\hat p(t,\omega,x)\big>\Big) dt, \quad \hat p(0,x)=p_0(x), \notag
\end{align}
in $[0,T] \times\mathbb{R}^n$ where $p(t,\omega,x)$ is the conditional probability density of $z_i^o(t,\omega)$ given $\mathcal F_t^{w_0}$. By the the McKean-Vlasov approximation (see Section \ref{S4:mm:MFCT}) it is possible to characterize the mean field behaviour of minor agents in terms of generic agent's density function $\hat p(t,\omega,x)$. The reason that the generic minor agent's FPK equation (\ref{CNP:SFPK1:Minor}) does not include the It\^o integral term with respect to $w_i$ is due to the fact that $p(t,\omega,x)$ is the conditional probability density given $\mathcal F_t^{w_0}$, and the independence of the Brownian motions $w_0$ and $w_i$, $1 \leq i \leq N$.
 
The density function $\hat p(t,\omega,x)$ generates the
random measure of the minor agent's mean field behaviour $\hat \mu_t(\omega)$ such that
$\hat \mu(t,\omega,dx)=\hat p(t,\omega,x)dx$ (a.s.), $0 \leq t \leq T$.  

We note that the major agent's SOCP (\ref{CNP:MajorinfDyn})-(\ref{CNP:MajorinfCos}) and minor agent's SOCP (\ref{CNP:MinorinfDyn})-(\ref{CNP:MinorinfCos}) may be written with respect to the random density $p(t,\omega,x)$ of the stochastic measure $\mu(t,\omega,dx)$ by $\mu(t,\omega,dx)=p(t,\omega,x)dx$ (a.s.), $0 \leq t \leq T$.

\subsection{The Mean Field Game Consistency Condition} \label{MFCC}
Based on the mean field game (MFG) or Nash certainty equivalence (NCE) consistency (see \cite{HuangMalhameCaines_06} and \cite{Lasry07}), we close the ``measure and control'' mapping loop by setting $\hat \mu_{t}(\omega) = \mu_{t}(\omega)$ a.s., $0 \leq t \leq T$, or $\hat p(t,\omega,x) = p(t,\omega,x)$ a.s. for $(t,x) \in [0,T]\times\mathbb{R}^n$. The
MFG consistency is demonstrated in: (i) the major agent's stochastic mean field game (SMFG) system
\begin{align}
\label{CNP:SHJB:Major}
& [\textrm{MFG-SHJB}]  \quad -d\phi_0(t,\omega,x)=\Big[H_0[t,\omega,x,D_x\phi_0(t,\omega,x)] \\
&  \hspace{1cm} 
+\big<\sigma_0[t,x,\mu_t(\omega)],D_x\psi_0(t,\omega,x)\big>+\frac{1}{2}\tr\big(a_0[t,\omega,x]D_{
xx}^2 \phi_0(t,\omega,x)\big)\Big]dt \notag\\
& \hspace{1cm}  -\psi^T_0(t,\omega,x)dw_0(t,\omega), \quad \phi_0(T,x) =0, \notag \\
\label{CNP:BR:Major}
& u^o_0(t,\omega,x) \equiv u^o_0(t,x|\{\mu_{s}(\omega)\}_{0\leq s\leq T}) \\
&  \hspace{1cm}  :=\arg\inf_{u_0\in U_0} \big\{\big<f_0[t,x,u_0,\mu_t(\omega)],
D_x\phi_0(t,\omega,x)\big> + L_0[t,x,u_0,\mu_t(\omega)]\big\},  \notag \\
& \label{CNP:SMV:Major}  [\textrm{MFG-SMV}] \qquad dz_0^o(t,\omega) = f_0[t,z_0^o,u^o_0(t,\omega,z_0^o),\mu_t(\omega)]dt
\\
& \hspace{1cm} + \sigma_0[t,z_0^o,\mu_t(\omega)]dw_0(t,\omega), \quad z_0^o(0)=z_0(0), \notag 
\end{align}
together with (ii) the minor agents' SMF system
\begin{align}
\label{CNP:SHJB:Minor}
& [\textrm{MFG-SHJB}] \quad -d\phi(t,\omega,x)=\Big[H[t,\omega,x,D_x\phi(t,\omega,x)] \\
& \hspace{1cm} +\frac{1}{2}\tr\big(a[t,\omega,x] D_{xx}^2 \phi(t,\omega,x)\big)\Big]dt -\psi^T(t,\omega,x)dw_0(t,\omega), \quad \phi(T,x) =0,  \notag \\
&  \label{CNP:BR:Minor} u^o(t,\omega,x) \equiv u^o(t,x|\{z_0^o(s,\omega),\mu_{s}(\omega)\}_{0\leq
s\leq
T}) \\
& \!\equiv\! \arg\inf_{u\in U} \big\{\big<f[t,x,u,z_0^o(t,\omega),\mu_t(\omega)],
D_x\phi(t,\omega,x)\big> \!+\! L[t,x,u,z_0^o(t,\omega),\mu_t(\omega)]\big\},  \notag \\
&  \label{CNP:SMV:Minor} [\textrm{MFG-SMV}] \qquad dz^o(t,\omega,\omega') =
f[t,z^o,u^o(t,\omega,z^o),z_0^o(t,\omega),\mu_t(\omega)]dt \\
& \hspace{1cm} + \sigma[t,z^o,z_0^o(t,\omega),\mu_t(\omega)]dw(t,\omega'), \notag
\end{align}
where $(t,x) \in [0,T] \times \mathbb{R}^n$, and $z^o(0)$ has the measure $\mu_0(dx)=dF(x)$ where
$F$ is defined in ({\bf A6.2}). We note that in the minor agents' SMFG system (\ref{CNP:SHJB:Minor})-(\ref{CNP:SMV:Minor}) we dropped index $i$ from the generic minor agent's equations (\ref{CNP:MinorinfDyn})-(\ref{CNP:Minor:CL}). The Major and Minor (MM) agent SMFG system is given by (\ref{CNP:SHJB:Major})-(\ref{CNP:SMV:Major}) and (\ref{CNP:SHJB:Minor})-(\ref{CNP:SMV:Minor}). 

The solution of the MM-SMFG system consists of 8-tuple $\mathcal F_t^{w_0}$-adapted random processes
\begin{align*}
&
\big(\phi_0(t,\omega,x),\psi_0(t,\omega,x),u^o_0(t,\omega,x),z_0^o(t,\omega),\phi(t
,\omega,x),\psi(t,\omega,x),u^o(t,\omega,x),z^o(t,\omega)\big),
\end{align*}
where $z^o(t,\omega)$ generates the conditional random law $\mu_t(\omega)$, i.e., $P(z^o(t,\omega)\leq\alpha|\mathcal F_t^{w_0})=\int_{-\infty}^\alpha \mu_t(\omega,dx)$ for all $\alpha \in \mathbb{R}^n$ and $0 \leq t \leq T$. Note that the MM-SMFG equations (\ref{CNP:SHJB:Major})-(\ref{CNP:SMV:Major}) and (\ref{CNP:SHJB:Minor})-(\ref{CNP:SMV:Minor}) are coupled together through $z_0^o(\cdot)$ and $\mu_{(\cdot)}$.

We observe that the solution to the MM-SMFG system is a ``stochastic mean field'' in contrast to the deterministic mean field of the standard MFG problems with only minor agents considered in \cite{HuangMalhameCaines_06,huang2007invariance,LasryLionsStationair2006,LasryLionsHorizonfini2006,Lasry07}. If the noise process of the major agent vanishes then the MM-SMFG system reduces to a deterministic MFG system (see (6)-(9) in \cite{huang2007invariance}).

For the analysis of next section we denote $\mu_t^0(\omega)$, $0 \leq t \leq T$, as the unit mass random measure concentrated at $z_0^o(t,\omega)$ (i.e., $\mu^0_t(\omega) = \delta_{z_0^o(t,\omega)}$).

\section{Existence and Uniqueness of Solutions to the Major and Minor Stochastic Mean Field Game System} \label{S6:mm:Exis} In this section we establish existence and uniqueness for the solution of the joint major and minor (MM) agents' SMFG system (\ref{CNP:SHJB:Major})-(\ref{CNP:SMV:Major}) and (\ref{CNP:SHJB:Minor})-(\ref{CNP:SMV:Minor}). The analysis is based on providing sufficient conditions for a map that goes from the random measure of minor agents $\mu_{(\cdot)}(\omega)$ back to itself, through the equations (\ref{CNP:SHJB:Major})-(\ref{CNP:SMV:Major}) and (\ref{CNP:SHJB:Minor})-(\ref{CNP:SMV:Minor}), to be a contraction operator on the space of random probability measures (see the diagram below).
\begin{align*}
& \begin{array}[c]{ccccc}
\mu_{(\cdot)}(\omega) &\stackrel{(\ref{CNP:SHJB:Major})}{\longrightarrow}& \big(\phi_0(\cdot,\omega,x),\psi_0(\cdot,\omega,x)\big)&\stackrel{(\ref{CNP:BR:Major})}{\longrightarrow}& u_0^o(\cdot,\omega,x) \\
\uparrow\scriptstyle{(\ref{CNP:SMV:Minor})}&&&&\downarrow\scriptstyle{(\ref{CNP:SMV:Major})}\\
u^o(\cdot,\omega,x)  &\stackrel{(\ref{CNP:BR:Minor})}{\longleftarrow}&\big(\phi(\cdot,\omega,x),\psi(\cdot,\omega,x)\big) &\stackrel{(\ref{CNP:SHJB:Minor})}{\longleftarrow} & \mu_{(\cdot)}^0(\omega) \equiv \delta_{z_0^o(t,\omega)}
\end{array}
\end{align*}

In this section we first introduce some preliminary material about the Wasserstein space of probability measures. Second, we analyze the SHJB and SMV equations of the major agent and minor agents in Sections \ref{CNP:MajorAgentAnalysis} and \ref{CNP:MinorAgentAnalysis}, respectively. Third, the analysis of the joint major and minor agents' SMFG system is carried out in Section \ref{analysisjoint} where the main result is given in Theorem \ref{mm:maintheorem} which provides sufficient conditions for a contraction operator map that goes from the random measure of minor agents $\mu_{(\cdot)}(\omega)$ back to itself.

On the Banach space $C([0,T];\mathbb{R}^n)$ we define the metric $\rho_T(x,y)= \sup_{0\leq t\leq T}|x(t)-y(t)|^2 \wedge 1$, where $\wedge$ denotes minimum. It can be shown that $C_{\rho}:=\big(C([0,T];\mathbb{R}^n),\rho_T\big)$ forms a separable complete metric space (i.e., a Polish space). Let $\mathcal M(C_{\rho})$ be the space of all Borel probability measures $\mu$ on $C([0,T];\mathbb{R}^n)$ such that $\int |x|^2d\mu(x)<\infty$. We also denote $\mathcal M(C_\rho\times C_\rho)$ as the space of probability measures on the product space $C([0,T];\mathbb{R}^n)\times C([0,T];\mathbb{R}^n)$. As in \cite{HuangMalhameCaines_06} the process $x$ is defined to be a generic random process with the sample space $C([0,T];\mathbb{R}^n)$, i.e., $x(t,\omega)=\omega(t)$ for $\omega \in C([0,T];\mathbb{R}^n)$. 

Based on the metric $\rho_T$, we introduce the Wasserstein metric on $\mathcal M(C_\rho)$:
\begin{align*}
& D_T^\rho(\mu,\nu) =\inf_{\gamma \in \Pi(\mu,\nu)} \Big[\int_{C_\rho\times C_\rho} \rho_T(x(\omega_1),x(\omega_2)) d\gamma(\omega_1,\omega_2)\Big]^{1/2},
\end{align*}
where $\Pi(\mu,\nu) \subset \mathcal M(C_\rho\times C_\rho)$ is the set of Borel probability measures $\gamma$ such that $\gamma(A \times C([0,T];\mathbb{R}^n))= \mu(A)$ and $\gamma(C([0,T];\mathbb{R}^n)\times A)= \nu(A)$ for any Borel set $A \in C([0,T];\mathbb{R}^n)$. The metric space $\mathcal M_\rho:=\big(\mathcal M(C_\rho), D_T^\rho\big)$ is a Polish space since $C_{\rho}\equiv\big(C([0,T];\mathbb{R}^n),\rho_T\big)$ is a Polish space.

We also introduce the class $\mathcal M_\rho^\beta$ of stochastic measures in the space $\mathcal M_\rho$ with a.s. H\"older continuity of exponent $\beta$, $0 < \beta<1$ (see Definition 3 in \cite{HuangMalhameCaines_06} for the non-stochastic case).  

\begin{definition} A stochastic probability measure $\mu_t(\omega)$, $0 \leq t \leq T$, in the space $\mathcal M_\rho$ is in $\mathcal M_\rho^\beta$ if $\mu$ is a.s. uniformly H\"older continuous with exponent $0 < \beta <1$, i.e., there exists $\beta \in (0,1)$ and constant $c$ such that for any bounded and Lipschitz continuos function $\phi$ on $\mathbb{R}^n$, 
\begin{align*}
& \big|\int_{\mathbb{R}^n} \phi(x)\mu_t(\omega,dx)-\int_{\mathbb{R}^n} \phi(x)\mu_s(\omega,dx)\big| \leq c(\omega)|t-s|^\beta, \quad a.s., 
\end{align*}
for all $0 \leq s < t \leq T$, where $c$ may depend upon the Lipschitz constant of $\phi$ and the sample point $\omega \in \Omega$. 
\end{definition}

As in \cite{HuangMalhameCaines_06}, we may take $\mu_t$, $0 \leq t \leq T$, to be a Dirac measure at any constant $x \in \mathbb{R}^n$ to show that the set $\mathcal M_\rho^\beta$ is nonempty. We introduce the following assumption.

({\bf A8}) For any $p \in \mathbb{R}^n$ and $\mu,\mu^0 := \delta_{z_0^o}  \in \mathcal M_\rho^\beta$, the sets
\begin{align*}
& S_0(t,\omega,x,p):=\arg\inf_{u_0\in U_0} H_0^{u_0}[t,\omega,x,u_0,p], \\
& S(t,\omega,x,p):=\arg\inf_{u\in U} H^u[t,\omega,x,u,p], 
\end{align*}
where $H_0^{u_0}$ and $H^u$ are respectively defined in (\ref{CNP:BR1:Major}) and (\ref{CNP:BR1:Minor}), are singletons and the resulting $u$ and $u_0$ as functions of $[t,\omega,x,p]$ are a.s. continuous in $t$, Lipschitz continuous in $(x,p)$, uniformly with respect to $t$ and $\mu,\mu^0 \in \mathcal M_\rho^\beta$. In addition, $u_0[t,\omega,0,0]$ and $u[t,\omega,0,0]$ are in the space $L^2_{\mathcal F_t}([0,T];\mathbb{R}^n)$.   

The first part of ({\bf A8}) may be satisfied under suitable convexity conditions with respect to $u_0$ and $u$ (see \cite{HuangMalhameCaines_06}).

\subsection{Analysis of the Major Agent's SMFG System}\label{CNP:MajorAgentAnalysis} Let $\mu_t(\omega)$, $0 \leq t \leq T$, be a fixed stochastic measure in the set $\mathcal M_\rho^\beta$ with $0 < \beta <1$ such that $\mu_0(dx):=dF(x)$ where $F$ is defined in ({\bf A2}). Then, the functionals of $\mu_{(\cdot)}(\omega)$ in (\ref{CNP:MajorinfDyn})-(\ref{CNP:MajorinfCos}) become random functions which we write as
\begin{align}
\label{CNP:f0*}
& f_0^*[t,\omega,z_0,u_0]:=f_0[t,z_0,u_0,\mu_t(\omega)], \quad \sigma_0^*[t,\omega,z_0]:=\sigma_0[t,z_0,\mu_t(\omega)], \\
& L_0^*[t,\omega,z_0,u_0] :=L_0[t,z_0,u_0,\mu_t(\omega)]. \notag
\end{align}

We have the following result which broadly follows Proposition 4 in \cite{HuangMalhameCaines_06}.

\begin{proposition} \label{CNP:Proposition:Major} Assume ({\bf A3}) holds for $U_0$. Let $\mu_t(\omega)$, $0 \leq t \leq T$, be a fixed stochastic measure in the set $\mathcal M_\rho^\beta$ with $0 < \beta <1$. For $f_0^*$, $\sigma_0^*$ and $L_0^*$ defined in (\ref{CNP:f0*}) it is the case that:
\begin{enumerate}[(i)]
\item Under ({\bf A4}) for $f_0$ and $\sigma_0$, the functions $f_0^*[t,\omega,z_0,u_0]$ and $\sigma_0^*[t,\omega,z_0]$ and their first order derivatives (w.r.t $z_0$) are a.s. continuous and bounded on $[0,T] \times \mathbb{R}^n \times U_0$ and $[0,T]\times \mathbb{R}^n$. $f_0^*[t,\omega,z_0,u_0]$ and $\sigma_0^*[t,\omega,z_0]$ are a.s. Lipschitz continuous in $z_0$. In addition, $f_0^*[t,\omega,0,0]$ is in the space $L^2_{\mathcal F_t}([0,T];\mathbb{R}^n)$ and $\sigma_0^*[t,\omega,0]$ is in the space $L^2_{\mathcal F_t}([0,T];\mathbb{R}^{n\times m})$. 

\item Under ({\bf A5}) for $f_0$, the function $f_0^*[t,\omega,z_0,u_0]$ is a.s. Lipschitz continuous in $u_0\in U_0$, i.e., there exist a constant $c>0$ such that
\begin{align*}
& \sup_{t \in [0,T],z_0 \in \mathbb{R}^n} \big|f_0^*[t,\omega,z_0,u_0]-f_0^*[t,\omega,z_0,u_0']\big| \leq c(\omega)|u_0-u_0'|, \quad (a.s.).
\end{align*}

\item Under ({\bf A6}) for $L_0$, the function $L_0^*[t,\omega,z_0,u_0]$ and its first order derivative (w.r.t $z_0$) is a.s. continuous and bounded on $[0,T] \times \mathbb{R}^n \times U_0$. $L_0^*[t,\omega,z_0,u_0]$ is a.s. Lipschitz continuous in $z_0$. In addition, $L_0^*[t,\omega,0,0]$ is in the space $L^2_{\mathcal F_t}([0,T];\mathbb{R}_+)$. 

\item Under ({\bf A8}) for $H_0^{u_0}$, the set of minimizers
\begin{align*}
&\arg\inf_{u_0\in U_0} \big\{\big<f_0^*[t,\omega,z_0,u_0],p\big> + L_0^*[t,\omega,z_0,u_0]\big\},
\end{align*}
is a singleton for any $p \in \mathbb{R}^n$, and the resulting $u_0$ as a function of $[t,\omega,z_0,p]$ is a.s. continuous in $t$, a.s. Lipschitz continuous in $(z_0,p)$, uniformly with respect to $t$. In addition, $u_0[t,\omega,0,0]$ is in the space $L^2_{\mathcal F_t}([0,T];\mathbb{R}^n)$. 
\end{enumerate}
\end{proposition}

{\it Proof}: ({\it i}) We only show the results for $f_0^*$, the analysis for $\sigma_0^*$ is similar. For $\omega \in \Omega$, we take $(t,z,u)$ and $(s,z',u')$ both from $[0,T]\times\mathbb{R}^n\times U_0$. We have
\begin{align*}
& \big|f_0^*[t,\omega,z,u]-f_0^*[s,\omega,z',u']\big|\equiv\big|f_0[t,z,u,\mu_t(\omega)]-f_0[s,z',u',\mu_s(\omega)]\big| \\
& ~ \leq \big|f_0[t,z,u,\mu_t(\omega)]-f_0[s,z',u',\mu_t(\omega)]\big| + \big|f_0[s,z',u',\mu_t(\omega)]-f_0[s,z',u',\mu_s(\omega)]\big| \\
& ~ \leq \big|f_0[t,z,u,\mu_t(\omega)]-f_0[s,z,u,\mu_t(\omega)]\big|+ \big|f_0[s,z,u,\mu_t(\omega)]-f_0[s,z',u',\mu_t(\omega)]\big|  \\  
& \quad  + \big|f_0[s,z',u',\mu_t(\omega)]-f_0[s,z',u',\mu_s(\omega)].
\end{align*}
By ({\bf A4}), $f_0[t,\omega,z,u]$ is continuous with respect to $(t,z,u)$ and therefore
\begin{align*}
& \big|f_0[t,z,u,\mu_t(\omega)]-f_0[s,z,u,\mu_t(\omega)]\big|+ \big|f_0[s,z,u,\mu_t(\omega)]-f_0[s,z',u',\mu_t(\omega)]\big|  \rightarrow 0,
\end{align*}
as $|t-s|+|z-z'|+|u-u'| \rightarrow 0$. Since $\mu_{(\cdot)}(\omega)$ is in the set $\mathcal M_\rho^\beta$, $0 < \beta <1$, and by ({\bf A4}) there exists a constant $k>0$ independent of $(s,z,u)$ such that
\begin{align*}
& \big|f_0[s,z,u,y]-f_0[s,z,u,y'] \big| \leq k|y-y'|,
\end{align*}
we get $\big|f_0[s,z',u',\mu_t(\omega)]-f_0[s,z',u',\mu_s(\omega)]\rightarrow 0$ as $|t-s|\rightarrow 0$. This concludes the a.s. continuity of $f_0^*[t,\omega,z_0,u_0]$ on $[0,T] \times \mathbb{R}^n \times U_0$. 

Using the Leibniz rule we have 
\begin{align*}
& D_{z_0} f_0^*[t,\omega,z_0,u_0] = \int  D_{z_0} f_0[t,z_0,u_0,x] \mu_t(\omega)(dx),  \quad a.s.,
\end{align*}
where the partial derivative exists due to the boundedness of the first order derivative (w.r.t $z_0$) of $f_0$ by ({\bf A4}). The a.s. continuity of $D_{z_0} f_0^*$ on $[0,T] \times \mathbb{R}^n \times U_0$ may be proved by a similar argument above for $f_0^*$. Other results of the Proposition follow directly from ({\bf A4}).

({\it ii}) This is a direct result of ({\bf A5}).

({\it iii}) The proofs are similar to the proofs for $f_0^*$ in part (i).

({\it iv}) This is a direct result of ({\bf A8}) for $S_0$ using the measure $\mu_{(\cdot)}(\omega) \in \mathcal M_\rho^\beta$, $0 < \beta <1$.  \qed

Employing the results of Section \ref{S3:mm:SHJB}, we analyze the SHJB equation (\ref{CNP:SHJB:Major}) where the probability measure $\mu_{(\cdot)}(\omega)$ is in the set $\mathcal M_\rho^\beta$, $0 < \beta <1$.
\begin{theorem} 
Assume ({\bf A3})-({\bf A7}) for $U_0$, $f_0$, $\sigma_0$ and $L_0$ hold, and the probability measure $\mu_{(\cdot)}(\omega)$ is in the set $\mathcal M_\rho^\beta$, $0 < \beta <1$. Then the SHJB equation for the major agent (\ref{CNP:SHJB:Major}) has a unique solution $(\phi_0(t,x),\psi_0(t,x))$ in $\big(L_{\mathcal F_t}^2([0,T];\mathbb{R}),L_{\mathcal F_t}^2([0,T];\mathbb{R}^m)\big)$. 
\end{theorem}
{\it Proof}: Proposition \ref{CNP:Proposition:Major} indicates that the SOCP of the major agent (\ref{CNP:MajorinfDyn})-(\ref{CNP:MajorinfCos}) satisfies the Assumptions ({\bf H1})-({\bf H3}) of Section \ref{S3:mm:SHJB} with $\varsigma[t,x]=0$. The result follows directly from Theorem \ref{CNP:PenTheorem}. \qed

Let $\mu_{(\cdot)}(\omega)\in\mathcal M_\rho^\beta$, $0 < \beta <1$, be given. We assume that the unique solution $(\phi_0,\psi_0)(t,x)$ to the SHJB equation (\ref{CNP:SHJB:Major}) satisfies the regularity properties: (i) for each $t$, $(\phi_0,\psi_0)(t,x)$ is a $C^2(\mathbb{R}^n)$ map from $\mathbb{R}^n$ into $\mathbb{R} \times \mathbb{R}^m$, (ii) for each $x$, $(\phi_0,\psi_0)$ and $(D_x \phi_0,D^2_{xx} \phi_0, D_x \psi_0)$ are continuous $F_t^W$-adapted stochastic processes. Then, $\phi_0(x,t)$ coincides with the value function (\ref{CNP:ValueFunction:Major}) \cite{peng1992stochastic}, and under ({\bf A8}) for $H_0^{u_0}$ we get the best response control process (\ref{CNP:BR1:Major}):
\begin{align}
\label{CNP:Analsysi:Major:BR}
& u^o_0(t,\omega,x) \equiv u^o_0(t,x|\{\mu_{s}(\omega)\}_{0\leq s\leq T}) :=\arg\inf_{u_0\in U_0}
H_0^{u_0}[t,\omega,x,u_0,D_x\phi_0(t,\omega,x)],
\end{align}
where $(t,x)\in[0,T]\times\mathbb{R}^n$.

We introduce the following assumption (see (H6) in \cite{HuangMalhameCaines_06}).

(\textbf{A9}) For any $\mu_{(\cdot)}(\omega)\in\mathcal M_\rho^\beta$, $0 < \beta <1$, the best response control $u^o_0(t,\omega,x)$ is a.s. continuous in $(t,x)$ and a.s. Lipschitz continuous in $x$.

We denote $C_{\textrm{Lip}(x)}([0,T]\times\Omega\times\mathbb{R}^n;H)$ be the class of a.s. continuous functions from $[0,T]\times\Omega\times\mathbb{R}^n$ to $H$, which are a.s. Lipschitz continuous in $x$ \cite{HuangMalhameCaines_06}. We introduce the following well-defined map:
\begin{align}
\label{Gamma0SHJB}
&\Upsilon_0\textsuperscript{SHJB}: M_\rho^\beta \longrightarrow C_{\textrm{Lip}(x)}([0,T]\times\Omega\times\mathbb{R}^n;U_0),  \qquad 0 < \beta < 1,\\
&\Upsilon_0\textsuperscript{SHJB} \big(\mu_{(\cdot)}(\omega)\big)= u^o_0(t,\omega,x)  \equiv u^o_0(t,x|\{\mu_{s}(\omega)\}_{0\leq s\leq T}). \notag
\end{align} 

We now analyze the major agent's SMV equation (\ref{CNP:SMV:Major}) with $\mu_{(\cdot)}(\omega)\in\mathcal M_\rho^\beta$ where $0 < \beta <1$, and $u^o_0(t,\omega,x)\in C_{\textrm{Lip}(x)}([0,T]\times\Omega\times\mathbb{R}^n;U_0)$ be given in (\ref{CNP:Analsysi:Major:BR}).

\begin{theorem} \label{CNP:Analysis:Major:SMV} Assume ({\bf A3})-({\bf A7}) for $U_0$, $f_0$ and $\sigma_0$, and ({\bf A9}) hold. Let $\mu_{(\cdot)}(\omega)\in\mathcal M_\rho^\beta$ where $0 < \beta <1$, and $u^o_0(t,\omega,x)$ be given in (\ref{CNP:Analsysi:Major:BR}). Then, there exists a unique solution $z_0^o$ on $[0,T]\times\Omega$ to the major agent's SMV equation (\ref{CNP:SMV:Major}). 
\end{theorem}

{\it Proof}: Proposition \ref{CNP:Proposition:Major} indicates that the major agent's SMV equation (\ref{CNP:SMV:Major}) satisfies the Assumption ({\bf RC}) in \cite{yong1999stochastic}, page 49. The result follows directly from Theorem 6.16, Chapter 1 of \cite{yong1999stochastic}, page 49.  \qed

\begin{theorem}  \label{CNP:Analysis:Major:SMV:measure} Assume ({\bf A3})-({\bf A7}) for $U_0$, $f_0$ and $\sigma_0$, and ({\bf A9}) hold. Let $\mu_{(\cdot)}(\omega)\in\mathcal M_\rho^\beta$ where $0 < \beta <1$, and $u^o_0(t,\omega,x)$ be given in (\ref{CNP:Analsysi:Major:BR}). Then, the probability measure $\mu_{(\cdot)}^0(\omega)$ as the unit mass measure concentrated at $z_0^o(t,\omega)$ (i.e., $\mu^0_t(\omega) = \delta_{z_0^o(t,\omega)}$) which is obtained from the major agent's SMV equation (\ref{CNP:SMV:Major}) is in the class $\mathcal M_\rho^\gamma$ where $0 < \gamma <1/2$.  
\end{theorem}

{\it Proof}: We take $0 \leq s < t \leq T$. Since $\mu^0_t(\omega) = \delta_{z_0^o(t,\omega)}$, for any bounded and Lipschitz continuos function $\phi$ on $\mathbb{R}^n$ with a Lipschitz constant $K>0$, we have
\begin{align*}
& \E\big|\int_{\mathbb{R}^n} \phi(x)\mu^0_t(\omega,dx)-\int_{\mathbb{R}^n} \phi(x)\mu^0_s(\omega,dx)\big|=\E\big|\phi(z_0^o(t,\omega))-\phi(z_0^o(s,\omega))\big| \\
& \qquad  \leq K~ \E\big| z_0^o(t,\omega)- z_0^o(s,\omega)\big|.
\end{align*}    
On the other hand, Theorem \ref{CNP:Analysis:Major:SMV} indicates that there exists a unique solution to the SMV equation (\ref{CNP:SMV:Major}) such that
\begin{align*}
& z_0^o(t,\omega)- z_0^o(s,\omega)=\int_s^t  f_0[\tau,z_0^o,u^o_0,\mu_\tau(\omega)]d\tau +\int_s^t \sigma_0[\tau,z_0^o,\mu_\tau(\omega)]dw_0(\tau).
\end{align*}
Boundedness of $f_0$ and $\sigma_0$ (see ({\bf A4})), the Cauchy-Schwarz inequality and the property of It\^o integral yield
\begin{align*}
& \E\big|z_0^o(t,\omega)- z_0^o(s,\omega)\big|^2 \leq 2 C_1^2 |t-s|^2 + 2C_2^2 |t-s|,
\end{align*} 
where $C_1$ and $C_2$ are upper bounds for $f_0$ and $\sigma_0$, respectively. Hence,
\begin{align*}
& \E\big|\int_{\mathbb{R}^n} \phi(x)\mu^0_t(\omega,dx)-\int_{\mathbb{R}^n} \phi(x)\mu^0_s(\omega,dx)\big| \leq \sqrt 2 K \big(C_1|t-s|+ C_2|t-s|^{1/2}\big) \\
& \qquad \qquad \qquad \qquad \qquad \qquad \qquad \qquad \qquad ~ \leq  \sqrt 2 K (C_1 \sqrt T+ C_2)|t-s|^{1/2}.
\end{align*}

By Kolmogorov's Theorem (Theorem 18.19, Page 266, \cite{leonid2007theory}), for each $0<\gamma <1/2$, $T > 0$, and almost every $\omega \in \Omega$, there exists a constant $c(\omega,\gamma,K,T)$ such that
\begin{align*}
& \big|\int_{\mathbb{R}^n} \phi(x)\mu^0_t(\omega,dx)-\int_{\mathbb{R}^n} \phi(x)\mu^0_s(\omega,dx)\big| \leq c(\omega,\gamma,K,T) |t-s|^{\gamma},
\end{align*}
for all $0 \leq s < t \leq T$. Hence, $\mu_{(\cdot)}^0(\omega)$ is in the class $\mathcal M_\rho^\gamma$ where $0 < \gamma <1/2$. \qed

By Theorems \ref{CNP:Analysis:Major:SMV} and \ref{CNP:Analysis:Major:SMV:measure} we may now introduce the following well-defined map:
\begin{align}
\label{Gamma0SMV}
& \Upsilon_0\textsuperscript{SMV}: M_\rho^\beta \times C_{Lip(x)}([0,T]\times\Omega\times\mathbb{R}^n;U_0)  \longrightarrow M_\rho^\gamma, \quad 0 < \beta < 1,~ 0< \gamma <1/2,\\
& \Upsilon_0\textsuperscript{SMV}\big(\mu_{(\cdot)}(\omega),u^o_0(t,\omega,x)\big)= \mu^0_{(\cdot)}(\omega) \equiv  \delta_{z_0^o(t,\omega)}. \notag
\end{align}

\subsection{Analysis of the Minor Agents' SMFG System}\label{CNP:MinorAgentAnalysis} Let $\mu_{(\cdot)}(\omega)\in \mathcal M_\rho^\beta$, $0 < \beta <1$, be the fixed stochastic measure assumed in Section \ref{CNP:MajorAgentAnalysis}. In this section we assume that $\mu_{(\cdot)}^0(\omega)\in \mathcal M_\rho^\gamma$, $0 < \gamma <1/2$, is the unit mass random measure concentrated at $z_0^o(\cdot,\omega)$ (i.e., $\mu^0_t(\omega) = \delta_{z_0^o(t,\omega)}$) obtained from the composite map:
\begin{align}
\label{Gamma0}
& \Upsilon_0: M_\rho^\beta \longrightarrow M_\rho^\gamma, \qquad 0 < \beta < 1, ~ 0 < \gamma < 1/2, \\
& \Upsilon_0 \big( \mu_{(\cdot)}(\omega)\big):= \Upsilon_0\textsuperscript{SMV} \Big(\mu_{(\cdot)}(\omega),\Upsilon_0\textsuperscript{SHJB}\big( \mu_{(\cdot)}(\omega)\big)\Big)= \mu^0_{(\cdot)}(\omega) \equiv  \delta_{z_0^o(t,\omega)},  \notag
\end{align}
where $\Upsilon_0\textsuperscript{SHJB}$ and $\Upsilon_0\textsuperscript{SMV}$ are given in (\ref{Gamma0SHJB}) and (\ref{Gamma0SMV}), respectively. 

Following arguments exactly parallel to those used in Section \ref{CNP:MajorAgentAnalysis}, we analyze the SHJB equation (\ref{CNP:SHJB:Minor}) where the probability measures $\mu_{(\cdot)}(\omega)\in\mathcal M_\rho^\beta$, $0 < \beta <1$ and $\mu^0_{(\cdot)}(\omega) \in \mathcal M_\rho^\gamma$, $0 < \gamma <1/2$.
\begin{theorem} 
Assume ({\bf A3})-({\bf A7}) for $U$, $f$, $\sigma$ and $L$ hold, and $\mu_{(\cdot)}(\omega)\in \mathcal M_\rho^\beta$, $0 < \beta <1$ and $\mu^0_{(\cdot)}(\omega)$ is in the set $\mathcal M_\rho^\gamma$, $0 < \gamma <1/2$. Then the SHJB equation for the generic minor agent (\ref{CNP:SHJB1:Minor}) has a unique solution $(\phi_i(t,x),\psi_i(t,x))$ in $\big(L_{\mathcal F_t}^2([0,T];\mathbb{R}),L_{\mathcal F_t}^2([0,T];\mathbb{R}^m)\big)$.
\end{theorem}

{\it Proof}: A similar argument to Proposition \ref{CNP:Proposition:Major} for the generic minor agent (see Proposition \ref{CNP:Proposition:Minor} in \cite{MNourian_SiamAppendix}) indicates that the SOCP of the generic minor agent (\ref{CNP:MinorinfDyn})-(\ref{CNP:MinorinfCos}) satisfies the Assumptions ({\bf H1})-({\bf H3}) of Section \ref{S3:mm:SHJB} with $\sigma[t,x]=0$. The result follows directly from Theorem \ref{CNP:PenTheorem}. \qed

For the probability measure $\mu_{(\cdot)}(\omega)\in\mathcal M_\rho^\beta$, $0 < \beta <1$, and $\mu^0_{(\cdot)}(\omega)\in\mathcal M_\rho^\gamma$, $0 < \gamma <1/2$, we assume that the unique solution $(\phi_i,\psi_i)(t,x)$ to the SHJB equation (\ref{CNP:SHJB1:Minor}) satisfies the regularity properties: (i) for each $t$, $(\phi_i,\psi_i)(t,x)$ is a $C^2(\mathbb{R}^n)$ map from $\mathbb{R}^n$ into $\mathbb{R} \times \mathbb{R}^m$, (ii) for each $x$, $(\phi_i,\psi_i)$ and $(D_x \phi_i,D^2_{xx} \phi_i, D_x \psi_i)$ are continuous $F_t^W$-adapted stochastic processes. Then, $\phi_i(x,t)$ coincides with the value function (\ref{CNP:ValueFunction:Minor}) \cite{peng1992stochastic}, and under ({\bf A8}) for $H^u$ we get the best response control process (\ref{CNP:BR1:Minor}):
\begin{align}
\label{CNP:Analsysi:Minor:BR}
& u^o_i(t,\omega,x) \equiv u^o_i(t,x|\{\mu^0_{s}(\omega),\mu_{s}(\omega)\}_{0\leq s\leq T}) \\
& ~ \qquad \qquad  :=\arg\inf_{u_i\in U} H^{u}[t,\omega,x,u_i,D_x\phi_i(t,\omega,x)], \notag
\end{align}
where $(t,x)\in[0,T]\times\mathbb{R}^n$.

We introduce the following assumption (see (\textbf{A9}) or (H6) in \cite{HuangMalhameCaines_06}).

(\textbf{A10}) For any $\mu_{(\cdot)}(\omega)\in\mathcal M_\rho^\beta$, $0 < \beta <1$, and $\mu^0_{(\cdot)}(\omega)\in\mathcal M_\rho^\gamma$, $0 < \gamma <1/2$, the best response control process $u^o_i(t,\omega,x)$ is a.s. continuous in $(t,x)$ and a.s. Lipschitz continuous in $x$.

We introduce the following well-defined map for the generic minor agent $i$:
\begin{align}
\label{GammaiSHJB}
& \!\!  \Upsilon_i\textsuperscript{SHJB}\!\!: M_\rho^\beta \times M_\rho^\gamma  \longrightarrow C_{\textrm{Lip}(x)}([0,T]\times\Omega\times\mathbb{R}^n;U), \quad 0 < \beta < 1, ~0 < \gamma <1/2,  \\
& \Upsilon_i\textsuperscript{SHJB} \big(\mu_{(\cdot)}(\omega),\mu^0_{(\cdot)}(\omega)\big)= u^o_i(t,\omega,x)  \equiv u^o_i(t,x|\{\mu^0_{s}(\omega),\mu_{s}(\omega)\}_{0\leq s\leq T}).  \notag
\end{align} 

For given probability measure $\mu^0_{(\cdot)}(\omega)\in\mathcal M_\rho^\gamma$, $0 < \gamma <1/2$, we analyze the generic minor agent's SMV equation (\ref{CNP:Minor:CL}):
\begin{align}
\label{CNP:Minor:CL2}
& ~ dz_i^o(t,\omega,\omega') = f[t,z_i^o,u^o_i(t,\omega,z_i^o),\mu_t^0(\omega),\mu_t(\omega)]dt \\
& \qquad \qquad \qquad + \sigma[t,z_i^o,\mu_t^0(\omega),\mu_t(\omega)]dw_i(t,\omega'),
\quad z_i^o(0)=z_i(0), \notag 
\end{align}
 where $u_i^o(t,\omega,x)\in C_{\textrm{Lip}(x)}([0,T]\times\Omega\times\mathbb{R}^n;U)$ is given in (\ref{CNP:Analsysi:Minor:BR}). We call the pair $\big(z_i^o(\cdot,\omega,\omega'),\mu_{(\cdot)}(\omega)\big)$ a consistent solution of the generic minor agent's SMV equation (\ref{CNP:Minor:CL2}) if $\big(z_i^o(\cdot,\omega,\omega'),\mu_{(\cdot)}(\omega)\big)$ solves (\ref{CNP:Minor:CL2}) and $\mu_{(\cdot)}(\omega)$ be the the law of the process $z_i^o(\cdot,\omega,\omega')$, i.e., $\mu_{(\cdot)}=\mathcal L\big(z_i^o(\cdot,\omega,\omega')\big)$. We define $\Lambda$ as the map which associates to $\mu_{(\cdot)}(\omega) \in \mathcal M_\rho^\beta$, $0 < \beta <1/2$, the law of the process $z_i^o(\cdot,\omega,\omega')$ in (\ref{CNP:Minor:CL2}):
\begin{align}
\label{CNP:FixedPoint:minoragent:Lambda}
& z_i^o(t,\omega,\omega') = z_i^o(0) + \int_0^t \Big(\int_{\mathbb{R}^n}  \int_{\mathbb{R}^n} f[s,z_i^o,u^o_i,y,z] d\mu_s^0(\omega)(y)  d\mu_s(\omega)(z)\Big) ds \\
& \qquad \qquad \qquad + \int_0^t \Big(\int_{\mathbb{R}^n}  \int_{\mathbb{R}^n} \sigma[s,z_i^o,y,z] d\mu_s^0(\omega)(y) d\mu_s(\omega)(z)\Big) dw_i(s,\omega'), \notag
\end{align}
where we observe that the law $\Lambda$ depends on the sample point $\omega \in \Omega$.

We now show that there exists a unique $\mu_{(\cdot)}(\omega) \in \mathcal M_\rho^\beta$, $0 < \beta <1$, such that $\mu(\omega)=\Lambda\big( \mu(\omega)\big)$. The proof of the following theorem, which is given in Appendix D in \cite{MNourian_SiamAppendix}, is based upon a fixed point argument with random parameters (see Theorem 6 in \cite{HuangMalhameCaines_06} and Theorem 1.1 in \cite{sznitman1991topics} for the standard fixed point argument). 

\begin{theorem} \label{CNP:Analysis:Minor:SMV} Assume ({\bf A3})-({\bf A7}) for $U$, $f$ and $\sigma$, and ({\bf A10}) hold. Let $\mu^0_{(\cdot)}(\omega)$ be in the set $\mathcal M_\rho^\gamma$ where $0 < \gamma <1/2$, and $u_i^o(t,\omega,x)$ be given in (\ref{CNP:Analsysi:Minor:BR}). Then, there exists a unique consistent solution pair $\big(z_i^o(\cdot,\omega,\omega'),\mu_{(\cdot)}(\omega)\big)$ to the generic minor agent's SMV equation (\ref{CNP:Minor:CL2}) where $\mu_{(\cdot)}(\omega)=\mathcal L\big(z_i^o(\cdot,\omega,\omega')\big)$. \qed
\end{theorem}

\begin{theorem}\label{CNP:Analysis:Minor:SMV:measure} Assume ({\bf A3})-({\bf A7}) for $U$, $f$ and $\sigma$, and ({\bf A10}) hold. Let $\mu^0_{(\cdot)}(\omega)$ be in the set $\mathcal M_\rho^\gamma$ where $0 < \gamma <1/2$. For given $u_i^o(t,\omega,x)$ in (\ref{CNP:Analsysi:Minor:BR}), let $\big(z_i^o(\cdot,\omega,\omega'),\mu_{(\cdot)}(\omega)\big)$ be the consistent solution pair of the SMV equation (\ref{CNP:Minor:CL2}). Then, the probability measure $\mu_{(\cdot)}(\omega)$ is in the class $\mathcal M_\rho^\beta$ where $0 < \beta <1$.  
\end{theorem}

{\it Proof}: We take $0 \leq s < t \leq T$. For any bounded and Lipschitz continuos function $\phi$ on $\mathbb{R}^n$ with a Lipschitz constant $K>0$, we have
\begin{align*}
& \E\big|\int_{\mathbb{R}^n} \phi(x)\mu_t(\omega,dx)-\int_{\mathbb{R}^n} \phi(x)\mu_s(\omega,dx)\big|=\E\big|\E_{\omega}\big(\phi(z_i^o(t,\omega,\omega'))-\phi(z_i^o(s,\omega,\omega'))\big)\big| \\
& \qquad  \leq K~ \E\big|\E_{\omega}\big(z_i^o(t,\omega,\omega')- z_i^o(s,\omega,\omega')\big)\big|.
\end{align*}    
On the other hand, Theorem \ref{CNP:Analysis:Minor:SMV} indicates that there exists a unique solution to the SMV equation (\ref{CNP:Minor:CL2}) such that
\begin{align*}
& \E_{\omega}\big(z_i^o(t,\omega,\omega')-z_i^o(s,\omega,\omega')\big)=\int_s^t  f[\tau,z_i^o,u^i_0,\mu^0_\tau(\omega),\mu_\tau(\omega)]d\tau,
\end{align*}
where we note that $\E_{\omega} \int_0^t \sigma[\tau,z_i^o,\mu^0_\tau(\omega),\mu_\tau(\omega)]dw_i(\tau,\omega')=0$ for $0 \leq t \leq T$. Boundedness of $f$ (see ({\bf A4})) yields
\begin{align*}
& \E\big|\E_{\omega}\big(z_i^o(t,\omega,\omega')-z_i^o(s,\omega,\omega')\big)\big| \leq C_1 |t-s|,
\end{align*} 
where $C_1$ is the upper bound for $f$.

By Kolmogorov's Theorem (Theorem 18.19, \cite{leonid2007theory}, Page 266), for each $0<\gamma <1$, $T > 0$, and almost every $\omega \in \Omega$, there exists a constant $c(\omega,\gamma,K,T)$ such that     
\begin{align*}
& \big|\int_{\mathbb{R}^n} \phi(x)\mu_t(\omega,dx)-\int_{\mathbb{R}^n} \phi(x)\mu_s(\omega,dx)\big| \leq c(\omega,\gamma,K,T) |t-s|^{\gamma},
\end{align*}
for all $0 \leq s < t \leq T$. Hence, $\mu_{(\cdot)}(\omega)$ is in the class $\mathcal M_\rho^\beta$ where $0 < \beta <1$. \qed

By Theorems \ref{CNP:Analysis:Minor:SMV} and \ref{CNP:Analysis:Minor:SMV:measure} we may now introduce the following well-defined map:
\begin{align}
 \label{GammaiSMV}
& \Upsilon_i\textsuperscript{SMV}:  M_\rho^\beta \times M_\rho^\gamma \times  C_{\textrm{Lip}(x)}([0,T]\times\Omega\times\mathbb{R}^n;U_0)  \longrightarrow M_\rho^\beta, ~ 0 < \beta <1, ~ 0 < \gamma < 1/2,\\
& \Upsilon_i\textsuperscript{SMV} \big( \mu_{(\cdot)}(\omega),\mu^0_{(\cdot)}(\omega),u^o_i(t,\omega,x)\big)= \mu_{(\cdot)}(\omega). \notag
\end{align}

\subsection{Analysis of the Joint Major and Minor Agents' SMFG System} \label{analysisjoint} Based on the analysis of Sections \ref{CNP:MajorAgentAnalysis} and \ref{CNP:MinorAgentAnalysis} we obtain the following well-defined map:
\begin{align}
\label{Gamma}
& \Upsilon: M_\rho^\beta \longrightarrow M_\rho^\beta, \qquad 0 < \beta < 1, \\
& \Upsilon \big(\mu_{(\cdot)}(\omega)\big)= \Upsilon_i\textsuperscript{SMV} \Big( \mu_{(\cdot)}(\omega),\Upsilon_0\big(\mu_{(\cdot)}(\omega)\big),\Upsilon_i\textsuperscript{SHJB} \big(\mu_{(\cdot)}(\omega)),\Upsilon_0\big(\mu_{(\cdot)}(\omega)\big)\big)\Big)=\mu_{(\cdot)}(\omega),  \notag
\end{align}
which is the composition of the maps $\Upsilon_0$, $\Upsilon_i^{\textrm{SHJB}}$ and $ \Upsilon^{\textrm{SMV}}_i$ introduced in (\ref{Gamma0}), (\ref{GammaiSHJB}) and (\ref{GammaiSMV}), respectively. Subsequently, the problem of existence and uniqueness of solution to the MM SMV system (\ref{CNP:SHJB:Major})-(\ref{CNP:SMV:Major}) and (\ref{CNP:SHJB:Minor})-(\ref{CNP:SMV:Minor}) is translated into a fixed point problem with random parameters for the map $\Upsilon$ on the Polish space $\mathcal M_\rho^\beta$, $0 < \beta < 1$. 

We introduce the following assumption without which one needs to work with the ``expectation'' of the Wasserstein metric $D_{(\cdot)}^\rho$ of stochastic measure.

(\textbf{A11}) We assume that the diffusion coefficient of the major agent $\sigma_0$ in (\ref{CNP:Major:GenDyn}) does not depend on its own state $z_0^N$ and the states of the minor agents $z_i^N$, $1 \leq i \leq N$. 

The proof of the following lemma is given in Appendix E in \cite{MNourian_SiamAppendix}.

\begin{lemma} \label{CNP:Lemma:Sensetivity} (i) Assume ({\bf A3})-({\bf A7}) for $U_0$, $f_0$ and $\sigma_0$, and ({\bf A11}) hold. Let $\mu_{(\cdot)}(\omega)$ be in the set $\mathcal M_\rho^\beta$ where $0 < \beta <1$. Then, for given $u_0, u_0' \in C_{\textrm{Lip}(x)}([0,T]\times\Omega\times\mathbb{R}^n;U_0)$ there exists a constant $c_0$ such that 
\begin{align} \label{CNP:ContractionLemmai}
& \Big( D_T^\rho \big(\mu^0(\omega),\nu^0(\omega)\big) \Big)^2 \leq c_0 \sup_{(t,x)\in [0,T]\times\mathbb{R}^n} \big| u_0(t,\omega,x) -u_0'(t,\omega,x)\big|^2, \qquad a.s.,
\end{align}
where $\mu^0(\omega),\nu^0(\omega)\in\mathcal M^\gamma_\rho$, $0 < \gamma <1/2$, are induced by the map $\Upsilon^{\textrm{SMV}}_0$ in (\ref{Gamma0SMV}) using the two control processes $u_0$ and $u_0'$, respectively.\\
(ii) Assume ({\bf A3})-({\bf A7}) for $U_0$, $f_0$ and $\sigma_0$, and ({\bf A11}) hold. Let $u_0^o$ be in the space $C_{\textrm{Lip}(x)}([0,T]\times\Omega\times\mathbb{R}^n;U_0)$. Then, for given $\mu(\omega),\nu(\omega) \in \mathcal M_\rho^\beta$, $0 < \beta <1$, there exists a constant $c_1$ such that 
\begin{align}
\label{CNP:ContractionLemmaii}
& \Big(D_T^\rho \big(\mu^0(\omega),\nu^0(\omega)\big) \Big)^2 \leq c_1\Big(D_T^\rho\big(\mu(\omega),\nu(\omega)\big)\Big)^2, \qquad a.s.,
\end{align}
where $\mu^0(\omega),\nu^0(\omega)\in\mathcal M^\gamma_\rho$, $0 < \gamma <1/2$, are induced by the map $\Upsilon^{\textrm{SMV}}_0$ in (\ref{Gamma0SMV}) using the stochastic measures $\mu(\omega)$ and $\nu(\omega)$, respectively.\\
(iii) Assume ({\bf A3})-({\bf A7}) for $U$, $f$ and $\sigma$ hold. Let $\mu^0_{(\cdot)}(\omega)$ be in the set $\mathcal M_\rho^\gamma$ where $0 < \gamma <1/2$. Then, for given $u, u' \in C_{\textrm{Lip}(x)}([0,T]\times\Omega\times\mathbb{R}^n;U)$ there exists a constant $c_2$ such that 
\begin{align} \label{CNP:ContractionLemmaiii}
&\Big(D_T^\rho \big(\mu(\omega),\nu(\omega)\big)\Big)^2 \leq c_2 \sup_{(t,x)\in [0,T]\times\mathbb{R}^n} \big| u(t,\omega,x) -u'(t,\omega,x)\big|^2, \qquad a.s.,
\end{align}
where $\mu(\omega),\nu(\omega)\in\mathcal M^\beta_\rho$, $0 < \beta <1$, are induced by the map $\Upsilon^{\textrm{SMV}}_i$ in (\ref{GammaiSMV}) using the two control processes $u$ and $u'$, respectively.\\
(iv) Assume ({\bf A3})-({\bf A7}) for $U$, $f$ and $\sigma$ hold. Let $u_i^o$ be in the space $C_{\textrm{Lip}(x)}([0,T]\times\Omega\times\mathbb{R}^n;U)$. Then, for given $\mu^0(\omega),\nu^0(\omega) \in \mathcal M_\rho^\gamma$, $0 < \gamma <1/2$, there exists a constant $c_3$ such that 
\begin{align}
\label{CNP:ContractionLemmaiv} & \Big(D_T^\rho\big(\mu(\omega),\nu(\omega)\big)\Big)^2 \leq c_3 \Big(D_T^\rho\big(\mu^0(\omega),\nu^0(\omega)\big)\Big)^2, \qquad a.s.,
\end{align}
where $\mu(\omega),\nu(\omega)\in\mathcal M^\beta_\rho$, $0 < \beta <1$, are induced by the map $\Upsilon^{\textrm{SMV}}_i$ in (\ref{GammaiSMV}) using the stochastic measures $\mu^0(\omega)$ and $\nu^0(\omega)$, respectively. \qed
\end{lemma}

We define the G\^ateaux derivative of the function $F(t,x,\mu)$ with respect to the measure $\mu(y)$ as \cite{kolokoltsov2011mean}
\begin{align*}
& \partial_{\mu(y)} F(t,x,\mu) = \lim_{\epsilon \rightarrow 0} \frac{F(t,x,\mu+\epsilon \delta(y))-F(t,x,\mu)}{\epsilon},
\end{align*}
where $\delta$ is the Dirac delta function. We introduce the following assumptions:

({\bf A12}) (i) In (\ref{CNP:MajorinfDyn})-(\ref{CNP:MajorinfCos}) the G\^ateaux derivative of $f_0$, $\sigma_0$ and $L_0$ with respect to $\mu$ exist, are $C^\f (\mathbb{R}^n)$ and a.s. uniformly bounded. (ii) In (\ref{CNP:MinorinfDyn})-(\ref{CNP:MinorinfCos}) the partial derivatives of $f$, $\sigma$ and $L$ with respect to $\mu^0$ and $\mu$ exist, are $C^\f (\mathbb{R}^n)$ and a.s. uniformly bounded.

The proof of the following lemma is based on the sensitivity analysis of the SHJB equations (\ref{CNP:SHJB:Major}) and (\ref{CNP:SHJB:Minor}) to the stochastic measures $\mu_{(\cdot)}(\omega)$ and $\mu^0_{(\cdot)}(\omega)$ developed in Appendix F in \cite{MNourian_SiamAppendix} (see also Section 6 in \cite{kolokoltsov2011mean}).

\begin{lemma} \label{CNP:Lemma:Sensetivity:Control} (i) Assume ({\bf A3})-({\bf A7}) for $U_0$, $f_0$, $\sigma_0$, $L_0$, and ({\bf A12})-(i) hold. Let $(\phi_0(t,x),\psi_0(t,x))$ be the unique solution pair to (\ref{CNP:SHJB:Major}) which is $C^\f (\mathbb{R}^n)$ and is a.s. uniformly bounded. In addition, we assume ({\bf A8}) holds for $S_0$ and the resulting $u_0$ is also a.s. Lipschitz continuous in $\mu$. Then, for $\mu_{(\cdot)}(\omega)$ and $\nu_{(\cdot)}(\omega)$ in the set $\mathcal M_\rho^\beta$, $0 < \beta <1$, there exists a constant $c_4$ such that 
\begin{align} \label{CNP:Contraction2Lemmai} 
&\sup_{(t,x)\in [0,T]\times\mathbb{R}^n} \big| u_0(t,\omega,x) -u_0'(t,\omega,x)\big|^2 \leq c_4 \Big(D_T^\rho\big(\mu(\omega),\nu(\omega)\big)\Big)^2, \qquad a.s.,
\end{align}
where $u_0,u_0'\in C_{\textrm{Lip}(x)}([0,T]\times\Omega\times\mathbb{R}^n;U_0) $ are induced by the map $\Upsilon^{\textrm{SHJB}}_0$ in (\ref{Gamma0SHJB}) using two stochastic measures $\mu_{(\cdot)}(\omega)$ and $\nu_{(\cdot)}(\omega)$, respectively.\\
(ii) Assume ({\bf A3})-({\bf A7}) for $U$, $f$, $\sigma$, $L$, and ({\bf A12})-(ii) hold. Let $(\phi(t,x),\psi(t,x))$ be the unique solution pair to (\ref{CNP:SHJB:Minor}) which is $C^\f (\mathbb{R}^n)$ and is a.s. uniformly bounded. In addition, we assume ({\bf A8}) holds for $S$ and the resulting $u$ is also a.s. Lipschitz continuous in $\mu$. Then, for $\mu^0_{(\cdot)}(\omega) \in \mathcal M_\rho^\gamma$, $0 < \gamma <1/2$, and $\mu_{(\cdot)}(\omega)$ and $\nu_{(\cdot)}(\omega)$ in the set $\mathcal M_\rho^\beta$, $0 < \beta <1$, there exists a constant $c_5$ such that 
\begin{align} \label{CNP:Contraction2Lemmaii} 
&\sup_{(t,x)\in [0,T]\times\mathbb{R}^n} \big| u(t,\omega,x) -u'(t,\omega,x)\big|^2 \leq c_5 \Big(D_T^\rho\big(\mu(\omega),\nu(\omega)\big) \Big)^2, \qquad a.s.,
\end{align}
where $u,u'\in C_{\textrm{Lip}(x)}([0,T]\times\Omega\times\mathbb{R}^n;U) $ are induced by the map $\Upsilon^{\textrm{SHJB}}_i$ in (\ref{GammaiSHJB}) using two stochastic measures $\mu_{(\cdot)}(\omega)$ and $\nu_{(\cdot)}(\omega)$, respectively.\\
(iii) Assume ({\bf A3})-({\bf A7}) for $U$, $f$, $\sigma$, $L$, and ({\bf A12})-(ii) hold. Let $(\phi(t,x),\psi(t,x))$ be the unique solution pair to (\ref{CNP:SHJB:Minor}) which is $C^\f (\mathbb{R}^n)$ and is a.s. uniformly bounded. In addition, we assume ({\bf A8}) holds for $S$ and the resulting $u$ is also a.s. Lipschitz continuous in $\mu^0$. Then, for $\mu_{(\cdot)}(\omega) \in \mathcal M_\rho^\beta$, $0 < \beta <1$, and $\mu^0_{(\cdot)}(\omega)$ and $\nu^0_{(\cdot)}(\omega)$ in the set $\mathcal M_\rho^\gamma$, $0 < \gamma <1/2$, there exists a constant $c_6$ such that 
\begin{align} \label{CNP:Contraction2Lemmaiii} 
&\sup_{(t,x)\in [0,T]\times\mathbb{R}^n} \big| u(t,\omega,x) -u'(t,\omega,x)\big|^2 \leq c_6 \Big(D^\rho_T \big(\mu^0(\omega),\nu^0(\omega)\big)\Big)^2, \qquad a.s.,
\end{align}
where $u,u'\in C_{\textrm{Lip}(x)}([0,T]\times\Omega\times\mathbb{R}^n;U) $ are induced by the map $\Upsilon^{\textrm{SHJB}}_i$ in (\ref{GammaiSHJB}) using the two stochastic measures $\mu^0_{(\cdot)}(\omega)$ and $\nu^0_{(\cdot)}(\omega)$, respectively.
\end{lemma}

{\it Proof}: (i) Assumption ({\bf A8}) for $S_0$ together with the fact that the resulting $u_0$ in ({\bf A8}) is also a.s. Lipschitz continuous in $\mu$ yields
\begin{align}
& |u_0(t,\omega,x)-u_0'(t,\omega,x)| \leq k_1 D_t^\rho\big(\mu(\omega),\nu(\omega)\big) \notag \\
& \hspace{5cm} + k_2 |D_x \phi_0^\mu(t,\omega,x)-D_x \phi_0^\nu(t,\omega,x)|,  \label{SAT:parti:1}
\end{align}
with positive constants $k_1, k_2$, where we indicate the dependence of $\phi_0$ on measures $\mu$ and $\nu$ by $\phi_0^\mu$ and $\phi_0^\nu$, respectively.

We consider the G\^ateaux derivative of $\phi_0$ with respect to the measure $\mu$. The assumptions of the theorem imply that the conditions for Proposition \ref{Proposition;alphadependent} in \cite{MNourian_SiamAppendix} hold. Therefore, Proposition \ref{Proposition;alphadependent} in \cite{MNourian_SiamAppendix} concludes that the G\^ateaux derivative of $D_x \phi_0$ with respect to measure $\mu$ is a.s. uniformly bounded. This together with the mean value theorem yields
\begin{align}  \label{SAT:parti:2}
&|D_x \phi_0^\mu(t,\omega,x)-D_x \phi_0^\nu(t,\omega,x)| \leq k_3 D_t^\rho\big(\mu(\omega),\nu(\omega)\big),
\end{align}
with positive constant $k_3$. (\ref{SAT:parti:1}) and (\ref{SAT:parti:2}) give
\begin{align*} 
& |u_0(t,\omega,x)-u_0'(t,\omega,x)| \leq k D_t^\rho\big(\mu(\omega),\nu(\omega)\big),
\end{align*}
with $k:=k_1+k_2k_3$, which yields the result. \qed

\begin{remark} In the standard mean field game model of \cite{HuangMalhameCaines_06} a similar condition to (\ref{CNP:Contraction2Lemmai})-(\ref{CNP:Contraction2Lemmaiii}) is taken as an assumption (see the feedback regularity condition (37) in \cite{HuangMalhameCaines_06}). Following the argument in Section 7.1 of \cite{HuangMalhameCaines_06}, one can show that the inequalities (\ref{CNP:Contraction2Lemmai})-(\ref{CNP:Contraction2Lemmaiii}) hold in the linear-quadratic-Gaussian (LQG) model with Lipschitz continuous nonlinear couplings.
\end{remark}

We recall the map $\Upsilon$ given in (\ref{Gamma}) which is the composition of the maps $\Upsilon_0$, $\Upsilon_i^{\textrm{SHJB}}$ and $ \Upsilon^{\textrm{SMV}}_i$ introduced in (\ref{Gamma0}), (\ref{GammaiSHJB}) and (\ref{GammaiSMV}), respectively (see the diagram below).
\begin{align*}
& \begin{array}[c]{ccc}
\mu_{(\cdot)}(\omega) &\stackrel{\Upsilon_0^{\textrm{SHJB}}}{\longrightarrow}& u_0^o(\cdot,\omega,x) \\
\uparrow\scriptstyle{ \Upsilon^{\textrm{SMV}}_i}&&\downarrow\scriptstyle{ \Upsilon^{\textrm{SMV}}_0}\\
u^o(\cdot,\omega,x)  &\stackrel{\Upsilon_i^{\textrm{SHJB}}}{\longleftarrow}& \mu_{(\cdot)}^0(\omega) \equiv \delta_{z_0^o(t,\omega)}
\end{array}
\end{align*}

\begin{theorem} \label{mm:maintheorem} ({\bf Main Result}) Let the assumptions of both Lemma \ref{CNP:Lemma:Sensetivity} and Lemma \ref{CNP:Lemma:Sensetivity:Control} hold. If the constants $\{c_i: 0 \leq i \leq 6\}$ for (\ref{CNP:ContractionLemmai})-(\ref{CNP:ContractionLemmaiv}) and (\ref{CNP:Contraction2Lemmai})-(\ref{CNP:Contraction2Lemmaiii}) satisfy the gain condition
\begin{align*}
& \max  \left\{c_2 c_5, c_2 c_6 c_0, c_2 c_6 c_1, c_3 c_1, c_3 c_0 c_4\right\} <1,
\end{align*}
then there exists a unique solution for the map $\Upsilon$, and hence a unique solution to the MM-SMFG system (\ref{CNP:SHJB:Major})-(\ref{CNP:SMV:Major}) and (\ref{CNP:SHJB:Minor})-(\ref{CNP:SMV:Minor}). 
\end{theorem}

{\it Proof}: The result follows from the Banach fixed point theorem for the map $\Upsilon$ given in (\ref{Gamma}) on the Polish space $\mathcal M_\rho^\beta$, $0 < \beta <1$. We note that the gain condition ensures that $\Upsilon$ is a contraction.  \qed

As in the classical FBSDEs, the gain condition in Theorem \ref{mm:maintheorem} is expected to hold for short time-horizon $T$. Another approach to the solution existence of the MM-SMFG system (\ref{CNP:SHJB:Major})-(\ref{CNP:SMV:Major}) and (\ref{CNP:SHJB:Minor})-(\ref{CNP:SMV:Minor}) is Schauder's fixed point argument which is the topic of future work. 

\section{$\epsilon$-Nash Equilibrium Property of the SMFG Control Laws} \label{S8:mm:Nash} We let
\begin{align*}
&
\big(\phi_0(t,\omega,x),\psi_0(t,\omega,x),u^o_0(t,\omega,x),z_0^o(t,\omega),\phi(t
,\omega,x),\psi(t,\omega,x),u^o(t,\omega,x),z^o(t,\omega)\big),
\end{align*}
be the unique solution of the MM-SMFG system (\ref{CNP:SHJB:Major})-(\ref{CNP:SMV:Major}) and (\ref{CNP:SHJB:Minor})-(\ref{CNP:SMV:Minor}) such that SMFG best response $u^o_0(t,\omega,x)$ and $u^o(t,\omega,x)$ are a.s. continuous in $(t,x)$ and a.s. Lipschitz continuous in $x$.

We now apply the SMFG best responses $u^o_0(t,\omega,x)$ and $u^o(t,\omega,x)$ into a finite $N+1$ major and minor population (\ref{CNP:Major:GenDyn})-(\ref{CNP:Minor:GenDyn}). This yields the following closed loop individual dynamics: 
\begin{align}
\label{CNP:Major:GenDyn:closed}
& dz_0^{o,N}(t) =\frac{1}{N} \sum_{j=1}^N f_0[t,z_0^{o,N}(t),u_0^o(t,z_0^{o,N}(t)),z_j^{o,N}(t)]dt \\
& \quad + \frac{1}{N} \sum_{j=1}^N \sigma_0[t,z_0^{o,N}(t),z_j^{o,N}(t)] dw_0(t),\quad z_0^{o,N}(0)=
z_0(0), ~ 0 \leq t \leq T, \notag \\
\label{CNP:Minor:GenDyn:closed}
& dz_i^{o,N}(t) = \frac{1}{N} \sum_{j=1}^N f[t,z_i^{o,N}(t),u^o(t,z_i^{o,N}(t)),z_0^{o,N}(t),z_j^{o,N}(t)] dt \\
& \quad +\frac{1}{N} \sum_{j=1}^N \sigma [t,z_i^{o,N}(t),z_0^{o,N}(t),z_j^{o,N}(t)] dw_i(t),\quad z_i^{o,N}(0)=
z_i(0), ~ 1 \leq i \leq N, \notag
\end{align} 

We set the admissible control set of agent $\mathcal A_j$, $0 \leq j \leq N$, as 
\begin{align*}
& \mathcal{U}_j = \Big\{u_j(\cdot,\omega):= u_j\big(\cdot,\omega, z_0(\cdot,\omega),\cdots,z_N(\cdot,\omega)\big) \in C_{\textrm{Lip}(z_0, \cdots, z_N)} : u_j (t,\omega)~
\textrm{is a}\\
& \quad \textrm{$\mathcal F_t^{w_0}$-measurable process adapted to sigma-field}~\sigma\big\{z_i(\tau,\omega): 0 \leq i \leq N, 0 \leq \tau \leq t\big\} \\
& \quad \textrm{such that}~ \E\int_0^T |u_j(t,\omega)|^2 dt <\infty\Big\}.
\end{align*}
We note that $\mathcal{U}_j$, $0 \leq j \leq N$, are the full information admissible control which are not restricted to be decentralized. 

\begin{definition}
Given $\epsilon > 0$, the admissible control laws $(u_{0}^o, \cdots,u_N^o)$ for $N+1$ agents
generates an $\epsilon$-Nash equilibrium with respect to the costs $J_{j}^N$, $0 \leq j \leq
N$, if $J_j^{N}(u^o_j;u^o_{-j})-\epsilon \leq \inf_{u_j \in \mathcal U_j} J_j^{N}(u_j;u^o_{-j}) \leq
J_j^{N}(u^o_j;u^o_{-j}),$ for any $0 \leq j \leq N$. \qed
\end{definition}

We now show that the SMFG best responses for a finite $N+1$ major and minor population system (\ref{CNP:Major:GenDyn:closed})-(\ref{CNP:Minor:GenDyn:closed}) is an $\epsilon$-Nash equilibrium with respect to the cost functions (\ref{CNP:Major:GenCost})-(\ref{CNP:Minor:GenCost}) in the case that minor agents are coupled to the major agent only through their cost functions (see the MM-MFG LQG model in \cite{NguyenHuang_CDC11}).

({\bf A13}) Assume the functions $f$ and $\sigma$ in (\ref{CNP:Minor:GenDyn}) (and hence in (\ref{CNP:Minor:GenDyn:closed})) do not contain the state of major agent $z_0^N$. 

Note that in the case of assumption ({\bf A13}) the major agent $\mathcal A_0$ has non-negligible influence on the minor agents through their cost functions (\ref{CNP:Minor:GenCost}). An analysis based on the anticipative variational calculations used in the MM-MFG LQG case \cite{NguyenHuangCDC2012} is required for establishing the $\epsilon$-Nash equilibrium property of the SMFG best responses in the general case. This extension is currently under investigation and will be reported in future work.

\begin{theorem} Assume ({\bf A1})-({\bf A6}) and ({\bf A13}) hold, and there exists a unique solution to the MM-SMFG system (\ref{CNP:SHJB:Major})-(\ref{CNP:SMV:Major}) and (\ref{CNP:SHJB:Minor})-(\ref{CNP:SMV:Minor}) such that the SMFG best response control processes $u^o_0(t,\omega,x)$ and $u^o(t,\omega,x)$ are a.s. continuous in $(t,x)$ and a.s. Lipschitz continuous in $x$. Then $(u^o_0, u^o_1, \cdots, u^o_N)$ where $u^o_i \equiv u^o$, $1 \leq i \leq N$, generates an $O(\epsilon_N+1/\sqrt N)$-Nash equilibrium with respect to the cost functions (\ref{CNP:Major:GenCost})-(\ref{CNP:Minor:GenCost}) such that $\lim_{N \rightarrow \infty} \epsilon_N=0$.
\end{theorem}

{\it Proof}: Under ({\bf A13}) we have the the following closed loop individual dynamics under the SMFG best response control processes: 
\begin{align*}
& dz_0^{o,N}(t) =\frac{1}{N} \sum_{j=1}^N f_0[t,z_0^{o,N}(t),u_0^o(t,z_0^{o,N}(t)),z_j^{o,N}(t)]dt \\
& \quad + \frac{1}{N} \sum_{j=1}^N \sigma_0[t,z_0^{o,N}(t),z_j^{o,N}(t)] dw_0(t),\quad z_0^{o,N}(0)=
z_0(0), ~ 0 \leq t \leq T, \notag \\
& dz_i^{o,N}(t) = \frac{1}{N} \sum_{j=1}^N f[t,z_i^{o,N}(t),u^o(t,z_i^{o,N}(t)),z_j^{o,N}(t)] dt \\
& \quad +\frac{1}{N} \sum_{j=1}^N \sigma [t,z_i^{o,N}(t),z_j^{o,N}(t)] dw_i(t),\quad z_i^{o,N}(0)=
z_i(0), ~ 1 \leq i \leq N. \notag
\end{align*} 
We also introduce the associated Mckean-Vlasov (MV) system
\begin{align} \label{SMV:ENash}
& d z_0^o(t) = f_0[t,z_0^o(t),u_0^o(t,z_0^o),\mu_t] dt + \sigma_0[t,z_0^o(t),\mu_t]dw_0(t), \\
& dz_i^o(t) = f [t,z_i^o(t),u^o(t,z_i^o),\mu_t] dt + \sigma[t,z_i^o,\mu_t]
dw_i(t), \notag
\end{align}
with the initial condition $z^o_j(0)=z_j(0)$, $0 \leq j \leq N$. In the above MV equation $\mu_{t}$, $0 \leq t \leq T$, is the conditional law of $z_i^o(t)$, $1 \leq i \leq N$, given $\mathcal F_t^{w_0}$ (i.e., $\mu_t := \mathcal L \big(z_i^o(t)|\mathcal F_t^{w_0}\big)$, $1 \leq i \leq N$). Theorem \ref{Theorem:MCT} implies that
\begin{align}
& \sup_{0 \leq j \leq N} \sup_{0 \leq t \leq T}  \E |z_j^{o,N}(t)-z_j^o(t)|= O(1/\sqrt N),
\label{MM:CDC:Nash:Major:I1}
\end{align}
where the right hand side may depend upon the terminal time $T$. 

Let $z(0)=\int_{\mathbb{R}^n} x dF (x)$ be the mean value of the minor agents' initial states (see ({\bf A2})). We denote
\begin{align*}
& (\epsilon_N)^2 = \Big|\int_{\mathbb{R}^N} x^T xdF_N(x) - 2 z^T(0) \int_{\mathbb{R}^N} x dF_N(x) + z^T(0)z(0)\Big|.
\end{align*}
It is evident from ({\bf A2}) that $\lim_{N \rightarrow \infty} \epsilon_N=0$. To prove the $\epsilon$-Nash equilibrium property we consider two cases as follows.

Case I (strategy change for the major agent $\mathcal A_0$): While the minor agents are using the SMFG best response control law $u^0(t,\omega,x)$, a strategy change from $u_0^0(t,\omega,x)$ to the $\mathcal F_t^{w_0}$-adapted process $u_0\big(t,\omega,x,z_{-0}^{o,N}(t,\omega)\big) \in \mathcal U_0$ for the major agent yields 
\begin{align*}
& dz_0^{N}(t) =\frac{1}{N} \sum_{j=1}^N f_0[t,z_0^{N}(t),u_0\big(t,z_0^{N}(t),z_{-0}^{o,N}(t)\big),z_j^{o,N}(t)]dt \\
& \quad + \frac{1}{N} \sum_{j=1}^N \sigma_0[t,z_0^{N}(t),z_j^{o,N}(t)] dw_0(t),\quad z_0^{N}(0)=
z_0(0), ~ 0 \leq t \leq T,
\end{align*}
where $z_{-0}^{o,N}\equiv(z_1^{o,N}, \cdots, z_N^{o,N})$. Since minor agents are coupled to the major agent only through their cost functions (see ({\bf A13})) the strategy change of the major agent does not affect the the minor agents' states $z_i^{o,N}$ and $z_i^o$, $1 \leq i \leq N$, above.
 
Let $\hat z_0^N(\cdot)$ be the solution to
\begin{align*}
& d\hat z_0^{N}(t) =\frac{1}{N} \sum_{j=1}^N f_0[t,\hat z_0^{N}(t),u_0\big(t,\hat z_0^{N}(t),z_{-0}^{o}(t)\big),z_j^{o}(t)]dt \\
& \quad + \frac{1}{N} \sum_{j=1}^N \sigma_0[t,\hat z_0^{N}(t),z_j^{o}(t)] dw_0(t),\quad \hat z_0^{N}(0)= z_0(0), ~ 0 \leq t \leq T,
\end{align*}
where $z_{-0}^{o}\equiv(z_1^o, \cdots, z_N^o)$ is given by the MV system above. Theorem \ref{Theorem:MCT} and the Gronwall's lemma imply that
\begin{align}
& \sup_{0 \leq t \leq T} \E |z_0^N(t)-\hat z_0^N (t)| = O(1/\sqrt N). \label{MM:CDC:Nash:Major:I3}
\end{align}
We also introduce
\begin{align*}
& d \hat z_0(t) = f_0[t,\hat z_0(t),u_0(t,\hat z_0(t),z_{-0}^{o}(t)),\mu_t] dt + \sigma_0[t,\hat z_0(t),\mu_t]dw_0(t), 
\end{align*}
with initial condition $\hat z_0(0)=z_0(0)$, where $\mu_{(\cdot)}$ is the minor agents' measure given by the MV system above. Again, by Theorem \ref{Theorem:MCT} and the Gronwall's lemma It can be shown that
\begin{align}
& \sup_{0 \leq t \leq T} \E |\hat z_0^N(t)-\hat z_0(t)|= O(1/\sqrt N). \label{MM:CDC:Nash:Major:I4}
\end{align}

({\bf A3}), ({\bf A6}), (\ref{MM:CDC:Nash:Major:I1})-(\ref{MM:CDC:Nash:Major:I4}) and Theorem \ref{Theorem:MCT} yield
\begin{align} \label{eps:Nash:Major1}
& J_{0}^{N}(u_0;u_{-0}^o) \equiv \E \int_0^T \Big((1/N)\sum_{j=1}^N L_0\big[t,z_0^{N}(t),u_0(t,z_0^{N},z_{-0}^{o,N}),z_j^{o,N}(t)\big]\Big)dt \\
& ~ \overset{(\ref{MM:CDC:Nash:Major:I1})}\geq \E \int_0^T  \Big((1/N)\sum_{j=1}^N L_0\big[t,z_0^{N}(t),u_0(t,z_0^{N},z_{-0}^{o}),z_j^o(t)\big]\Big)dt-O( \epsilon_N+1/\sqrt N) \notag \\
& ~ \overset{(\ref{MM:CDC:Nash:Major:I3})} \geq \E \int_0^T  \Big((1/N)\sum_{j=1}^N L_0\big[t,\hat z_0^{N}(t),u_0(t,\hat z_0^{N},z_{-0}^{o}),z_j^o(t)\big]\Big)dt-O( \epsilon_N+1/\sqrt N) \notag \\
& ~ \overset{(\ref{MM:CDC:Nash:Major:I4})} \geq \E \int_0^T  \Big((1/N)\sum_{j=1}^N L_0\big[t,\hat z_0(t),u_0(t,\hat z_0,z_{-0}^{o}),z_j^o(t)\big]\Big)dt-O( \epsilon_N+1/\sqrt N) \notag \\
& ~  \overset{(\ref{MCT:Result})} \geq \E \int_0^T L_0\big[t,\hat z_0(t),u_0(t,\hat z_0,z_{-0}^{o}),\mu_t\big] dt-O( \epsilon_N+1/\sqrt N), \notag
\end{align}
where the appearance of the $\epsilon_N$ term in the first inequality of (\ref{eps:Nash:Major1}) is due to the fact that here the sequence of minor agents' initials $\{z_j^o(0): 1 \leq j \leq N\}$ in the SMV system (\ref{SMV:ENash}) is generated by independent randomized observations on the distribution $F$ given in ({\bf A2}). 

Furthermore, by the construction of the major agent's SMFG system (\ref{CNP:SHJB:Major})-(\ref{CNP:SMV:Major}) (see the major agent's SOCP (\ref{CNP:MajorinfDyn})-(\ref{CNP:MajorinfCos})) we have
\begin{align}
\label{eps:Nash:Major2}
& \E \int_0^T L_0\big[t,\hat z_0(t),u_0(t,\hat z_0,z_{-0}^{o}),\mu_t\big] dt \geq \E \int_0^T L_0\big[t,z_0^o(t),u_0^o(t,z_0^o),\mu_t\big] dt. 
\end{align}
But, Theorem \ref{Theorem:MCT} and (\ref{MM:CDC:Nash:Major:I1}) imply
\begin{align}
\label{eps:Nash:Major3}
&  \E \int_0^T L_0\big[t,z_0^o(t),u_0^o(t,z_0^o),\mu_t\big] dt \\
&~ \overset{(\ref{MCT:Result})} \geq \E \int_0^T  \Big((1/N)\sum_{j=1}^N L_0\big[t,z_0^o(t),u_0(t,z_0^o),z_j^o(t)\big]\Big)dt-O( \epsilon_N+1/\sqrt N) \notag \\
& ~ \overset{(\ref{MM:CDC:Nash:Major:I1})} \geq \E \int_0^T  \Big((1/N)\sum_{j=1}^N L_0\big[t,z_0^{o,N}(t),u_0(t,z_0^{o,N}),z_j^{o,N}(t)\big]\Big)dt-O( \epsilon_N+1/\sqrt N) \notag \\
& ~ \equiv J_{0}^{N}(u_0^o;u_{-0}^o)-O( \epsilon_N+1/\sqrt N). \notag
\end{align}
It follows from (\ref{eps:Nash:Major1})-(\ref{eps:Nash:Major3}) that  $J_0^{N}(u^o_0;u^o_{-0})-O( \epsilon_N+1/\sqrt N) \leq \inf_{u_0 \in \mathcal U_0} J_0^{N}(u_0;u^o_{-0})$.

Case II (strategy change for the minor agents): Without loss of generality, we assume that the first minor agent changes its MF best response control strategy $u^o(t,\omega,x)$ to $u_1\big(t,\omega,x,z_{-1}(t,\omega)\big) \in \mathcal U_1$. This leads to
\begin{align*}
& d z_0^N(t) =\frac{1}{N} \sum_{j=1}^N f_0[t, z_0^N, u_0^o(t,z_0^N), z_j^N] dt +
\frac{1}{N} \sum_{j=1}^N \sigma_0[t,z_0^N,z_j^N]dw_0(t), \\
& d z_1^N(t) =\frac{1}{N} \sum_{j=1}^N f[t, z_1^N,u_1(t,z_1^N,z_{-1}^N),z_j^N]dt + \frac{1}{N} \sum_{j=1}^N \sigma[t,z_1^N, z_j^N] dw_1(t),\\
& d z_2^N(t) =\frac{1}{N} \sum_{j=1}^N f[t, z_2^N,u^o(t,z_2^N),z_j^N]dt + \frac{1}{N} \sum_{j=1}^N \sigma[t,z_2^N, z_j^N] dw_2(t),\\
& \qquad \vdots \\
& d z_N^N(t) =\frac{1}{N} \sum_{j=1}^N f[t, z_N^N,u^o(t,z_N^N),z_j^N]dt + \frac{1}{N} \sum_{j=1}^N \sigma[t,z_N^N, z_j^N] dw_N(t).\\
\end{align*}
By the same argument as in proving Theorem \ref{Theorem:MCT} (see Appendix B in \cite{MNourian_SiamAppendix}) it can be shown that
\begin{align*}
& \sup_{j=0,2,\cdots,N} ~ \sup_{0 \leq t \leq T}  \E |z_j^{o,N}(t)-z_j^N(t)|= O(1/\sqrt N), \\
& \sup_{j=0,2,\cdots,N} ~ \sup_{0 \leq t \leq T}  \E |z_j^{o}(t)-z_j^N(t)|= O(1/\sqrt N).
\end{align*}
Let $\hat z_1^N(\cdot)$ be the solution to
\begin{align*}
& d\hat z_1^{N}(t) =\frac{1}{N} \sum_{j=1}^N f[t,\hat z_1^{N}(t),u_1\big(t,\hat z_1^{N}(t),z_{-1}^{o}(t)\big),z_j^{o}(t)]dt \\
& \qquad + \frac{1}{N} \sum_{j=1}^N \sigma[t,\hat z_1^{N}(t),z_j^{o}(t)] dw_1(t),\quad \hat z_1^{N}(0)= z_1(0), ~ 0 \leq t \leq T,
\end{align*}
where $z_{-1}^{o}\equiv(z_1^o, \cdots, z_N^o)$ is given by the MV system above. Theorem \ref{Theorem:MCT} and the Gronwall's lemma implies that
\begin{align}
& \sup_{0 \leq t \leq T} \E |z_1^N(t)-\hat z_1^N (t)| = O(1/\sqrt N). \label{MM:CDC:Nash:Minor:I3}
\end{align}
We also introduce
\begin{align*}
& d \hat z_1(t) = f[t,\hat z_1(t),u_1(t,\hat z_1(t),z_{-1}^{o}(t)),\mu_t] dt + \sigma[t,\hat z_1(t),\mu_t]dw_1(t), 
\end{align*}
with initial condition $\hat z_1(0)=z_1(0)$, where $\mu_{(\cdot)}$ is the minor agents' measure given by the MV system above. Again, by Theorem \ref{Theorem:MCT} and the Gronwall's lemma It can be shown that
\begin{align}
& \sup_{0 \leq t \leq T} \E |\hat z_1^N(t)-\hat z_1(t)|= O(1/\sqrt N). \label{MM:CDC:Nash:Minor:I4}
\end{align}

Using (\ref{MM:CDC:Nash:Major:I1}) and (\ref{MM:CDC:Nash:Minor:I3})-(\ref{MM:CDC:Nash:Minor:I4}), and by the same argument as in (\ref{eps:Nash:Major1})-(\ref{eps:Nash:Major3}) we can show that $J_1^{N}(u^o_1;u^o_{-1})-O( \epsilon_N+1/\sqrt N) \leq \inf_{u \in \mathcal U_1} J_1^{N}(u_1;u^o_{-1})$. \qed

\section{Conclusion} \label{S9:mm:Conc} This paper studies a stochastic mean field game (SMFG) system for a class of dynamic games involving nonlinear stochastic dynamical systems with major and minor (MM) agents. The SMFG system consists of coupled (i) backward in time stochastic Hamilton-Jacobi-Bellman (SHJB) equations, and (ii) forward in time stochastic McKean-Vlasov (SMV) or stochastic Fokker-Planck-Kolmogorov (SFPK) equations. Existence and uniqueness of the solution to the MM-SMFG system is established by a fixed point argument in the Wasserstein space of random probability measures. In the case that minor agents are coupled to the major agent only through their cost functions, the $\epsilon_N$-Nash equilibrium property of the SMFG best responses is shown for a finite $N$ population system where $\epsilon_N=O(1/\sqrt N)$. As a particular but important case, the results of Nguyen and Huang \cite{NguyenHuang_CDC11} for MM-SMFG linear-quadratic-Gaussian (LQG) systems with homogeneous population are retrieved, and, in addition, the results of this paper are illustrated with a major and minor agent version of a game model of the synchronization of coupled nonlinear oscillators (see Appendices G and H in \cite{MNourian_SiamAppendix}).

\bibliographystyle{plain}
\bibliography{MMref}

\newpage
\setcounter{page}{1}
\begin{center}
{\bf Appendices}\footnote{This document supplies appendices of the paper ``$\epsilon$-Nash Mean Field Game Theory for Nonlinear Stochastic Dynamical Systems with Major and Minor Agents'' by Mojtaba Nourian and Peter E. Caines, provisionally accepted in SIAM J. Control Optim (first submission: Aug. 2012, revised May 2013). Available online at \url{http://arxiv.org/abs/1209.5684}.}
\end{center}

\section*{Appendix A: Proof of Theorem \ref{Theorem:MCT} (McKean-Vlasov Convergence Result)}
\renewcommand{\theequation}{A.\arabic{equation}}
\setcounter{equation}{0}
\renewcommand{\thetheorem}{A.\arabic{theorem}}
\setcounter{theorem}{0}

We will show
\begin{align*}
& \sup_{0 \leq j \leq N} \sup_{0 \leq t \leq T} \E|\hat z^N_j(t)-\bar z_j(t)|^2=O(1/N),
\end{align*}
which implies the result of the theorem by the Cauchy-Schwarz inequality. First by the inequality $(x+y)^2 \leq 2x^2+2y^2$, we have
\begin{align}
& \E |\hat z^N_0(t)-\bar z_0(t)|^2 \notag \\
& \hspace{1cm} \leq 2 \E \Big|\int_0^t
\Big(\frac{1}{N} \sum_{j=1}^N f_0[s,\hat z^N_0,\varphi_0(s,\hat z^N_0),\hat z^N_j]-f_0[s,\bar
z_0,\varphi_0(s,\bar z_0),\mu_s]\Big)ds \Big|^2 \notag \\
& \hspace{1.2cm}+ 2 \E \Big|\int_0^t \Big(\frac{1}{N} \sum_{j=1}^N \sigma_0
[s,\hat z^N_0,\hat z^N_j]-\sigma_0[s,\bar z_0,\mu_s]\Big) dw_0(s) \Big|^2.\notag
\end{align}
By the Cauchy-Schwarz inequality and the properties
of It\^{o} integrals we then obtain 
\begin{align}
\label{MM:MCT:Ineq}
& \E |\hat z^N_0(t)-\bar z_0(t)|^2 \\
& \hspace{0.5cm} \leq 2 t \E\Big(\int_0^t \Big|\frac{1}{N} \sum_{j=1}^N
f_0[s,\hat z^N_0,\varphi_0(s,\hat z^N_0),\hat z^N_j]-f_0[s,\bar z_0,\varphi_0(s,\bar z_0),\mu_s]
\Big|^2ds
\Big) \notag \\
&\hspace{0.7cm} + 2 \E \Big(\int_0^t \Big|\frac{1}{N} \sum_{j=1}^N \sigma_0
[s,\hat z^N_0,\hat z^N_j]-\sigma_0[s,\bar z_0,\mu_s]\Big|^2 ds
\Big). \notag 
\end{align}
Clearly,
\begin{align}
\label{MM:MCT:Ineq111}
& \frac{1}{N} \sum_{j=1}^N f_0[s,\hat z^N_0,\varphi_0(s,\hat z^N_0),\hat z^N_j]-f_0[s,\bar
z_0,\varphi_0(s,\bar z_0),\mu_s]  \\
& \qquad \qquad = \Big(\frac{1}{N} \sum_{j=1}^N f_0[s,\hat z^N_0,\varphi_0(s,\hat
z^N_0),\hat z^N_j]-\frac{1}{N} \sum_{j=1}^N f_0[s,\bar z_0,\varphi_0(s,\bar z_0),\hat z^N_j]\Big)
\notag \\
& \qquad \qquad ~~ + \Big(\frac{1}{N} \sum_{j=1}^N f_0[s,\bar z_0,\varphi_0(s,\bar
z_0),\hat z^N_j]-\frac{1}{N} \sum_{j=1}^N f_0 [s,\bar z_0,\varphi_0(s,\bar z_0),\bar z_j]\Big)
\notag \\
& \qquad \qquad ~~ + \Big(\frac{1}{N} \sum_{j=1}^N f_0[s,\bar z_0,\varphi_0(s,\bar z_0),\bar z_j]
-f_0[s,\bar z_0,\varphi_0(s,\bar z_0),\mu_s]\Big), \notag
\end{align}
and
\begin{align}
& \frac{1}{N} \sum_{j=1}^N \sigma_0[s,\hat z^N_0,\hat z^N_j]-\sigma_0[s,\bar z_0,\mu_s] =
\Big(\frac{1}{N} \sum_{j=1}^N \sigma_0[s,\hat z^N_0,\hat z^N_j]-\frac{1}{N}
\sum_{j=1}^N\sigma_0[s,\bar
z_0,\hat z^N_j]\Big) \notag \\
& ~~ + \Big(\frac{1}{N} \sum_{j=1}^N\sigma_0[s,\bar z_0,\hat z^N_j]-\frac{1}{N}
\sum_{j=1}^N\sigma_0[s,\bar z_0,\bar z_j]\Big)+ \Big(\frac{1}{N} \sum_{j=1}^N\sigma_0[s,\bar
z_0,\bar z_j]-\sigma_0[s,\bar z_0,\mu_s]\Big). \notag
\end{align}

Applying the inequality $(x+y+z)^2 \leq 3(x^2+y^2+z^2)$, and the Lipschitz continuity conditions of $f_0$ and $\varphi_0$ to (\ref{MM:MCT:Ineq111}) we obtain
\begin{align}
&  \label{MM:MCT:Ineq1} \E\Big(\int_0^t \Big|\frac{1}{N} \sum_{j=1}^N
f_0[s,\hat z^N_0,\varphi_0(s,\hat z^N_0),\hat z^N_j]-f_0[s,\bar z_0,\varphi_0(s,\bar z_0),\mu_s]
\Big|^2ds \\
& \quad \leq 3 C \int_0^t \E\big|\hat z^N_0(s)-\bar z_0(s)\big|^2 ds + 3 C \int_0^t
\E\big|\frac{1}{N}\sum_{j=1}^N\hat z^N_j(s)-\bar z_j(s)\big|^2 ds  \notag \\
& \quad ~~ + 3C \int_0^t \E\Big|\frac{1}{N} \sum_{j=1}^N f_0[s,\bar z_0,\varphi_0(s,\bar z_0),\bar
z_j] -f_0[s,\bar z_0,\varphi_0(s,\bar z_0),\mu_s]\Big|^2 ds, \notag
\end{align}
where $C>0$ is a constant independent of $N$. Due to the centring of $g_s[s,\bar z_0,x]:= f_0[s,\bar z_0,\varphi_0(s,\bar z_0),x]-f_0[s,\bar
z_0,\varphi_0(s,\bar z_0),\mu_s]$ with respect to $x$ and the independence of $\bar z_j$ and $\bar
z_{j'}$ when $j \neq j'$, there are no cross terms in the expansion of the last term in (\ref{MM:MCT:Ineq1}), i.e.,
$\E\big(g_s[s,\bar z_0,\bar z_j]g_s[s,\bar z_0,\bar z_{j'}]\big)=\E \E_{\mathcal F_t^{w_0}}\big(g_s[s,\bar z_0,\bar z_j]g_s[s,\bar z_0,\bar z_{j'}]\big)=0$ for $j \neq j'$ (see
\cite{sznitman1991topics}, Page 175). This property together with (\ref{MM:MCT:Ineq1}), the
boundedness of $f_0$ and the inequality $(\sum_{i=1}^N
x_i)^2 \leq N \sum_{i=1}^N x_i^2$ yields
\begin{align}
\label{MM:MCT:Ineq2}
& \E\Big(\int_0^t \Big|\frac{1}{N} \sum_{j=1}^N
f_0[s,\hat z^N_0,\varphi_0(s,\hat z^N_0),\hat z^N_j]-f_0[s,\bar z_0,\varphi_0(s,\bar z_0),\mu_s]
\Big|^2ds \\
& \quad \leq 3 C \int_0^t \E\big|\hat z^N_0(s)-\bar z_0(s)\big|^2 ds +  \frac{3C}{N} \int_0^t
\sum_{j=1}^N \E\big|\hat z^N_j(s)-\bar z_j(s)|^2 ds + \frac{k_1(t)}{N}, \notag
\end{align}
where $k_1(t) \geq 0$ is an increasing function independent of $N$. Similarly, for the second term on the right hand side of (\ref{MM:MCT:Ineq}) we have
\begin{align}
\label{MM:MCT:Ineq3}
& \E \Big(\int_0^t \Big|\frac{1}{N} \sum_{j=1}^N \sigma_0 [s,\hat z^N_0,\hat z^N_j]-\sigma_0[s,\bar
z_0,\mu_s]\Big|^2 ds\Big) \\
& \quad \leq 3 C \int_0^t \E|\hat z^N_0(s)-\bar z_0(s)|^2 ds +  \frac{3C}{N} \int_0^t \sum_{j=1}^N
\E|\hat z^N_j(s)-\bar z_j(s)|^2 ds + \frac{k_1(t)}{N}. \notag 
\end{align}

The inequalities (\ref{MM:MCT:Ineq}), (\ref{MM:MCT:Ineq2}) and (\ref{MM:MCT:Ineq3}) imply that
\begin{align}
\label{MM:MCT:Major}
& \sup_{0 \leq t \leq T} \E |\hat z^N_0(t)-\bar z_0(t)|^2 \leq 6C(T+1) \int_0^T \E|\hat z^N_0(s)-\bar
z_0(s)|^2 ds \\
& \qquad \qquad \qquad + \frac{6C(T+1)}{N} \int_0^T \sum_{j=1}^N \E|\hat z^N_j(s)-\bar
z_j(s)|^2 ds + \frac{2(T+1)k_1(T)}{N}. \notag
\end{align}

Second, by taking a similar approach for the $i^{\textrm{th}}$ minor agent ($1 \leq i \leq N$) we
get
\begin{align}
\label{MM:MCT:Minor}
& \sup_{0 \leq t \leq T} \E |\hat z^N_i(t)-\bar z_i(t)|^2 \leq 8C(T+1) \int_0^T
\E|\hat z^N_i(s)-\bar z_i(s)|^2 ds + \frac{k(T)}{N} \\
& \quad + 8C(T+1) \Big(\int_0^T \E|\hat z^N_0(s)-\bar z_0(s)|^2 ds + \frac{1}{N} \int_0^T
\sum_{j=1}^N \E|\hat z^N_j(s)-\bar z_j(s)|^2 ds\Big), \notag 
\end{align}
where $k(T)>0$ is independent of $N$.

The inequalities (\ref{MM:MCT:Major}) and (\ref{MM:MCT:Minor}) yield
\begin{align}
\label{MM:MCT:sum:mm}
& g^N(T):= \sup_{0 \leq t \leq T} \E |\hat z^N_0(t)-\bar z_0(t)|^2 + \frac{1}{N} \sum_{j=1}^N \sup_{0
\leq t \leq T} \E |\hat z^N_j(t)-\bar z_j(t)|^2 \\
& \quad  \leq 22C(T+1) \int_0^T \Big(\E|\hat z^N_0(s)-\bar z_0(s)|^2 + \frac{1}{N}\sum_{j=1}^N
\E|\hat z^N_j(s)-\bar z_j(s)|^2 \Big) ds  \notag \\
& \qquad + \frac{k_0(T)+k(T)}{N}  \leq 22C(T+1) \int_0^T g(s) ds + \frac{k_0(T)+k(T)}{N}. \notag
\end{align}

It follows from Gronwall's Lemma that
\begin{align}
& g^N (T) \leq \frac{k_0(T)+k(T)}{N} \Big(\exp\big(22C(T+1)T\big)\Big) = O(1/N),
\label{MM:MCT:sum:mm1}
\end{align}
where the right hand side may only depend upon the terminal time $T$. This yields
\begin{align*}
& \sup_{0 \leq t \leq T} \E |\hat z^N_0(t)-\bar z_0(t)|^2 = O(1/N).
\end{align*}
The inequalities (\ref{MM:MCT:Minor}) and (\ref{MM:MCT:sum:mm1}) combined with Gronwall's Lemma
imply that
\begin{align*}
& \sup_{1 \leq i \leq N} \sup_{0 \leq t \leq T} \E |\hat z^N_i(t)-\bar z_i(t)|^2 = O(1/N).
\end{align*}
This completes the proof. \qed

\section*{Appendix B: Extended It\^o-Kunita Formula}
\renewcommand{\theequation}{B.\arabic{equation}}
\setcounter{equation}{0}
\renewcommand{\thetheorem}{B.\arabic{theorem}}
\setcounter{theorem}{0}

We recall an extended version of the It\^o-Kunita formula \cite{kunita1981some} for the composition of stochastic processes (see Theorem 2.3 in \cite{peng1992stochastic}).

\begin{theorem} Let $\phi(t,x)$ be a stochastic process a.s. continuous in $(t,x)$ such that (i) for each $t$, $\phi(t,\cdot)$ is a $C^2(\mathbb{R}^n)$ map a.s., (ii) for each $x$, $\phi(\cdot,x)$ is a continuous semi-martingale represented by
\begin{align*}
&d \phi(t,x) = -\Gamma(t,x) dt + \sum_{k=1}^m \psi_k(t,x) dW_k(t),\qquad (t,x) \in [0,T] \times \mathbb{R}^n,\end{align*}
where $\Gamma(t,x)$ and $\psi_k(t,x)$, $1 \leq k \leq m$, are $\mathcal F_t^W$-adapted stochastic processes which are continuous in $(t,x)$ a.s., such that for each $t$, $\Gamma(t,\cdot)$ is a $C^1(\mathbb{R}^n)$ map a.s., and $\psi_k(t,\cdot)$, $1 \leq k \leq m$, are $C^2(\mathbb{R}^n)$ maps (a.s.).

Let $x(\cdot) = \big(x^1(\cdot), \cdots, x^n(\cdot)\big)$ be a continuous semi-martingale of the form
\begin{align*}
& dx^i(t) = f_i(t) dt + \sum_{k=1}^m \sigma_{ik}(t) dW_k(t)+ \sum_{k=1}^m \varsigma_{ik}(t) dB_k(t), \quad 1 \leq i \leq n,
\end{align*}
where $f_i$, $\sigma_i=(\sigma_{i1},\cdots,\sigma_{im})$ and $\varsigma_i=(\varsigma_{i1},\cdots,\varsigma_{im})$, $1 \leq i \leq n$, are $\mathcal F_t^W$-adapted stochastic processes such that (i) $f_i$ is an integrable process a.s., and (ii) $\sigma_i$ and $\varsigma_i$ are square integrable processes (a.s.). 

Then the composition map $\phi(\cdot,x(\cdot))$ is also a continuous semi-martingale which has the form   
\begin{align} \label{ItoKunita}
& d\phi\big(t,x(t)\big) = -\Gamma\big(t,x(t)\big) dt + \sum_{k=1}^m \psi_k\big(t,x(t)\big) dW_k(t) + \sum_{i=1}^n \partial_{x_i} \phi \big(t,x(t)\big) f_i(t) dt  \\
& \quad + \sum_{i=1}^n \sum_{k=1}^m \partial_{x_i} \phi \big(t,x(t)\big) \sigma_{ik} (t) dW_k(t) + \sum_{i=1}^n \sum_{k=1}^m \partial_{x_i} \phi \big(t,x(t)\big) \varsigma_{ik} (t) dB_k(t)  \notag \\
& \quad + \sum_{i=1}^n \sum_{k=1}^m  \partial_{x_i} \psi_k \big(t,x(t)\big) \sigma_{ik} (t) dt +\frac{1}{2} \sum_{i,j=1}^n \sum_{k=1}^m \partial^2_{x_ix_j} \phi\big(t,x(t)\big)\sigma_{ik} (t) \sigma_{jk}(t) dt  \notag \\
& \quad +\frac{1}{2} \sum_{i,j=1}^n \sum_{k=1}^m \partial^2_{x_ix_j} \phi\big(t,x(t)\big)\varsigma_{ik} (t) \varsigma_{jk}(t) dt. \notag
\end{align} \qed
\end{theorem}

\section*{Appendix C}
\renewcommand{\theequation}{C.\arabic{equation}}
\setcounter{equation}{0}
\renewcommand{\thetheorem}{C.\arabic{theorem}}
\setcounter{theorem}{0}

We may write the functionals of $\mu^0_{(\cdot)}(\omega)$ and $\mu_{(\cdot)}(\omega)$ in (\ref{CNP:MinorinfDyn})-(\ref{CNP:MinorinfCos}) as random functions: \begin{align}
\label{CNP:f*}
& f^*[t,\omega,z_i,u_i]:=f[t,z_i,u_i,\mu^0_t(\omega),\mu_t(\omega)], \quad \sigma^*[t,\omega,z_i]:=\sigma[t,z_i,\mu^0_t(\omega),\mu_t(\omega)], \\
& L^*[t,\omega,z_i,u_i] :=L[t,z_i,u_i,\mu^0_t(\omega),\mu_t(\omega)]. \notag
\end{align}

We have the following proposition where its proof closely resembles that of Proposition \ref{CNP:Proposition:Major} (see Proposition 4 in \cite{HuangMalhameCaines_06}).

\begin{proposition} \label{CNP:Proposition:Minor} Assume ({\bf A3}) holds for $U$. Let $\mu_t(\omega)$, $0 \leq t \leq T$, be a fixed stochastic measure in the set $\mathcal M_\rho^\beta$ with $0 < \beta <1$, and $ \mu^0_{(\cdot)}(\omega)=\Upsilon_0\big( \mu_{(\cdot)}(\omega)\big)\in \mathcal M_\rho^\gamma$, $0 < \gamma <1/2$, be the obtained probability measure of the major agent in Section \ref{CNP:MajorAgentAnalysis}. For $f^*$, $\sigma^*$ and $L^*$ defined in (\ref{CNP:f*}) we have:
\begin{enumerate}[(i)]
\item Under ({\bf A4}) for $f$ and $\sigma$, the functions $f^*[t,\omega,z_i,u_i]$ and $\sigma^*[t,\omega,z_i]$ and their first order derivatives (w.r.t $z_i$) are a.s. continuous and bounded on $[0,T] \times \mathbb{R}^n \times U$ and $[0,T]\times \mathbb{R}^n$. $f^*[t,\omega,z_i,u_i]$ and $\sigma^*[t,\omega,z_i]$ are a.s. Lipschitz continuous in $z_i$. In addition, $f^*[t,\omega,0,0]$ is in the space $L^2_{\mathcal F_t}([0,T];\mathbb{R}^n)$ and $\sigma^*[t,\omega,0]$ is in the space $L^2_{\mathcal F_t}([0,T];\mathbb{R}^{n\times m})$. 

\item Under ({\bf A5}) for $f$, the function $f^*[t,\omega,z_i,u_i]$ is a.s. Lipschitz continuous in $u_i\in U$, i.e., there exist a constant $c>0$ such that
\begin{align*}
& \sup_{t \in [0,T],z_i \in \mathbb{R}^n} \big|f^*[t,\omega,z_i,u_i]-f^*[t,\omega,z_i,u_i']\big| \leq c(\omega)|u_i-u_i'|, \quad (a.s.).
\end{align*}

\item Under ({\bf A6}) for $L$, the function $L^*[t,\omega,z_i,u_i]$ and its first order derivative (w.r.t $z_i$) is a.s. continuous and bounded on $[0,T] \times \mathbb{R}^n \times U$. It is a.s. Lipschitz continuous in $z_i$. In addition, $L^*[t,\omega,0,0] \in L^2_{\mathcal F_t}([0,T];\mathbb{R}_+)$. 

\item Under ({\bf A8}) for $H^{u}$, the set of minimizers
\begin{align*}
&\arg\inf_{u_i\in U} \big\{\big<f^*[t,\omega,z_i,u_i],p\big> + L^*[t,\omega,z_i,u_i]\big\},
\end{align*}
is a singleton for any $p \in \mathbb{R}^n$, and the resulting $u_i$ as a function of $[t,\omega,z_i,p]$ is a.s. continuous in $t$, a.s. Lipschitz continuous in $(z_i,p)$, uniformly with respect to $t$. In addition, $u_i[t,\omega,0,0]$ is in the space $L^2_{\mathcal F_t}([0,T];\mathbb{R}^n)$.
\end{enumerate} \qed
\end{proposition}

\section*{Appendix D: Proof of Theorems \ref{CNP:Analysis:Minor:SMV}}
\renewcommand{\theequation}{D.\arabic{equation}}
\setcounter{equation}{0}
\renewcommand{\thetheorem}{D.\arabic{theorem}}
\setcounter{theorem}{0}

Let $\omega \in \Omega$ be fixed. For given probability measure $\mu_{(\cdot)}(\omega)\in\mathcal M_\rho^\beta$, $0 < \beta <1$, we can show that the law of the process $z_i^o(\cdot,\omega,\omega')$ given in (\ref{CNP:FixedPoint:minoragent:Lambda}), $\Lambda\big(z_i^o(\cdot,\omega,\omega')\big)$, belongs to $\mathcal M_\rho^\beta$, $0 < \beta <1$ (see Theorem \ref{CNP:Analysis:Minor:SMV:measure}). 

We take $\mu_{(\cdot)}(\omega),~\nu_{(\cdot)}(\omega)\in\mathcal M_\rho^\beta$, $0 < \beta <1$. Let $z_i^o(\cdot,\omega,\omega')$ be defined by (\ref{CNP:FixedPoint:minoragent:Lambda}), and similarly $x_i^o(\cdot,\omega,\omega')$ be defined by (\ref{CNP:FixedPoint:minoragent:Lambda}) after replacing $\mu_{(\cdot)}(\omega)$ by $\nu_{(\cdot)}(\omega)$. We have
\begin{align}
\label{CNP:RFP:Inq1}
& \E_{\mathcal F_t^{w_0}} \sup_{0 \leq s \leq t} \big|z_i^o(s,\omega)-x_i^o(s,\omega)\big|^2\\
& \qquad  \leq 2 t \int_0^t \Big|\int_{\mathbb{R}^n\times\mathbb{R}^n} f[s,z_i^o,u^o_i,y,z]d \mu_s^0(\omega)(y)d \mu_s(\omega)(z) \notag \\
& \qquad \qquad \qquad \qquad -\int_{\mathbb{R}^n\times\mathbb{R}^n} f[s,x_i^o,u^o_i,y,z] d\mu_s^0(\omega)(y)d\nu_s(\omega)(z)\Big|^2ds \notag \\
& \qquad \quad + 2 \int_0^t \Big|\int_{\mathbb{R}^n\times\mathbb{R}^n} \sigma[s,z_i^o,y,z]d \mu_s^0(\omega)(y)d \mu_s(\omega)(z) \notag \\
& \qquad \qquad \qquad \qquad -\int_{\mathbb{R}^n\times\mathbb{R}^n} \sigma[s,x_i^o,y,z]d \mu_s^0(\omega)(y)d\nu_s(\omega)(z)\Big|^2 ds. \notag
\end{align}
But,
\begin{align*}
& \Big| \int f[s,z_i^o,u^o_i,y,z] d\mu_s^0(\omega)(y) d\mu_s(\omega)(z)-\int f[s,x_i^o,u^o_i,y,z] d\mu_s^0(\omega)(y)d\nu_s(\omega)(z)\Big|^2 \\
& \qquad \leq 2 C \Big( |z_i^o(s)-x_i^o(s)|^2 +\int_{C_\rho \times C_\rho} |z_s(\omega_1)-z_s(\omega_2)|^2 d\gamma(\omega_1,\omega_2) \Big),
\end{align*}
where $C$ is obtained from the boundedness and Lipschitz continuity of both $f$ and $u^o$, and $\gamma \in \mathcal M(C_\rho \times C_\rho)$ is any coupling of $\mu$ and $\nu$ where $\gamma(A \times C([0,T];\mathbb{R}^n))= \mu(A)$ and $\gamma(C([0,T];\mathbb{R}^n)\times A)= \nu(A)$ for any Borel set $A \in C([0,T];\mathbb{R}^n)$. Taking the infimum over all such $\gamma$ couplings and then using the definition of metrics $\rho_{(\cdot)}$ and $D^\rho_{(\cdot)}$ yields
\begin{align}
\label{CNP:RFP:Inq2}
& \Big| \int f[s,z_i^o,u^o_i,y,z] d\mu_s^0(\omega)(y) d\mu_s(\omega)(z)-\int f[s,x_i^o,u^o_i,y,z] d\mu_s^0(\omega)(y)d\nu_s(\omega)(z)\Big|^2 \\
& \qquad \leq 2 C \Big( \rho_s \big(z_i^o(s),x_i^o(s)\big)+ \big(D_s^\rho(\mu,\nu)\big)^2\Big). \notag
\end{align}

Similarly we have
\begin{align}
\label{CNP:RFP:Inq3}
& \Big| \int\sigma[s,z_i^o,y,z]d \mu_s^0(\omega)(y)d \mu_s(\omega)(z)-\int\sigma[s,x_i^o,y,z]d \mu_s^0(\omega)(y)d\nu_s(\omega)(z)\Big|^2 \\
& \qquad \leq 2 C_1 \Big( \rho_s \big(z_i^o(s),x_i^o(s)\big)+ \big(D_s^\rho(\mu,\nu)\big)^2\Big), \notag
\end{align}
where $C_1$ is obtained from the boundedness and Lipschitz continuity of both $\sigma$.

It follows from (\ref{CNP:RFP:Inq1})-(\ref{CNP:RFP:Inq3}) that
\begin{align}
\label{CNP:RFP:Inq4}
& \E_{\mathcal F_t^{w_0}} \rho_t\big(z_i^o(\omega),x_i^o(\omega)\big) \equiv \E_{\mathcal F_t^{w_0}} \sup_{0 \leq s \leq t} \big|z_i^o(s,\omega)-x_i^o(s,\omega)\big|^2 \wedge 1 \\
& \qquad \leq 2 (Ct + C_1) \int_0^t \Big( \rho_s\big(z_i^o(\omega),x_i^o(\omega)\big)+ \big(D_s^\rho\big(\mu(\omega),\nu(\omega)\big)\big)^2\Big) ds, \notag
\end{align}
which by Gronwall's lemma yields
\begin{align*}
& \E_{\mathcal F_t^{w_0}} \rho_t\big(z_i^o(\omega),x_i^o(\omega)\big) \leq 2 (CT + C_1) \exp\big(2 (CT + C_1)\big)  \int_0^t \Big(D_s^\rho\big(\mu(\omega),\nu(\omega)\big)\Big)^2 ds. 
\end{align*}
This together with the definition of the Wasserstein metric $D_{(\cdot)}^\rho$ leads to the contraction inequality:
\begin{align*}
& \Big(D_t^\rho\big(\mu(\omega),\nu(\omega)\big)\Big)^2 \leq 2 (CT + C_1) \exp\big(2 (CT + C_1)\big) \int_0^t \Big(D_s^\rho\big(\mu(\omega),\nu(\omega)\big)\Big)^2 ds.
\end{align*}
By following a similar argument as in \cite{sznitman1991topics} (Theorem 1.1), we can show that $\{\Lambda^{k}(\mu(\omega)): k \geq 1\}$ forms a Cauchy sequence a.s. in the complete metric space $\mathcal M_\rho^\beta$, $0 < \beta <1$, and converges a.s. to a unique (a.s.) fixed point of $\Lambda$.  \qed

\section*{Appendix E: Proof of Lemma \ref{CNP:Lemma:Sensetivity}}
\renewcommand{\theequation}{E.\arabic{equation}}
\setcounter{equation}{0}
\renewcommand{\thetheorem}{E.\arabic{theorem}}
\setcounter{theorem}{0}

(i) (\ref{CNP:SMV:Major}) gives
\begin{align*}
& z_0(s,\omega) = z_0(0) + \int_0^s \Big( \int_{\mathbb{R}^n} f_0[\tau,z_0,u_0,y] d\mu_\tau(\omega)(y)\Big) d\tau + \int_0^s \sigma_0[\tau] dw_0(\tau,\omega), \\
& z_0'(s,\omega) = z_0(0) + \int_0^s \Big( \int_{\mathbb{R}^n} f_0[\tau,z_0',u_0',y] d\mu_\tau(\omega)(y)\Big) d\tau + \int_0^s \sigma_0[\tau] dw_0(\tau,\omega), 
\end{align*}
corresponding to the control processes $u_0$ and $u_0'$ in $C_{\textrm{Lip}(x)}([0,T]\times\Omega\times\mathbb{R}^n;U_0)$. By the Lipschitz continuity of $f_0$ (see ({\bf A4}) and ({\bf A5})) there are positive constants $C_0$ and $C_1$ such that
\begin{align*}
& |z_0(s,\omega)-z_0'(s,\omega)|^2 \leq 2C_0 s \int_0^s |z_0(\tau,\omega)-z_0'(\tau,\omega)|^2d\tau \\
& \qquad \qquad \qquad \qquad \qquad + 2 C_1 s^2  \sup_{(t,x)\in [0,T]\times\mathbb{R}^n} \big| u_0(t,\omega,x) -u_0'(t,\omega,x)\big|^2.
\end{align*}
The Gronwall's lemma yields
\begin{align*}
& \rho_t\big(z_0(\omega),z_0'(\omega)\big)\leq 2 C_1 t^2 \exp (2C_0 t) \sup_{t,x} \big| u_0(t,\omega,x) -u_0'(t,\omega,x)\big|^2.
\end{align*}
This together with the fact that $\mu^0_t(\omega) = \delta_{z_0(t,\omega)}$ and $\nu^0_t(\omega) = \delta_{z_0'(t,\omega)}$, and the definition of the Wasserstein metric $D_{(\cdot)}^\rho$ leads to (\ref{CNP:ContractionLemmai}) where $c_0:=2 C_1 T^2 \exp (2C_0 T)$.

(ii) We have
\begin{align*}
& z_0(s,\omega) = z_0(0) + \int_0^s \Big( \int_{\mathbb{R}^n} f_0[\tau,z_0,u_0^o,y] d\mu_\tau(\omega)(y)\Big) d\tau + \int_0^s \sigma_0[\tau] dw_0(\tau,\omega), \\
& z_0'(s,\omega) = z_0(0) + \int_0^s \Big( \int_{\mathbb{R}^n} f_0[\tau,z_0',u_0^o,y] d\nu_\tau(\omega)(y)\Big) d\tau + \int_0^s \sigma_0[\tau] dw_0(\tau,\omega), 
\end{align*}
corresponding to the stochastic measures $\mu(\omega),\nu(\omega) \in \mathcal M_\rho^\beta$, $0 < \beta <1$. By the Lipschitz continuity of $f_0$ (see ({\bf A4}) and ({\bf A5})) and $u_0^o$ there are positive constants $C_0$ and $C_1$ such that
\begin{align*}
& |z_0(s,\omega)-z_0'(s,\omega)|^2 \leq 2C_0 s \int_0^s |z_0(\tau,\omega)-z_0'(\tau,\omega)|^2d\tau \\
& \qquad \qquad \qquad \qquad \qquad + 2 C_1 s^2 \Big(D_T^\rho\big(\mu(\omega),\nu(\omega)\big)\Big)^2.
\end{align*}
The Gronwall's lemma yields
\begin{align*}
& \rho_t\big(z_0(\omega),z_0'(\omega)\big)\leq 2 C_1 t^2 \exp (2C_0 t) \Big(D_T^\rho\big(\mu(\omega),\nu(\omega)\big)\Big)^2.
\end{align*}
This together with the fact that $\mu^0_t(\omega) = \delta_{z_0(t,\omega)}$ and $\nu^0_t(\omega) = \delta_{z_0'(t,\omega)}$, and the definition of the Wasserstein metric $D_{(\cdot)}^\rho$ leads to (\ref{CNP:ContractionLemmaii}) where $c_1:=2 C_1 T^2 \exp (2C_0 T)$.

 (iii) (\ref{CNP:SMV:Minor}) gives
\begin{align*}
& z_i(s,\omega,\omega') = z_i(0) + \int_0^t \Big(\int_{\mathbb{R}^n}  \int_{\mathbb{R}^n} f[s,z_i,u,y,z] d\mu_s^0(\omega)(y)  d\mu_s(\omega)(z)\Big) ds \\
& \qquad \qquad \qquad + \int_0^t \Big(\int_{\mathbb{R}^n}  \int_{\mathbb{R}^n} \sigma[s,z_i,y,z] d\mu_s^0(\omega)(y) d\mu_s(\omega)(z)\Big) dw_i(s,\omega'), \notag \\
& z_i'(s,\omega,\omega') = z_i(0) + \int_0^t \Big(\int_{\mathbb{R}^n}  \int_{\mathbb{R}^n} f[s,z_i',u',y,z] d\mu_s^0(\omega)(y)  d\nu_s(\omega)(z)\Big) ds \\
& \qquad \qquad \qquad + \int_0^t \Big(\int_{\mathbb{R}^n}  \int_{\mathbb{R}^n} \sigma[s,z_i',y,z] d\mu_s^0(\omega)(y) d\nu_s(\omega)(z)\Big) dw_i(s,\omega'), \notag
\end{align*}
corresponding to the control processes $u$ and $u'$ in $C_{\textrm{Lip}(x)}([0,T]\times\Omega\times\mathbb{R}^n;U)$. By the Lipschitz continuity of $f$ and $\sigma$ (see ({\bf A4}) and ({\bf A5})) there are positive constants $C_0, C_1$ and $C_2$ such that
\begin{align*}
& \E_{\omega}|z_i(s,\omega,\omega')-z_i'(s,\omega,\omega')|^2 \leq 2 (3C_0 s+2 C_1) \E_{\omega} \int_0^s |z_0(\tau,\omega)-z_0'(\tau,\omega)|^2d\tau \\
&\qquad \qquad \qquad \qquad + 2 (3C_0 s+2 C_1) \E_{\omega} \int_0^s  \Big(D_\tau^\rho\big(\mu(\omega),\nu(\omega)\big)\Big)^2d\tau \\
&\qquad \qquad \qquad \qquad  + 6 C_2 s^2 \sup_{t,x} \E_{\omega} \big| u(t,\omega,x) -u'(t,\omega,x)\big|^2.
\end{align*}
The Gronwall's lemma yields
\begin{align*}
& \rho_t\big(z_i(s,\omega),z_i'(s,\omega)\big)\leq 2 (3C_0 t+2 C_1)  \exp\big(2 (3C_0 t+2 C_1)\big) \int_0^t  \Big(D_\tau^\rho\big(\mu(\omega),\nu(\omega)\big)\Big)^2d\tau  \\
& \qquad \qquad \qquad  + 6 C_2 t^2 \exp\big(2 (3C_0 t+2 C_1)\big) \sup_{t,x} \big| u(t,\omega,x) -u'(t,\omega,x)\big|^2.
\end{align*}
This together with the definition of the Wasserstein metric $D_{(\cdot)}^\rho$ leads to
\begin{align*}
& \Big(D_T^\rho\big(\mu(\omega),\nu(\omega)\big)\Big)^2 \leq K(T) \int_0^T  \Big(D_\tau^\rho\big(\mu(\omega),\nu(\omega)\big)\Big)^2d\tau \\
& \qquad \qquad \qquad \qquad \qquad  + K'(T) \sup_{t,x} \big| u(t,\omega,x) -u'(t,\omega,x)\big|^2,
\end{align*}
where $K(T):=2 (3C_0 T+2 C_1)  \exp\big(2 (3C_0 T+2 C_1)\big)$ and $K'(T):=6 C_2 T^2 \exp\big(2 (3C_0 T+2 C_1)\big)$. Applying the Gronwall's lemma gives (\ref{CNP:ContractionLemmaiii}) with $c_2:=K'(T) \exp(K(T))$. 

(iv) The proof of this part closely resembles that of Part (iii). \qed 

\section*{Appendix F: The Sensitivity Analysis of the SHJB Equations} \label{sectionsensitivity}
\renewcommand{\theequation}{F.\arabic{equation}}
\setcounter{equation}{0}
\renewcommand{\thetheorem}{F.\arabic{theorem}}
\setcounter{theorem}{0}

In this section we study the sensitivity of the major and minor agents' SHJB equations (\ref{CNP:SHJB:Major}) and (\ref{CNP:SHJB:Minor}) to the stochastic measures $\mu_{(\cdot)}(\omega)$ and $\mu^0_{(\cdot)}(\omega)$ in order to show the feedback regularity conditions. The analysis of this section is based on the framework of Section 6 of \cite{kolokoltsov2011mean}.

First we consider a family of stochastic optimal control problems (SOCP) (\ref{CNP:GenAgeDyn})-(\ref{CNP:GenAgeCos}) parameterized by $\alpha \in \mathbb{R}$. In this $\alpha$-parameterized formulation called $(\textrm{SOCP})_\alpha$: (i) the dynamics of the states $z^\alpha(t,\omega)$, denoted by $(\ref{CNP:GenAgeDyn})_\alpha$, are of the form (\ref{CNP:GenAgeDyn}) with $f[t,\omega,z,u]$, $\sigma[t,\omega,z]$ and $\varsigma[t,\omega,z]$ replaced by $f^\alpha[t,\omega,z^\alpha,u^\alpha]$, $\sigma^\alpha[t,\omega,z^\alpha]$ and $\varsigma^\alpha[t,\omega,z^\alpha]$, respectively, and (ii) the cost functions $J^\alpha(u^\alpha)$, denoted by $(\ref{CNP:GenAgeCos})_\alpha$, are of the form (\ref{CNP:GenAgeCos}) with $L[t,\omega,z,u]$ replaced by $L^\alpha[t,\omega,z^\alpha,u^\alpha]$.    

The value functions $\phi^\alpha (\cdot,x(\cdot))$ correspond to the $(\textrm{SOCP})_\alpha$ are defined similar to (\ref{PengValueFun}) with $L[t,\omega,z,u]$ replaced by $L^\alpha[t,\omega,z^\alpha,u^\alpha]$. Based on \cite{peng1992stochastic} we shall restrict to the case where $\phi^\alpha(\cdot,x(\cdot))$ are semi-martingales of the form (\ref{SMR}) with $\Gamma(\cdot,x(\cdot))$ and $\psi(\cdot,x(\cdot))$ are replaced by $\Gamma^\alpha(\cdot,x(\cdot))$ and $\psi^\alpha(\cdot,x(\cdot))$, respectively.

If the $\alpha$-parameterized family of processes $\phi^\alpha(t,x)$, $\Gamma^\alpha(t,x)$ and $\psi^\alpha(t,x)$ are a.s. continuous in $(x,t)$ and are smooth enough with respect to $x$, then by using the analysis in \cite{peng1992stochastic} we can show that the pairs $\big(\phi^\alpha(s,x),\psi^\alpha(s,x)\big)$ satisfy the following backward in time $\alpha$-parameterized stochastic Hamilton-Jacobi-Bellman $(\textrm{SHJB})_\alpha$ equations: 
\begin{align}
& \label{CNP:SHJB:alpha}
-d\phi^\alpha(t,\omega,x)=\Big[H^\alpha[t,\omega,x,D_x\phi^\alpha(t,\omega,x)]+\big<\sigma^\alpha[t,\omega,x],D_x\psi^\alpha(t,\omega,x)\big> \\ 
& +\frac{1}{2}\tr\big(a^\alpha[t,\omega,x]D_{xx}^2 \phi^\alpha(t,\omega,x)\big)\Big]dt
-(\psi^\alpha)^T(t,\omega,x)dW(t,\omega),~ \phi^\alpha(T,x) =0, \notag
\end{align}
where $a^\alpha[t,\omega,x]:=\sigma^\alpha[t,\omega,x]\big(\sigma^\alpha[t,\omega,x]\big)^T+\varsigma^\alpha[t,\omega,x]\big(\varsigma^\alpha[t,\omega,x]\big)^T$
, and the stochastic Hamiltonians $H^\alpha:[0,T] \times\Omega\times\mathbb{R}^n\times\mathbb{R}^n\rightarrow \mathbb{R}$ are given by
\begin{align*}
& H^\alpha[t,\omega,x,p]:=\inf_{u^\alpha \in \mathcal U}\big\{\big<f^\alpha[t,\omega,x,u], p\big> + L^\alpha[t,\omega,x,u]\big\}.
\end{align*}

Suppose the assumptions ({\bf H1})-({\bf H3}) hold for $(f^\alpha,L^\alpha,\sigma^\alpha,\varsigma^\alpha)$. Then the $(\textrm{SHJB})_\alpha$ equations (\ref{CNP:SHJB:alpha}) have unique solutions (see Theorem \ref{CNP:PenTheorem} or Theorem 4.1 in \cite{peng1992stochastic}):
\begin{align*}
& (\phi^\alpha(t,x),\psi^\alpha(t,x)) \in\big(L_{\mathcal F_t}^2([0,T];\mathbb{R}),L_{\mathcal F_t}^2([0,T];\mathbb{R}^m)\big), \qquad\forall  \alpha \in \mathbb{R}.
\end{align*}

The forward in time $\mathcal F_t^W$-adapted optimal control processes of the $(\textrm{SOCP})_\alpha$
$(\ref{CNP:GenAgeDyn})_\alpha$-$(\ref{CNP:GenAgeCos})_\alpha$ are given by (see \cite{peng1992stochastic}) 
\begin{align}
&  \label{CNP:OC:alpha} u^{\alpha,o} (t,\omega,x) :=\arg\inf_{u^\alpha\in U} H^{\alpha,u}[t,\omega,x,D_x\phi^\alpha(t,\omega,x),u^\alpha]  \\
& \qquad \qquad ~~~~ = \arg\inf_{u^\alpha \in U}\big\{\big<f^\alpha[t,\omega,x,u^\alpha],D_x\phi^\alpha(t,\omega,x)\big> +
L^\alpha[t,\omega,x,u^\alpha]\big\}. \notag
\end{align} 

We set 
\begin{align*}
& g^\alpha[t,\omega,x,\phi^\alpha(t,\omega,x),\psi^\alpha(t,\omega,x)] := H^\alpha[t,\omega,x,D_x\phi^\alpha(t,\omega,x)] \\
& \hspace {6cm} +\big<\sigma^\alpha[t,\omega,x],D_x\psi^\alpha(t,\omega,x)\big>, \\
& A^\alpha (t,\omega,x) (\cdot) := \frac{1}{2}\tr\big(a^\alpha[t,\omega,x]D_{xx}^2 (\cdot) \big),
\end{align*}
where $A^\alpha$ in $[0,T]\times\Omega\times\mathbb{R}^n$ is an operator on $C^2(\mathbb{R}^n)$. We may now rewrite the backward in time $\alpha$-parameterized $(\textrm{SHJB})_\alpha$ equations (\ref{CNP:SHJB:alpha}) as
 \begin{align}
& \label{CNP:SHJB:alpha:operator}
d\phi^\alpha(t,\omega,x) + A^\alpha (t,\omega,x) \big(\phi^\alpha(t,\omega,x)\big) dt  \\ 
& \hspace{1cm} = - g^\alpha[t,\omega,x, \phi^\alpha(t,\omega,x),\psi^\alpha(t,\omega,x)] dt +(\psi^\alpha)^T(t,\omega,x)dW(t,\omega), \notag
\end{align}
with $ \phi^\alpha(T,x) =0$. 

At this point we introduce the mild form of (\ref{CNP:SHJB:alpha:operator}) because this form is more suitable for the sensitivity analysis of this section. We note that it is sufficient to consider the mild solution in the analysis of existence and uniqueness of solutions to the SMFG system. 

If the pair $ (\phi^\alpha(t,x),\psi^\alpha(t,x))$ is a smooth solution to (\ref{CNP:SHJB:alpha:operator}) that satisfies the following mild form by a Duhamel Principle \cite{kolokoltsov2011mean}:
 \begin{align}
&  \label{CNP:SHJB:alpha:mild}
\phi^\alpha(t,\omega,x) = \int_t^T \exp\Big(\int_t^s A^\alpha (\tau,\omega,x) d \tau\Big) \Big(g^\alpha[s,\omega,x,\phi^\alpha(s,\omega,x),\psi^\alpha(s,\omega,x)]\Big) ds \\ 
& \hspace{2cm} - \int_t^T \exp\Big(\int_t^s A^\alpha (\tau,\omega,x) d \tau\Big) \Big((\psi^\alpha)^T(s,\omega,x)\Big) dW(s,\omega). \notag
\end{align}

We define the operators:
 \begin{align*}
& \Phi^\alpha(t,s,\omega,x)(\cdot) = \exp\Big(\int_t^s A^\alpha (\tau,\omega,x) (\cdot) d \tau \Big) \equiv \exp\Big(\int_t^s  \frac{1}{2}\tr\big(a^\alpha[\tau,\omega,x]D_{xx}^2 (\cdot) \big) d \tau, \\
& \Psi^\alpha(t,s,\omega,x)(\cdot) = \int_t^s \partial_\alpha A^\alpha (\tau,\omega,x) (\cdot) d \tau \equiv \int_t^s  \frac{1}{2}\tr\big(\partial_\alpha a^\alpha[\tau,\omega,x]D_{xx}^2 (\cdot)\big) d \tau,
\end{align*}
in $[0,T]\times\Omega\times\mathbb{R}^n$ which are maps on $C^\f(\mathbb{R}^n)$ and $C^2(\mathbb{R}^n)$, respectively.

Differentiating (\ref{CNP:SHJB:alpha:mild}) with respect to $\alpha$ gives
 \begin{align}
&  \label{CNP:SHJB:alpha:mild:alphaDiff}
\partial_\alpha \phi^\alpha(t,\omega,x) = \int_t^T\big( \Phi^\alpha(t,s,\omega,x)\big) \big(\Psi^\alpha(t,s,\omega,x) \big) \\
& \hspace{4cm} \Big(g^\alpha[s,\omega,x,\phi^\alpha(s,\omega,x),\psi^\alpha(s,\omega,x)]\Big) ds \notag \\ 
& + \int_t^T\big( \Phi^\alpha(t,s,\omega,x)\big) \Big(\partial_\alpha g^\alpha[s,\omega,x,\phi^\alpha(s,\omega,x),\psi^\alpha(s,\omega,x)]\Big) ds \notag \\ 
& - \int_t^T\big( \Phi^\alpha(t,s,\omega,x)\big)  \big( \Psi^\alpha(t,s,\omega,x)\big) \Big((\psi^\alpha)^T(s,\omega,x)\Big)dW(s,\omega) \notag\\
&  - \int_t^T \big( \Phi^\alpha(t,s,\omega,x)\big) \Big((\partial_\alpha \psi^\alpha)^T(s,\omega,x)\Big)dW(s,\omega),\notag
\end{align}
where 
\begin{align*}
&  \partial_\alpha g^\alpha[t,\omega,x,\phi^\alpha(t,\omega,x),\psi^\alpha(t,\omega,x)] \equiv \partial_\alpha H^\alpha[t,\omega,x,D_x\phi^\alpha(t,\omega,x)]  \\
&  \hspace {1cm} + \partial_p H^\alpha[t,\omega,x,D_x\phi^\alpha(t,\omega,x)] D_x \big(\partial_\alpha \phi^\alpha(t,\omega,x) \big) \\
& \hspace {1cm} +\big<\partial_\alpha \sigma^\alpha[t,\omega,x],D_x\psi^\alpha(t,\omega,x)\big> + \big<\sigma^\alpha[t,\omega,x],D_x \big( \partial_\alpha\psi^\alpha(t,\omega,x) \big)\big> .
\end{align*}

We may rewrite (\ref{CNP:SHJB:alpha:mild:alphaDiff}) as
 \begin{align}
& \label{CNP:SHJB:alpha:mild:alphaDiff2}
\partial_\alpha \phi^\alpha(t,\omega,x) = \int_t^T\big( \Phi^\alpha(t,s,\omega,x)\big) A_1^\alpha (s,\omega,x) \big( \partial_\alpha  \phi^\alpha(t,\omega,x) \big) ds \\  
& \qquad \qquad  +  \int_t^T\big( \Phi^\alpha(t,s,\omega,x)\big) \Big(h^\alpha_1[t,s,\omega,x,\partial_\alpha \psi^\alpha]\Big) ds \notag \\ 
& \qquad \qquad  - \int_t^T \big( \Phi^\alpha(t,s,\omega,x)\big) \Big((\partial_\alpha \psi^\alpha)^T(s,\omega,x)\Big)dW(s,\omega), \notag \\
& \qquad \qquad - \int_t^T\big( \Phi^\alpha(t,s,\omega,x)\big) \Big(h^\alpha_2 [t,s,\omega,x]\Big) dW(s,\omega), \notag
\end{align}
where 
\begin{align*}
& A_1^\alpha (s,\omega,x) (\cdot) := \partial_p H^\alpha[s,\omega,x,D_x\phi^\alpha(s,\omega,x)] D_x (\cdot), \\
& h^\alpha_1[t,s,\omega,x,\partial_\alpha \psi^\alpha] := \big( \Psi^\alpha(t,s,\omega,x)\big) \Big(g^\alpha[s,\omega,x,\phi^\alpha(s,\omega,x),\psi^\alpha(s,\omega,x)]\Big) \\
& \qquad + \partial_\alpha H^\alpha[s,\omega,x,D_x\phi^\alpha(s,\omega,x)] +\big<\partial_\alpha \sigma^\alpha[s,\omega,x],D_x\psi^\alpha(s,\omega,x)\big> \\
& \qquad + \big<\sigma^\alpha[s,\omega,x],D_x \big( \partial_\alpha\psi^\alpha\big)\big>,\\
& h^\alpha_2 [t,s,\omega,x] :=  \big( \Psi^\alpha(t,s,\omega,x)\big) \Big((\psi^\alpha)^T(s,\omega,x)\Big).
\end{align*} 

We introduce the following assumption:

({\bf H5}) $\partial_\alpha f^\alpha[t,x,u]$, $\partial_\alpha L^\alpha [t,x,u]$, $\partial_\alpha \sigma^\alpha[t,x]$ and $\partial_\alpha \varsigma^\alpha[t,x]$ exist and are $C^\f (\mathbb{R}^n)$. Assume ({\bf H1})-({\bf H3}) hold where $(f,L,\sigma,\varsigma)$ are replaced by $(\partial_\alpha f^\alpha,\partial_\alpha L^\alpha,\partial_\alpha \sigma^\alpha,\partial_\alpha \varsigma^\alpha)$, and all the boundedness assumptions are uniformly.

\begin{proposition} \label{Proposition;alphadependent} Assume ({\bf H11})-({\bf H3}) hold for $(f^\alpha,L^\alpha,\sigma^\alpha,\varsigma^\alpha)$. Let the pair $(\phi^\alpha(t,x),\psi^\alpha(t,x))$ be the unique solution to (\ref{CNP:SHJB:alpha}) which are $C^\f (\mathbb{R}^n)$ and a.s. uniformly bounded. In addition, we assume ({\bf H5}) holds. Then, the equation (\ref{CNP:SHJB:alpha:mild:alphaDiff}) has a unique solution 
\begin{align*} 
& (\partial_\alpha \phi(t,x),\partial_\alpha \psi(t,x))\in \big(L_{\mathcal F_t}^2([0,T];\mathbb{R}),L_{\mathcal F_t}^2([0,T];\mathbb{R}^m)\big)
\end{align*}
such that $\sup_{0 \leq t \leq T}|D_x \partial_\alpha  \phi(t,\cdot)| < \infty$ (a.s.).
\end{proposition} 

{\it Proof}: The proof of existence and uniqueness of solution to (\ref{CNP:SHJB:alpha:mild:alphaDiff2}) follows from Theorem 4.1 in \cite{hu1991adapted} (see the proof of Theorem 4.1 in \cite{peng1992stochastic}, see also \cite{ma1997adapted,ma1999linear,hu2002semi} or Chapter 5 of \cite{ma1999forward}). By taking the conditional expectation $\E_{\mathcal F_t^{w_0}}$ of the square of both sides of (\ref{CNP:SHJB:alpha:mild:alphaDiff2}) and the boundedness assumptions in the theorem, one can show $\sup_{0 \leq t \leq T}|\partial_\alpha \phi(t,\cdot)| < \infty$ (a.s.) (see the proof of Theorem 2.1 in \cite{peng1992stochastic}). Using this in equation (\ref{CNP:SHJB:alpha:mild:alphaDiff2}) implies the boundedness of $D_x \partial_\alpha \phi(t,\cdot)$.  \qed

\section*{Appendix G: The Major and Minor (MM) SMFG Linear-Qudratic-Gaussian (LQG) System} \label{Appendix:mm:LQG} 
\renewcommand{\theequation}{G.\arabic{equation}}
\setcounter{equation}{0}
\renewcommand{\thetheorem}{G.\arabic{theorem}}
\setcounter{theorem}{0}

We consider the MM LQG dynamic game problem of \cite{HuangSIAM2010}. In this case all functions in
(\ref{CNP:Major:GenDyn})-(\ref{CNP:Minor:GenCost}) are given by (see Remark \ref{RemarkforLQGcase})
\begin{align*}
& f_0 [t,z_0^N(t),u_0^N(t),z_j^N(t)] = A_0 z_0^N(t) + B_0 u_0^N(t) +F_0 z_j^N(t),  \\
& f [t,z_i^N(t),u_i^N(t),z_0^N(t),z_j^N(t)] = A z_i^N(t) + B u_i^N(t) + F z_j^N(t)+G z_0^N(t), \\
& \sigma_0[t,z_0^N(t),z_j^N(t)]=S_0, \qquad  \hspace{2cm}  \sigma [t,z_i^N(t),z_0^N(t),z_j^N(t)] =S,\\
& L_0[t,z_0^N(t),u_0^N(t),z_j^N(t)] = \Big[z_0^N(t)- \Big(H_0 \big(\frac{1}{N} \sum_{j=1}^N z_j^N(t) \big) +\eta_0\Big)\Big]^T Q_0   \\
&  \qquad \hspace{1cm} \times \Big[z_0^N(t)- \Big(H_0 \big(\frac{1}{N} \sum_{j=1}^N z_j^N(t)\big)+\eta_0\Big)\Big]^T+ (u_0^N(t))^T R_0 u_0^N(t), \\
& L[t,z_i^N(t),u_i^N(t),z_0^N(t),z_j^N(t)] = \Big[z_i^N(t) -\Big(H z_0^N(t) +\hat H  \big(\frac{1}{N} \sum_{j=1}^N z_j^N(t) \big) +\eta\Big)\Big]^T Q \\
& \qquad \hspace{1cm} \times\Big[z_i^N(t) -\Big(H z_0^N(t) +\hat H  \big(\frac{1}{N} \sum_{j=1}^N z_j^N(t) \big)  +\eta\Big)\Big] + (u_i^N(t))^T R u_i^N(t),
\end{align*}
with the deterministic constant matrices: (i) $A_0, F_0, A, F, G, H_0, H$ and $\hat H$ in $\mathbb{R}^{n\times n}$, (ii) $B_0$ and $B$ in $\mathbb{R}^{n\times k}$, (iii) $S_0$ and $S$ in $\mathbb{R}^{n\times m}$, (iv) the symmetric nonnegative definite matrices $Q_0$ and $Q$ in $\mathbb{R}^{n\times n}$, (v) the symmetric positive definite matrices $R_0$ and $R$ in $\mathbb{R}^{k\times k}$, and the deterministic constant vectors $\eta$ and $\eta_0$ are in $\mathbb{R}^{n}$. 

In this formulation the major agent's SMFG system (\ref{CNP:SHJB:Major})-(\ref{CNP:SMV:Major}) is of the form 
\begin{align}
\label{CNP:SHJB:Major:LQG}
& \hspace{-0.2cm} -d\phi_0(t,\omega,x)=\Big[\big<A_0 x -\frac{1}{4} B_0 R_0^{-1} B_0^T D_x\phi_0(t,\omega,x) +F_0 z^o (t,\omega),D_x\phi_0(t,\omega,x)\big>\\
& \quad + \big<x- (H_0 z^o (t,\omega) +\eta_0),Q_0 \big(x- (H_0 z^o (t,\omega) +\eta_0) \big)\big> \notag \\
& \quad +\big<S_0,D_x\psi_0(t,\omega,x)\big> +\frac{1}{2}\tr\big((S_0^T S_0) D_{xx}^2 \phi_0(t,\omega,x)\big)\Big]dt \notag \\
& \quad  -\psi^T_0(t,\omega,x)dw_0(t,\omega), \hspace{5cm} \phi_0(T,x) =0, \notag\\
\label{CNP:BR:Major:LQG}
&  u^o_0(t,\omega,x) = -\frac{1}{2} R_0^{-1} B_0^T D_x\phi_0(t,\omega,x), \\ 
\label{CNP:SMV:Major:LQG}
& dz_0^o(t,\omega) = \Big[ A_0 z_0^o(t,\omega) + B_0 u^o_0(t,\omega,z_0^o) +F_0 z^o (t,\omega) \Big] dt
\\
& \quad + S_0 dw_0(t,\omega), \hspace{6.2cm} z_0^o(0)=z_0(0), \notag
\end{align}
and the minor agents' SMFG system (\ref{CNP:SHJB:Minor})-(\ref{CNP:SMV:Minor}) is given by
\begin{align}
\label{CNP:SHJB:Minor:LQG}
& \hspace{-0.3cm} -d\phi(t,\omega,x)=\Big[\big<A x -\frac{1}{4} B R^{-1} B^T D_x\phi(t,\omega,x) +F x + Gz_0^o (t,\omega),D_x\phi(t,\omega,x)\big> \\
& \quad  + \big<x- (H z_0^o (t,\omega) + \hat H x +\eta),Q \big(x- (H z_0^o (t,\omega) + \hat H x +\eta) \big)\big> \notag \\
&\quad   +\frac{1}{2}\tr\big((S^T S) D_{xx}^2 \phi(t,\omega,x)\big)\Big]dt -\psi^T(t,\omega,x)dw(t,\omega), \hspace{0.4cm} \phi_0(T,x) =0, \notag\\
\label{CNP:BR:Minor:LQG}
& u^o(t,\omega,x) = -\frac{1}{2}R^{-1} B^T D_x\phi(t,\omega,x), \\
\label{CNP:SMV:Minor:LQG}
& dz^o(t,\omega) = \Big[ A z^o(t,\omega) + B u^o(t,\omega,z^o) + F_0 z^o (t,\omega) + G z_0^o (t,\omega)  \Big] dt
\\
& \quad  + S dw(t,\omega), \hspace{6.5cm} z_0^o(0)=z_0(0). \notag
\end{align}

Let $\Pi_0(\cdot) \geq 0$ be the unique solution of the deterministic Riccati equation
\begin{align*}
& \partial_t \Pi_0(t)+ \Pi_0(t) A_0 + A_0^T \Pi_0(t) - \Pi_0(t) B_0 R_0^{-1} B_0^T \Pi_0(t) + Q_0 =0, \quad  \Pi_0(T)=0.
\end{align*}
We denote $\mathbb{A}_0(\cdot)= A_0 - B_0 R_0^{-1}B_0^T \Pi_0(\cdot)$. It can be verified that the pair $(\phi_0,\psi_0)(t,\omega,x)$ in (\ref{CNP:SHJB:Major}) is given by
\begin{align*}
& \phi_0(t,\omega,x)=x^T \Pi_0(t) x + 2 x^T s_0 (t,\omega)  + g_0(t,\omega), \notag \\
& \psi_0^T(t,\omega,x)=2 x^T q_0 (t,\omega)  + h_0(t,\omega),
\end{align*}
where $(s_0,q_0)(t,\omega)$ and $(g_0,h_0)(t,\omega)$ are unique solutions of the following Backward Stochastic Differential Equations (BSDEs):
\begin{align*}
& -ds_0(t,\omega) =\Big[\mathbb{A}_0^T(t) s_0(t,\omega) + \big( \Pi_0(t) F_0 - Q_0 H_0\big)z^o(t,\omega) - Q_0 \eta_0\Big]dt \\
& \hspace{2cm} -q_0(t,\omega)dw_0(t,\omega), \hspace{5cm}  s_0(T)=0,\\
& -dg_0(t,\omega) =\Big[- s_0^T(t,\omega) B_0 R_0^{-1} B_0^T  s_0(t,\omega) + 2 F_0 z^o(t,\omega) + 2~ \tr \big(S_0^T q_0(t,\omega)\big)\\
& \hspace{2.4cm} + \big(H_0z^o(t,\omega) + \eta_0\big)^T Q_0 \big(H_0z^o(t,\omega) + \eta_0\big)+ \tr \big(S_0^T S_0 \Pi_0(t)\big) \Big]dt \\
& \hspace{2cm} - h_0(t,\omega) dw_0(t,\omega), \hspace{5cm} g_0(T)=0.
\end{align*}

We may now express the major agent's SMFG LQG system (\ref{CNP:SHJB:Major:LQG})-(\ref{CNP:SMV:Major:LQG}) in the following form:
\begin{align*}
& -ds_0(t,\omega) =\Big[\mathbb{A}_0^T(t) s_0(t,\omega) + \big( \Pi_0(t) F_0 - Q_0 H_0\big)z^o(t,\omega) - Q_0 \eta_0\Big]dt \\
& \hspace{2cm} -q_0(t,\omega)dw_0(t,\omega), \hspace{4.2cm} s_0(T)=0,  \\
& u^o_0(t,\omega) = - R_0^{-1} B_0^T \big( \Pi_0(t) z_0^o(t,\omega) + s_0(t,\omega)\big), \\
& dz_0^o(t,\omega) = \Big[ \mathbb{A}_0 (t) z_0^o(t,\omega)- B_0 R_0^{-1} B_0^T  \Pi_0(t) s_0(t,\omega) +F_0 z^o (t,\omega) \Big] dt
\\
& \hspace{2cm} + S_0 dw_0(t,\omega), \hspace{5cm} z_0^o(0)=z_0(0),
\end{align*}
where $z^o(t,\omega)$ is the mean field behaviour of the minor agents (see the minor agents' SMFG LQG system below). 

In a similar way, let $\Pi(\cdot) \geq 0$ be the unique solution of the deterministic Riccati equation
\begin{align*}
& \partial_t \Pi(t)+ \Pi(t) A + A^T \Pi(t) - \Pi(t) B R^{-1} B^T \Pi(t) + Q =0, \quad  \Pi(T)=0.
\end{align*}
We denote $\mathbb{A}(\cdot)= A - B R^{-1}B^T \Pi(\cdot)$. It can be verified that the pair $(\phi,\psi)(t,\omega,x)$ in (\ref{CNP:SHJB:Minor}) is given by
\begin{align*}
& \phi(t,\omega,x)=x^T \Pi(t) x + 2 x^T s (t,\omega)  + g(t,\omega), \notag \\
& \psi^T(t,\omega,x)=2 x^T q (t,\omega)  + h(t,\omega),
\end{align*}
where $(s,q)(t,\omega)$ and $(g,h)(t,\omega)$ are unique solutions of the following BSDEs:
\begin{align*}
& -ds(t,\omega) =\Big[\mathbb{A}^T(t) s(t,\omega) + \big( \Pi(t) F - Q \hat H\big)z^o(t,\omega) + \big( \Pi(t) G - Q  H\big)z^o_0(t,\omega) \\
& \hspace{2cm} - Q\eta \Big]dt - q(t,\omega)dw_0(t,\omega), \hspace{3.6cm}  s(T)=0,\\
& -dg(t,\omega) =\Big[- s^T(t,\omega) B R^{-1} B^T  s(t,\omega) + 2 F z^o(t,\omega) + 2 G z^o_0(t,\omega)\\
& \hspace{2.4cm} + \big(\hat H z^0(t,\omega) + H z^o_0 (t,\omega) + \eta\big)^T Q_0 \big(\hat H z^0(t,\omega) + H z^o_0 (t,\omega) + \eta\big) \\
& \hspace{2.4cm} + \tr \big(S^T S \Pi(t)\big) \Big]dt - h(t,\omega) dw_0(t,\omega), \hspace{1.8cm} g(T)=0.
\end{align*}

We may now express the minor agents' SMFG LQG system (\ref{CNP:SHJB:Minor:LQG})-(\ref{CNP:SMV:Minor:LQG}) in the following form:
\begin{align*}
& -ds(t,\omega) =\Big[\mathbb{A}^T(t) s(t,\omega) + \big( \Pi(t) F - Q \hat H\big)z^o(t,\omega) + \big( \Pi(t) G - Q  H\big)z^o_0(t,\omega) \\
& \hspace{2cm} - Q\eta \Big]dt - q(t,\omega)dw_0(t,\omega), \hspace{3.6cm}  s(T)=0,\\
& u^o(t,\omega) = - R^{-1} B^T \big( \Pi(t) z^o(t,\omega) + s(t,\omega)\big), \\
& dz^o(t,\omega) = \Big[ \big(\mathbb{A} (t) + F \big) z^o(t,\omega) - B R^{-1} B^T  \Pi(t) s(t,\omega) + G z^o_0 (t,\omega) \Big] dt
\\
& \hspace{2cm} + S dw(t,\omega), \hspace{6cm} z^o(0)=z(0).
\end{align*}

So we retrieve the MM-SMFG system for LQG dynamic games model of \cite{NguyenHuang_CDC11} for minor agents with uniform parameters (see equations (10)-(11) and (22)-(23) in \cite{NguyenHuang_CDC11}, see also \cite{HuangSIAM2010}). The reader is referred to \cite{NguyenHuang_CDC11} for an explicit representation of a solution to the SMFG LQG system under some appropriate conditions.

We note that key assumption for solution existence and uniqueness of MM-SMFG system is that all drift and cost functions and their derivatives are bounded (see Section \ref{Sub:ass}) which clearly does not hold for the MM-SMFG LQG problem (as in classical LQG control). In this case, a generalized Four-Step Scheme (see Section 5.2 in Chapter 7 of \cite{yong1999stochastic}) seems to give not only weaker general conditions but also presents explicit solutions to the MM-SMFG LQG case. This is currently under investigation and will be reported in future work. 

\section*{Appendix H: A Nonlinear Example}
\renewcommand{\theequation}{H.\arabic{equation}}
\setcounter{equation}{0}
\renewcommand{\thetheorem}{H.\arabic{theorem}}
\setcounter{theorem}{0}

In this section we present a major and minor version of the synchronization of coupled nonlinear oscillators game model \cite{yin2010synchronization}. Consider a population of $N+1$ oscillators with dynamics
\begin{align}
& d\theta_j^N(t)= u_j^N (t) dt + \sigma dw_j(t) \quad (\textrm{mod}~2\pi) & 0 \leq j \leq N, \quad t \geq 0, \label{osillatordyn}
\end{align}
where $\theta_j(t) \in [0,2\pi]$ is the phase of the $j^{\textrm{th}}$ oscillator at time $t$, $u_{j}(\cdot)$ is the control input, $\sigma$ is a non-negative scalar, and $\{w_{j}:0 \leq  j \leq N\}$ denotes a sequence of independent standard scalar Wiener processes (see \cite{yin2010synchronization}). It is assumed that the initial states $\{\theta_j(0)\}$ are chosen independently on $[0,2\pi]$. The objective of the $j^{\textrm{th}}$ oscillator is to minimize its own cost function
\begin{align}
& J_{0}^{N}(u_0^N,u_{-0}^N) := \E \int_0^{T} \Big(\frac{1}{N} \sum_{k=1}^N \sin^2\big[\theta_0^N(t)-\theta_k^N(t)\big] + r \big(u_0^N(t)\big)^2\Big)dt, \label{osillatorcos:major} \\
& J_{i}^{N}(u_i^N,u_{-i}^N) := \E \int_0^{T} \Big(\frac{1}{N} \sum_{k=1}^N \sin^2\Big[\theta_i^N(t)-\big(\lambda\theta_0^N(t) +(1-\lambda) \theta_k^N(t)\big)\Big]  \label{osillatorcos:minor} \\
& \hspace{4cm} + r \big(u_i^N(t)\big)^2\Big)dt,  \hspace{2cm} 1 \leq i \leq N, \notag
\end{align}
where $r$ is a positive scalar and $\lambda \in (0,1)$.

Similar arguments in previous section yield the following major agent's SMFG system (\ref{CNP:SHJB:Major})-(\ref{CNP:SMV:Major}):
\begin{align*}
& -d\phi_0(t,\omega,x)=\Big[-\frac{1}{4r} \big(\partial_x\phi_0(t,\omega,x))^2 + m_0(t,\omega,x) + \sigma \partial_x \psi_0(t,\omega,x) \\
& \hspace{2cm} +\frac{\sigma^2}{2} \partial_{xx}^2 \phi_0(t,\omega,x)\Big]dt -\psi_0(t,\omega,x)dw_0(t,\omega), \hspace{0.6cm} \phi_0(T,x) =0, \\
&  u^o_0(t,\omega,x) = -\frac{1}{2r} \partial_x\phi_0(t,\omega,x), \\ 
& dp^0_s(t,\omega,x) =  \Big[\frac{1}{2r} \partial_x \Big(\big(\partial_x\phi_0(t,\omega,x)\big) p_s^0(t,\omega,x) \Big) + \frac{\sigma^2}{2} \partial_{xx}^2 p_s^0(t,\omega,x) \Big] dt
\\
&  \hspace{2cm} -\sigma \partial_x  p_s^0(t,\omega,x)  dw_0(t,\omega), \hspace{3.5cm} p_s^0(s,x) = \delta_{\theta_0^o(s)}(dx), \\
& m_0(t,\omega,x) = \int_0^{2\pi} \sin^2(x-\theta) p(t,\omega,\theta) d \theta,
\end{align*}
where $m_0(t,\omega,x)$ is called the infinite population cost-coupling of the major agent, and $\theta_0^o(\cdot)$ is the solution of the closed-loop equation 
\begin{align*}
& d\theta_0^o(t)= u_0^o (t,\theta_0^o(t)) dt + \sigma dw_0(t) \quad (\textrm{mod}~2\pi) & \quad t \geq 0.
\end{align*}

In a similar way, the minor agents' SMFG system (\ref{CNP:SFPK1:Minor}) and (\ref{CNP:SHJB:Minor})-(\ref{CNP:BR:Minor}) is given by
\begin{align*}
& \hspace{-0.1cm} -d\phi(t,\omega,x)=\Big[-\frac{1}{4r} \big(\partial_x\phi(t,\omega,x))^2 + m(t,\omega,x)+\frac{\sigma^2}{2} \partial_{xx}^2 \phi(t,\omega,x)\Big]dt\\
& \hspace{2cm}  -\psi (t,\omega,x)dw(t,\omega), \hspace{5cm} \phi(T,x) =0, \\
&  u^o(t,\omega,x) = -\frac{1}{2r} \partial_x\phi(t,\omega,x), \\ 
& dp(t,\omega,x) =  \Big[\frac{1}{2r} \partial_x \Big(\big(\partial_x\phi(t,\omega,x)\big) p(t,\omega,x) \Big) + \frac{\sigma^2}{2} \partial_{xx}^2 p(t,\omega,x) \Big] dt, \quad p (0,x) \\
& m(t,\omega,x) = \int_0^{2\pi} \int_0^{2\pi} \sin^2\big(x-(\lambda \theta_0+(1-\lambda)\theta)\big) p_0^0(t,\omega,\theta_0) p(t,\omega,\theta)  d \theta_0 d \theta,
\end{align*}
where $m(t,\omega,x)$ is called the infinite population cost-coupling of the major agent. The reader is referred to the deterministic mean field system (14a)-(14c) in \cite{yin2010synchronization} for the synchronization of coupled nonlinear oscillators game model with only minor agents. 
\end{document}